\documentclass[12pt,thp]{report}

\usepackage{amssymb,amsfonts,amsmath,amsthm,amscd}
\usepackage[brazil]{babel}
\usepackage[latin1]{inputenc}

\newcommand{\ai} {\vbox to 7pt{\hbox to 7pt{\vrule height 7pt width 7pt}}}
\newcommand{\dps} {\displaystyle}

\newcommand{\N}{\mbox{$I\!\!N$}}
\newcommand{\Q}{\mbox{$I\!\!\!Q$}} 
\newcommand{\R}{\mbox{$I\!\!R$}}
\newcommand{\C}{\mbox{$|\!\!C$}}

\newcommand{\Z}{{\rm Z}\!\!{\rm Z}}

\newcommand{\SA} {\overline{S^A(\R^{n})}}
\newcommand{\sa} {S^A(\R^n)}
\newcommand{\CBA} {CB^{\infty}(\R^n,A)}
\newcommand{\CB} {CB^{\infty}(\R^n)}

\newcommand{\s} {S(\R^n)}

\newcommand{\op} {{\cal{B}}^{*}(\SA)}
\newcommand{\dc}{d\!\!\! /}
\newcommand{\supp}{{\rm sup}}
\newcommand{\epsi}{\varepsilon}

\newtheorem{predef}{{\bf Defini\c c\~ao}}[chapter]
\newenvironment{definition}{\setcounter{equation}{0}
\begin{predef}\rm}{\end{predef}}

\newtheorem{preexam}[predef]{{\bf Exemplo}}

\newtheorem{preprop}[predef]{{\bf Proposi\c c\~ao}}
\newenvironment{proposition}{\setcounter{equation}{0}
\begin{preprop}\rm}{\end{preprop}}

\newtheorem{prelema}[predef]{{\bf }}

\newtheorem{prelemma}[predef]{{\bf Lema}} 
\newenvironment{lemma}{\setcounter{equation}{0}
\begin{prelemma}\rm}{\end{prelemma}}

\newtheorem{pretheo}[predef]{{\bf Teorema}}
\newenvironment{theorem}{\setcounter{equation}{1}
\begin{pretheo}\rm}{\end{pretheo}}

\newtheorem{pretheom}[predef]{{\bf Teorema de Rohlin}}

\newtheorem{pretheoCD}[predef]{{\bf Teorema da Converg\^encia Dominada}}

\newtheorem{pretheoF}[predef]{{\bf Teorema de Fubini}}

\newtheorem{pretheoFC}[predef]{{\bf Teorema Fundamental do C\'alculo}}

\newtheorem{prenotac}[predef]{{\bf Nota\c c\~ao}}

\newtheorem{prenota}[predef]{{\bf Nota}}

\newtheorem{preobs}[predef]{{\bf Observa\c c\~ao}}

\newtheorem{precorol}[predef]{{\bf Corol\'ario}}

\newtheorem{precorollary}[predef]{{\bf Corol\'ario}}

\newtheorem{prenumremark}[predef]{{\bf remark}}

\newtheorem{preapendice}[predef]{{\bf AP\^ENDICE}}

\newenvironment{remark}{\par\noindent{\bf Demonstra\c c\~ao.\ }}{\par}

\newcommand{\citecorol}[1]{(\ref{#1})~}

\newcommand{\citelemma}[1]{(\ref{#1})~}
\newcommand{\citelema}[1]{(\ref{#1})~}
\newcommand{\citelem}[1]{(\ref{#1})~}
\newcommand{\citenota}[1]{(\ref{#1})~}
\newcommand{\citetheo}[1]{(\ref{#1})~}

\newcommand{\citeprop}[1]{(\ref{#1})~}
\newcommand{\citeobs}[1]{(\ref{#1})~}

\makeindex

\begin{document}

\vspace*{2cm}
\large

\thispagestyle{empty}

\centerline{\bf Resultados motivados por uma
caracteriza\c c\~ao}
\centerline{\bf de operadores
pseudo-diferenciais}
\centerline{\bf conjecturada por
Rieffel}


\vspace*{0.4cm}
\normalsize

\large
\vspace {1.1 cm}
\centerline{\bf\em Marcela Irene Merklen Olivera}

\normalsize
\vspace{2 cm}
\centerline{TESE APRESENTADA}
\centerline{AO}
\centerline{INSTITUTO DE MATEM\'ATICA E ESTAT\'ISTICA}
\centerline{DA}
\centerline{UNIVERSIDADE DE S\~AO PAULO}
\centerline{PARA}
\centerline{OBTEN\c C\~AO DO GRAU DE DOUTOR}
\centerline{EM}
\centerline{MATEM\'ATICA }

\large
\vspace {1.5 cm}
\centerline{Orientador: {\bf Prof.~Dr.~Severino Toscano do R\^ego Melo}}

\hspace{0.06cm}Co-orientadora: {\bf Prof.~Dra.~Cristina Cerri}

\vspace{1.3cm}

\normalsize

\centerline{\it Durante a elabora\c c\~ao deste trabalho a autora recebeu apoio financeiro do CNPq}


\normalsize

\pagebreak
\thispagestyle{empty}\vspace*{2cm}
\large

\centerline{\bf Resultados motivados por uma
caracteriza\c c\~ao}
\centerline{\bf de operadores
pseudo-diferenciais}
\centerline{\bf conjecturada por
Rieffel}

\vspace*{.5cm}
\normalsize
\normalsize
\vspace{2cm}
\begin{flushright}Este exemplar corresponde \`a reda\c c\~ao \\final da tese devidamente corregida e \\defendida por Marcela Irene Merklen Olivera \\e aprovada pela comiss\~ao julgadora.

\vspace{2cm}S\~ao Paulo,  outubro de 2002.\end{flushright}
\vspace{3.5cm}
Comiss\~ao Julgadora:
\begin{itemize}\item Prof.~Dr.~Severino Toscano do Rego Melo (Presidente) - IME-USP\item Prof.~Dr.~Ruy Exel Filho  - UFSC\item Prof.~Dr.~Milton Lopes - IMECC-UNICAMP\item Prof.~Dr.~Jos\'e Ruidival Santos Filho - UFSCar\item Profa.~Dra.~Beatriz Abadie - UNIV. REP. URUGUAI\end{itemize}
\newpage
\vspace*{16cm}
\pagestyle{empty}

\begin{flushright}
{\it a Alejandra, Ana, H\'ector,}\\
{\it minha fam\' \i lia,}\\
{\it com muito amor }
\end{flushright}
\vfill \break

\thispagestyle{empty}

\centerline{\bf Agradecimentos}

\vspace{1cm}
N\~ao h\'a palavras que possam expressar minha gratid\~ao por poder ter tido
a oportunidade de trabalhar com o professor Dr. Severino Toscano do Rego Melo,
meu orientador.

Gostaria de agradecer \`a professora Cristina Cerri, minha co-orientadora, que
tantas vezes me ajudou a clarear as id\'eias.

\`A professora Marina Pizzotti, por sua ajuda em Teoria da Integra\c c\~ao.

A tantos outros professores que, de um modo ou de otro, me ajudaram a concluir 
este trabalho.

Aos amigos t\~ao queridos  Pablo, Hern\'an, Liane, Daniel, Marcelo, Irene, Z\'e, 
Carolina, Marcel, Raul, Mait\'e, Aldo, B\'arbara, que me deram seu carinho 
e seu apoio todos estes anos.

A Cecilia Tosar, minha amiga, minha companheira, por estar sempre por perto e me ajudar
tanto, não só na digitação final da tese, mas também em várias ocasiões em minha vida.

A Fernando Abadie, um grande amigo, que me ajudou tantas vezes discutindo
matem\'atica comigo. E \`a sua esposa, Daniella, que sempre teve palavras
doces para mim.

A meus amigos no Uruguay, Gerardo, Yolanda, Mauro, Fernando que, a pesar
de morar longe, sempre estiveram pr\'oximos.

A meu pai, professor Dr. H\'ector Merklen, que me ajudou discutindo problemas
e a quem devo a digita\c c\~ao da tese.

A minha m\~ae, professora Dra. Ana Mar\'\i a Olivera Yam\'\i n, quem sempre me acompanhou e me ajudou.

A minha irm\~a, Alejandra, que sempre foi carinhosa comigo.

A minha fam\'\i lia, em fim, sem a qual n\~ao sou nada, n\~ao perten\c co a nada, que 
tanto me abrigou e me ajudou nestes anos e em toda minha vida.

Finalmente, a todos aqueles que n\~ao tive oportunidade de mencionar, mas 
com quem tive  oportunidade de compartilhar minha vida nestes anos, muito obrigada.

\newpage

\thispagestyle{empty}

\vspace{1,5cm}

\begin{center} {\bf Resumo}\end{center}

\medskip Trabalhamos com fun\c c\~oes definidas em $\R^n$ que tomam valores 
numa $C^*$-\'algebra $A$.
Consideramos o conjunto $\sa$ das fun\c c\~oes de Schwartz, (de decrescimento
r\'apido), com norma dada por $\|f\|_2=
\|\int f(x)^*f(x)dx\|^{\frac{1}{2}}$. Denotamos por $CB^\infty(\R^{2n},A)$ o conjunto
das fun\c c\~oes $C^\infty$ com todas as suas derivadas limitadas.
Provamos que os operadores pseudo-diferenciais 
com s\'\i mbolo em $CB^\infty(\R^{2n},A)$ s\~ao cont\'\i nuos em $\sa$
com a norma $\|\cdot\|_2$, fazendo uma generaliza\c c\~ao de \cite{10}.
Rieffel prova em \cite{1} que $CB^\infty(\R^n,A)$ age em $\sa$ por
meio de um produto deformado, induzido por uma matriz anti-sim\'etrica, $J$,
como segue: $L_Fg(x)=F\times_Jg(x) = \int e^{2\pi iuv}F(x+Ju)g(x+v)du dv$,
(integral oscilat\'oria).

Dizemos que um operador $S$ \'e Heisenberg-suave se as aplica\c c\~oes $z \mapsto
T_{-z}ST_z$ e $\zeta \mapsto M_{-\zeta}SM_\zeta, \,\, z,\zeta\in\R^n$, s\~ao
$C^\infty$; onde $T_zg(x)=g(x-z)$ e $M_\zeta g(x)=e^{i\zeta x}g(x)$. No final
do cap\'\i tulo 4 de \cite{1}, Rieffel prop\~oe uma conjectura: que todos os
operadores "adjunt\'aveis" \, em $\sa$, Heisenberg-suaves, que comutam com
a representa\c c\~ao regular \`a direita de $CB^\infty(\R^n,A)$, $R_Gf = f\times_JG$,
s\~ao os operadores do tipo $L_F$. Provamos este resultado para o caso
$A=\C$, \cite{14}, usando a caracteriza\c c\~ao de Cordes \cite{17} dos operadores Heisenberg-suaves em $L^2(\R^n)$ como sendo os operadores pseudo-diferenciais
com s\'\i mbolo em $CB^\infty(\R^{2n})$. Tamb\'em \'e provado neste trabalho que, 
se vale uma generaliza\c c\~ao natural da caracteriza\c c\~ao de Cordes, 
a conjectura de Rieffel \'e verdadeira.

\vspace{1,5cm}
\centerline{\bf Abstract}

\vspace{1cm}
We work with functions defined on $\R^n$ with values in a $C^*$-algebra $A$.
We consider the set $\sa$ of Schwartz functions (rapidly decreasing), with norm
given by $\|f\|_2 = \|\int f(x)^*f(x)dx\|^{\frac{1}{2}}$. We denote $CB^\infty(\R^{2n},A)$
the set of functions which are $C^\infty$ and have all their derivatives bounded.
We prove that pseudo-differential operators with symbol in $CB^\infty(\R^{2n},A)$
are continuous on $\sa$ with the norm $\|\cdot\|_2$, thus generalizing the result 
in \cite{10}. Rieffel proves in \cite{1} that $CB^\infty(\R^n,A)$ acts on $\sa$
through a deformed product induced by an anti-symmetric matrix, J, as
follows: $L_Fg(x) = F\times_Jg(x) = \int e^{2\pi iuv}F(x+Ju)g(x+v)dudv$ 
(an oscillatory integral).

We say that an operator $S$ is Heisenberg-smooth if the maps $z \mapsto
T_{-z}ST_z$ and $\zeta \mapsto M_{-\zeta}SM_\zeta, \,\, z,\zeta \in \R^n$ are
$C^\infty$; where $T_zg(x) = g(x-z)$ and where $M_\zeta g(x)=e^{i\zeta x}g(x)$. At the end of chapter 4 of \cite{1}, Rieffel proposes a conjecture: that all
"adjointable" operators in $\sa$ that are Heisenberg-smooth and that
commute with the right-regular representation of $CB^\infty(\R^n,A)$, 
$R_Gf=f\times_JG$, are operators of type $L_F$. We proved this result
for the case $A=\C$ in \cite{14}, using Cordes´ characterization of 
Heisenberg-smooth operators on $L^2(\R^n)$ as being the pseudo-differential
operators with symbol in $CB^\infty(\R^{2n})$. It is also proved in this
thesis that, if a natural generalization of Cordes´ characterization is valid,
then the Rieffel conjecture is true.

\pagebreak

\large
\thispagestyle{empty}

\centerline {\bf Sum\'ario}

\normalsize

\baselineskip=7mm

\vskip 2cm

\par\noindent Introdu\c c\~ao \dotfill 1

\vskip 0.5cm

\par\noindent Cap\'\i tulo 1 \dotfill 4

\vskip 0.5cm

\par\noindent Cap\'\i tulo 2 \dotfill 19

\vskip 0.5cm

\par\noindent Cap\'\i tulo 3 \dotfill 26

\vskip 0.5cm

\par\noindent Cap\'\i tulo 4. \dotfill 36

\vskip 0.5cm

\par\noindent Ap\^endice A \dotfill 49

\vskip 0.5cm

\par\noindent Ap\^endice B \dotfill 59

\vskip 0.6cm

\par\noindent Refer\^encias Bibliograficas \dotfill 61

\newpage

\vskip 4.5cm


\vfill\break

\setcounter{page}{1}

\begin{center} {\bf Introdu\c c\~ao}\end{center}


\vspace{1cm}
Consideramos uma $C^*$-\'algebra, $A$, e fun\c c\~oes de $\R^n$ em $A$. Em particular,
trabalhamos com o conjunto das fun\c c\~oes de Schwartz, $\sa$, que consiste nas fun\c c\~oes
$f:\R^n\rightarrow A$, $C^\infty$, tais que $\sup_{x\in\R^n}\|x^\alpha\partial^\beta f(x)\|$ \'e limitado para quaisquer multi-\'\i ndices $\alpha$ e $\beta$, \cite[cap\'\i tulo 3]{1}. Neste conjunto
definimos o "produto interno" 
 $<f,g> = \int f(x)^*g(x)dx$, \cite[(4.1)]{1}. $\sa$, com a norma
induzida por este produto interno, $\|f\|_2=\|<f,f>\|^\frac{1}{2}$, \'e um pr\'e-m\'odulo de Hilbert
(ver defini\c c\~ao \ref{1umquarto}).
Denotamos por $\SA$ o m\'odulo de Hilbert correspondente, isto \'e o completamento
de $\sa$ com respeito a esta norma.

\medskip Observamos que o espa\c co de Fr\'echet $CB^\infty(\R^n,A)$, das fun\c c\~oes 
$C^\infty$ limitadas, com todas as derivadas limitadas, age em $\sa$ do seguinte modo:
dadas $F\in CB^\infty(\R^n,A)$ e $g\in\sa$, $L_F(g)=F\times_J g\in\sa,$ como
vemos em \cite[(3.3)]{1}. Aqui Rieffel define o produto deformado, $\times_J$, para
uma matriz, $J, n\times n$, anti-sim\'etrica, em $CB^\infty(\R^n,A)$, 
como segue: $F\times_J G(x) =
\int\int F(x+Ju)G(x+v)e(u\cdot v)dudv$, onde $e(u\cdot v) = exp(2\pi iu\cdot v)$
e a integral acima \'e uma integral oscilat\'oria, \cite[(1.3)]{1}. O operador $L_F$ definido
acima possui "adjunto", \cite[(4.2)]{1}, (com respeito ao "produto interno"  acima definido),
e \'e cont\'\i nuo no m\'odulo de Hilbert $\SA$, \cite[(4.7)]{1}.

\medskip Seja $B^*(\SA)$ o conjunto dos operadores lineares cont\'\i nuos em $\SA$
que possuem adjunto. Como no caso escalar, \cite[cap\'\i tulo 8]{2}, consideramos
em $B^*(\SA)$ a a\c c\~ao do grupo de Heisenberg como segue:
$ (z,\zeta,\varphi) \mapsto E_{z,\zeta,\varphi}\Lambda E^{-1}_{z,\zeta,\varphi},\hspace{0.2cm}
z,\zeta\in\R^n,\varphi\in\R,$ para $\Lambda\in B^*(\SA)$, onde $E_{z,\zeta,\varphi}=e^{i\varphi}M_\zeta T_z$, com $T_z u(x)=u(x-z)$, 
$M_\zeta u(x)=e^{\zeta x}u(x), \hspace{0,3cm} u\in\sa$.
Notamos que $E_{z,\zeta\varphi}^{-1}\Lambda E_{z,\zeta,\varphi}$ independe de $\varphi$.
Tomamos, ent\~ao, $\varphi=0$ e usamos a nota\c c\~ao $E_{z,\zeta}=E_{z,\zeta,0}$ \, e
\,$\Lambda_{z,\zeta}=E_{z,\zeta}^{-1}\Lambda E_{z,\zeta}$. Dizemos que o operador
$\Lambda$ \'e Heisenberg-suave se a aplica\c c\~ao $(z,\zeta) \mapsto \Lambda_{z,\zeta}
\,\, ((z,\zeta)\in\R^{2n})$ \'e $C^\infty$.

\medskip Assim como definimos a a\c c\~ao \`a esquerda de $CB^\infty(\R^n,A)$ em $\sa$
por $L_F(g)=F\times_Jg$, podemos definir a a\c c\~ao \`a direita de $CB^\infty(\R^n,A)$ em
$\sa$ por $R_G(f)=f\times_JG, \,\, G\in CB^\infty(\R^n,A)$ e $f\in\sa$. $R_G$ \'e um
operador cont\'\i nuo no espa\c co de Banach $\SA$, embora n\~ao seja homomorfismo
de m\'odulos \`a direita.


\bigskip
No final do cap\'\i tulo 4 de \cite{1}, Rieffel prop\~oe uma quest\~ao que podemos enunciar
do seguinte modo: o conjunto dos operadores em $B^*(\SA)$, Heisenberg-suaves,
que comutam com $R_G$, para toda $G\in CB^\infty(\R^n,A)$, \'e o conjunto dos
operadores $L_F, \,\, F\in CB^\infty(\R^n,A)$.

\medskip No primeiro cap\'\i tulo deste trabalho demonstramos esta conjectura
para o caso em que $A$ \'e a $C^*$-\'algebra dos n\'umeros complexos, \cite{14}.
Para isto, usamos a caracteriza\c c\~ao feita por Cordes em \cite{17} dos operadores lineares
cont\'\i nuos em $L^2(\R^n)$ (no caso em que $A=\C$, $\SA=L^2(\R^n)$)
que s\~ao Heisenberg-suaves. Cordes prova que estes s\~ao os operadores $a(x,D)$,
dados por um s\'\i mbolo, uma fun\c c\~ao $a\in CB^\infty(\R^{2n})$, do seguinte modo:
$a(x,D)u(x) = \int\int(2\pi)^{-n}e^{i(x-y)\xi}a(x,\xi)u(y)dyd\xi$, onde $u\in C^\infty_0(\R^n)$.
Provamos, ent\~ao, que, se $a(x,D)$ comuta com $R_G$ para toda $G$ em 
$CB^\infty(\R^{n})$, ent\~ao se $F(x)=a(x,0)$, $a(x,D)=L_F$.

\medskip No cap\'\i tulo 2, precisamos supor que $A$ \'e separ\'avel para provar um
resultado interessante.  Na verdade, se $A$ tiver unidade, a separabilidade
n\~ao \'e necess\'aria, pois s\'o a necessitamos para ter uma aproxima\c c\~ao da
unidade em $A$ dada por uma seq\"u\^encia. 
Dizemos que um operador $\Lambda$ em $B^*(\SA)$ \'e
suave por transla\c c\~oes se a aplica\c c\~ao $z \mapsto \Lambda_z, \,\, z\in\R^n$
\'e $C^\infty$, onde temos que $\Lambda_z = T^{-1}_z\Lambda T_z$, usando
a nota\c c\~ao anterior. Neste cap\'\i tulo, provamos que se $\Lambda$ \'e suave por 
transla\c c\~oes e comuta com $R_g$ para qualquer $g$ em $\sa$, ent\~ao
$\Lambda$ \'e Heisenberg-suave.

\medskip O cap\'\i tulo 3 consiste em uma generaliza\c c\~ao do teorema de Calder\'on-Vaillancourt, \cite{10}, baseada no teorema 3.14 de \cite{11}. Provamos que um
operador no m\'odulo $\SA$ dado por um s\'\i mbolo (como no caso escalar,
\cite{2}) \'e cont\'\i nuo. Aqui estamos nos referindo aos operadores $a(x,D), \,\, a\in CB^\infty(\R^{2n},A)$ dados por $a(x,D)u(x) = \int\int (2\pi)^{-n}
e^{i(x-y)\xi}a(x,\xi)$ $u(y)dyd\xi,$ para $u\in\sa$. A integral acima \'e uma integral oscilat\'oria. Novamente supomos que $A$ \'e separ\'avel,  para poder aplicar
o Teorema de Converg\^encia Dominada, e outras propriedades, para a integral de Bochner, \cite[(E.6)]{13}, e tamb\'em para usar que a transformada de Fourier \'e um operador
"unit\'ario" em $\sa$ (ver ap\^endice B) e que $\sa$ \'e denso em $L^2(\R^n,A)$
(ver ap\^endice A).

\medskip Usando os resultados dos cap\'\i tulos 2 e 3, no cap\'\i tulo 4, provamos que uma gene~-raliza\c c\~ao natural, 
para uma $C^*$-\'algebra separ\'avel qualquer, $A$, 
da caracteriza\c c\~ao de Cordes 
implicaria que a conjectura de Rieffel \'e verdadeira. 
Isto generaliza os resultados de \cite{14}; l\'a n\'os provamos
que a caracteriza\c c\~ao de Cordes implica a conjectura de Rieffel, no caso $A=\C$,
atrav\'es de um argumento que usa o n\'ucleo de Schwartz de um operador 
pseudo-diferencial. No cap\'\i tulo 4, os argumentos que envolvem distribui\c c\~oes 
temperadas no cap\'\i tulo 1, ser\~ao substitu\'\i dos pela proposi\c c\~ao 4.11, que se trata, na verdade,
de uma extens\~ao imediata, para o caso $A$ qualquer, 
de um trecho da demonstra\c c\~ao na caracteriza\c c\~ao de Cordes.

\chapter{ }

\vspace{1cm}

Neste cap\'\i tulo vamos apresentar a conjectura proposta por Rieffel no final de \cite[Cap\'\i tulo 4]{1} para uma $C^*$-\'algebra qualquer e veremos uma resposta para o caso da $C^*$-\'algebra dos n\'umeros complexos.

\bigskip
\bigskip

\leftskip=0cm \rightskip=0cm

Seja $A$ uma $C^{*}$-\'algebra.  Consideremos o espa\c co das fun\c c\~oes cont\'\i nuas e limitadas de $\R^{n}$ em $A$, que possuem derivadas de todas as ordens que s\~ao cont\'\i nuas e limitadas.

$$\CBA = \{ F: \R^{n} \longrightarrow A ; \,\, \dps \| \partial^{\alpha} F \dps \|_ { \infty} \leq \infty, \forall \alpha \}, $$

\noindent onde $\alpha$ \'e um multi-\'\i ndice, $\alpha = (\alpha_1,\alpha_2,\cdots,\alpha_n) \in \N^n$ e

$$\partial^{\alpha}F = \partial^{\alpha_1}\partial^{\alpha_2}\cdots \partial^{\alpha_n}F, $$

$$\dps \| \partial^{\alpha}F \dps \|_{\infty} = \sup_{x\in\R^n} \dps \|\partial^{\alpha}F(x) \dps \|, $$ onde esta \'ultima \'e a norma em $A$.

\medskip Consideremos em $\CBA$ as seminormas $$ \dps \|F\|_j = \sup_{\dps | \alpha|\leq j} \dps \|\partial^\alpha F \dps \|_ \infty, $$ \noindent onde $|\alpha| = \alpha_1+\alpha_2+\cdots+\alpha_n$. 

Deste modo vemos que $\CBA$ \'e um espa\c co de Fr\'echet.

\vspace{1cm}

\medskip Dada $J$, uma matriz anti-sim\'etrica $n\times n$, definimos o {\em produto deformado} $\times_J$ em $\CBA$ como segue:

\[ F \times_J G (x) = \int F(x+Ju)G(x+v)e^{iu\cdot v}\dc u\dc v, \] onde $\dc u = 
(2\pi)^{-\frac{n}{2}}du, \,\,\, \dc v = (2\pi)^{-\frac{n}{2}}dv.$ Notemos que a rigor, considerando
a defini\c c\~ao de Rieffel do cap\'\i tulo 3 de \cite{1}, nosso $\times_J$ \'e seu 
$\times_{2\pi J}$.

A integral acima \'e uma integral oscilat\'oria. Isto \'e, se temos $\psi_m\in C^\infty_c(\R^n)$, que satisfaz 

\medskip
$(i)$ \hspace{0,7cm} $ \supp_{x\in\R^n}|\partial^\alpha\psi_m(x)| < c_\alpha$, para cada multi-\'\i ndice $\alpha$, para todo
$m \in \N$, onde $c_\alpha \in \R^+$ \'e uma constante que depende s\'o de $\alpha$.

\medskip
$(ii)$ \hspace{0,7cm} $\psi_m(x) = 1$, se $x$ pertence \`a bola de centro na origem e raio $r_m$, $B(O,r_m)$,
onde $r_m$ cresce para infinito;

\medskip
e tamb\'em $\psi'_k \in C^\infty_c(\R^n)$, que satisfaz $(i)$ e $(ii)$, ent\~ao
\[ \dps \int F(x+Ju)G(x+v)e^{iu\cdot v}\dc u \dc v = \lim_{k,m \rightarrow \infty} \dps\int\dps\int F(x+Ju)G(x+v)\psi_m(u)\psi'_k(v)e^{iu\cdot v}\dc u \dc v .\]
Prova-se que o limite existe e \'e independente da escolha de $\psi_m$ e $\psi'_k$ (\cite[(2.4)]{1}).

\bigskip
O produto deformado est\'a bem definido em $\CBA$ (\cite[2.2]{1}).

\begin{preprop}\label{asterisco}{\rm   Dadas $F, G \in CB^\infty(\R^n,A)$, temos
que 
\[ F \times_J G(x) = \dps\int e^{iu(x-v)}F(x-Ju)G(v)\dc v\dc u.\]                                   }\end{preprop}

\begin{remark} Pelo visto acima, temos que
\[ F \times_J G(x) = \dps\lim_{m,k}\dps\int\dps\int e^{iuw}F(x+Ju)G(x+w)\psi_m(u)
\psi'_k(w)\dc u\dc w = \]
\[ = \dps\lim_{m,k\rightarrow\infty}\dps\int\dps\int e^{iu(v-x)}F(x+Ju)G(v)\psi_m(u)
\psi'_k(v-x)\dc u\dc v = \] \[ = \dps\lim_{m,k\rightarrow\infty}\dps\int\dps\int e^{-iu(v-x)}
F(x-Ju)G(v)\psi_m(-u)\psi'_k(v-x)\dc v\dc u\] pela proposi\c c\~ao \ref{A.16}  e
pelo Teorema de Fubini. Para $x$ fixo, $\psi_m(-u)$ e $\psi'_k(v-x)$ satisfazem as 
condi\c c\~oes   (i) e (ii) acima, ent\~ao temos:
\[ F \times_J G(x) = \dps\int e^{iu(x-v)}F(x-Ju)G(v)\dc u\dc v. \hspace{0,5cm} \ai\]\end{remark}
\vspace{1cm}
\medskip
Interessa-nos $\CBA$ agindo num espa\c co com produto interno. Assim, vamos considerar as fun\c c\~oes de Schwartz.  Seja $S^{A}(\R^n)$ o subespa\c co de $\CBA$ das fun\c c\~oes tais que o produto de suas derivadas por polin\^{o}mios em $\R^n$ s\~ao limitadas. 

\[ \sa  = \{ f:\R^n \rightarrow A / \,\, \dps \|x^{\alpha}\partial^{\beta}f \dps \|_{\infty} < \infty \} \] para $\alpha$ e $\beta$ multi-\'\i ndices.

\medskip
Dadas $F \in \CBA$ e $g \in \sa$, prova-se que $F\times_J g$ e $g \times_J F$  est\~ao em $\sa$ (\cite[3.3]{1}).

\begin{predef}\label{1umquarto} {\rm Um {\em pr\'e-m\'odulo de Hilbert} sobre $A$
\'e um espa\c co vetorial sobre $\C$, $E$, que \'e tamb\'em um m\'odulo  \`a direita
equipado com uma aplica\c c\~ao $<\cdot,\cdot>:E \times E \rightarrow A$ que \'e linear
na segunda vari\'avel e satisfaz, para $a \in A, x, y \in E$, as seguintes propriedades:

\medskip
\noindent (i) $<x,ya> = <x,y>a $

\noindent (ii) $<x,y>^* = <y,x>$

\noindent (iii) $<x,x> \, \geq \, 0$

\noindent (iv) se $<x,x> = 0$, ent\~ao $x = 0. $ 

\bigskip
Em $E$ definimos $\|x\| = \|<x,x>\|^{\frac{1}{2}}, \hspace{0,5cm} x \in E$. }
\end{predef}

Prova-se que assim $E$ \'e um espa\c co normado (\cite[lema 1.1.2]{7}). Se $E$ \'e 
completo com respeito a esta norma, dizemos que $E$ \'e um {\em m\'odulo de Hilbert.}

\begin{prenota}\label{1ummedio} {\rm Podemos ler mais sobre m\'odulos de
Hilbert em \cite[cap\'\i tulo 1]{7}. } \end{prenota}

Consideremos em $\sa$ o produto interno \[ <f,g> = \int f(x)^*g(x)dx.\] \'E f\'acil ver que $\sa$ com este produto interno \'e um pr\'e-m\'odulo de Hilbert sobre $A$. A norma correspondente \'e 

\[ \dps \|f\|_2 = \dps \| <f,f> \dps \|^{\frac{1}{2}},  \hspace{1cm} f \in \sa . \]

\medskip Denotemos por $L$ a a\c c\~ao \`a esquerda de $\CBA$ em $\sa$ dada pelo produto deformado como segue:

\[ L_Fg = F \times_J g  \hspace{1cm}  F\in \CBA, \,\, g \in \sa. \]

\medskip
Vemos em \cite[(4.7)]{1}  que $L_F$ \'e um operador cont\'\i nuo em $\sa$ como
pr\'e-m\'odulo de Hilbert. No cap\'\i tulo 3 deste trabalho, veremos uma 
demonstra\c c\~ao do Teorema de Calder\'on-Vaillancourt \cite{10} para operadores
pseudo-diferenciais com s\'\i mbolo tomando valores em $A$, 
que tamb\'em implica
que $L_F$ \'e cont\'\i nuo.

\medskip
Al\'em disso, 
o operador $L_F$ \'e "adjunt\'avel", isto \'e, existe $(L_F)^*$, operador limitado em $\sa$,
tal que $<L_Fu,v>=<u,(L_F)^*v>$, para todo $u, v \in \sa$. De fato,
denotando por $F^*$ a fun\c c\~ao $F^*(x)=F(x)^*$, demonstra-se (ver \cite[(4.2)]{1})
que $(L_F)* = L_{F^*}$.

\medskip Do mesmo modo, denotemos por $R$ a a\c c\~ao \`a direita de $\CBA$ em $\sa$:

\[ R_Gf = f \times_J G   \hspace{1cm} f \in \sa , \,\, G \in \CBA .\]

Observemos que $R_G$ n\~ao \'e um homomorfismo no m\'odulo de Hilbert $\sa$, pois,
dado $a\in A$, 
\[R_G(fa)(x) = \dps\int f(x+Ju)aG(x+v)e^{i u\cdot v}\dc u\dc v\]
\[R_G(f)a(x) = \dps\int f(x+Ju)G(u+v)ae^{i u\cdot v}\dc u\dc v.\] Ou seja, se $A$ n\~ao
\'e comutativa, n\~ao vale necess\'ariamente que $R_G(f)a = R_G(fa)$.

Como feito em \cite[(4.4)]{1}, temos que, dada $g\in\sa$, o operador $R_g$ \'e cont\'\i nuo no
espa\c co de Banach $\SA$, o completamento de $\sa$ com a norma $\dps \| \cdot \dps \|_2$.

\vspace{1cm}
\medskip Denotemos com $\op$ os operadores em $\SA$ que possuem adjunto. 
Como feito por Cordes, para o caso $A = \C$ (\cite[Cap.8,Sec.2]{2}), consideremos a a\c c\~ao do grupo de Heisenberg em $\op$.
O Grupo de Heisenberg \'e o espa\c co $\R^{2n+1}$ equipado com a opera\c c\~ao

\[ (z,\zeta,\varphi)\circ (z',\zeta',\varphi') = (z+z',\zeta +\zeta', \varphi +\varphi'- \zeta'z), \] unidade =$(0,0,0)$.

\noindent Para $z,\zeta  \in \R^n$, sejam os operadores $T_z,M_\zeta \in \op$ dados por

\[ T_zf(x) = f(x-z) \]
\[ M_\zeta f(x)= e^{i\zeta x}f(x) \]  para $f \in \sa$. Notemos que $T_z$ e $M_\zeta$
s\~ao "unit\'arios"  \,\,com respeito ao "produto interno"  \,\,$<\cdot,\cdot>$.

Observemos que vale

\medskip \begin{center}
$ 
 T_{-z}M_\zeta T_zM_{-\zeta} =  e^{iz\zeta }  \hspace{0,5cm}
 M^*_\zeta = M_{-\zeta}  \hspace{0,5cm}
T^*_z = T_{-z}  $  \end{center}

\medskip
Assim, o sub-grupo dos operadores unit\'arios em $\op$ gerado pelos operadores $T_z$ e
$M_\zeta$, $GH$, consiste no conjunto dos operadores

\begin{center} $ E_{z,\zeta,\varphi} = e^{i\varphi}M_\zeta T_z ,  \hspace{0,7cm} \varphi \in \R, \,\, z,\zeta
\in \R^n .$ \end{center}

\medskip Deste modo, temos uma aplica\c c\~ao do grupo de Heisenberg em $GH$, bem definida
(n\~ao injetora), que a cada $(z,\zeta,\varphi) \in \R^{2n+1}$ associa o operador $E_{z,\zeta,\varphi}$.
A partir de agora, estaremos pensando no grupo $GH$ quando mencionarmos o grupo de Heisenberg.

\bigskip
Consideremos ent\~ao a a\c c\~ao do grupo de Heisenberg em $\op$ por conjuga\c c\~ao.
Para $T \in \op $ seja $T_{z,\zeta} = E^{-1}_{z,\zeta,\varphi}TE_{z,\zeta,\varphi}$.

Observemos que $T_{z,\zeta}$ independe de $\varphi$. De agora em diante tomemos $\varphi=0$ e denotemos por $E_{z,\zeta}$ o operador $E_{z,\zeta,0}$.
No final de \cite[Cap.4]{1}, Rieffel faz uma conjectura que, no presente contexto, pode ser
enunciada como segue:

\medskip \leftskip=1cm \rightskip=1cm  {\em O conjunto dos operadores $L_F$ com $F$ em $\CBA$ \'e o conjunto dos operadores adjunt\'aveis em $\SA$ que s\~ao suaves pela a\c c\~ao do grupo de Heisenberg e comutam com $R_G$ para toda $G$ em $\CBA$, ou seja os operadores $T \in \op$ tais que a aplica\c c\~ao $(z,\zeta) \rightarrow T_{z,\zeta}$ \'e $C^\infty$ para $(z,\zeta) \in \R^{2n}$ e $[T,R_G]=0, \,\,\, \forall G \in \CBA$}.

\leftskip=0cm \rightskip=0cm

\begin{predef}\label{1ummeio} {\rm Dado $T \in \op$, dizemos que $T$ \'e {\em Heisenberg-suave} se a aplica\c c\~ao $(z,\zeta) \rightarrow T_{z,\zeta}$ \'e $C^\infty$, para $(z,\zeta) \in \R^{2n}$.} \end{predef}

\begin{preprop}\label{1tresquartos} {\rm Dada $F \in \CBA$, o operador $L_F$ \'e 
Heisenberg-suave.     }\end{preprop}
\begin{remark} Dada $g \in \sa$, temos, pela proposi\c c\~ao \ref{asterisco}
\[ L_Fg(x) = \dps\int e^{iu(x-v)}F(x-Ju)g(v)\dc v \dc u .\]

Para $z,\zeta \in \R^n$, por argumentos similares aos usados na proposi\c c\~ao \ref{asterisco},

\[ E_{z,\zeta}^{-1}L_F E_{z,\zeta}g(x) = T_{-z}M_{-\zeta}L_F E_{z,\zeta}g(x) 
 =  e^{-i\zeta(x+z)}L_F E_{z,\zeta}g(x+z) = \] 
\[= e^{-i\zeta(x+z)}L_F E_{x,\zeta}g(x+z) = e^{-i\zeta(x+z)}\dps\int e^{iu(x+z-v)}
F(x+z-Ju)E_{z,\zeta}g(v)\dc v\dc u = \]
\[ = e^{-i\zeta(x+z)}\dps\int e^{iu(x+z-v)}F(x+z-Ju)e^{i\zeta v}g(v-z)\dc v\dc u =\]
\[ = e^{-i\zeta(x+z)}\dps\int e^{iu(x-w)}F(x+z-Ju)e^{i\zeta(w+z)}g(w)\dc w\dc u = \]
\[ = e^{-i\zeta x}\dps\int e^{iux}e^{-iw(u-\zeta)}F(x+z-Ju)g(w)\dc w\dc u =\]
\[ = e^{-i\zeta x}\dps\int e^{i(\xi+\zeta)x}e^{-iw\xi}F(x+z-J(\xi+\zeta))g(w)\dc w\dc \xi =\]
\[ = \dps\int e^{i\xi(x-w)}F(x+z-J(\xi+\zeta))g(w)\dc w\dc \xi.\]

(Notemos que esta manipula\c c\~ao formal \'e 
v\'alida tamb\'em para integrais oscilat\'orias.)

\bigskip Consideremos o quociente:

\[ \frac{(L_F)_{x+te_j,\zeta} - (L_F)_{z,\zeta}}{t} =  
\frac{E^{-1}_{z,\zeta}(L_F)_{te_j,0}E_{z,\zeta} -  E^{-1}_{z,\zeta}L_FE_{z,\zeta}}{t} = 
 E^{-1}_{z,\zeta}\left[\frac{(L_F)_{te_j,0} - L_F}{t}\right]E_{z,\zeta},\] onde $e_j=(0,\cdots,1,0,\cdots,0) \in \R^n $.

E, pelo visto acima,

\[ \left[\frac{(L_F)_{te_j,0} - L_F}{t}\right](g)(x) = \] \[ = \dps\int e^{i\xi(x-w)}
\left[\frac{F(x+te_j - J\xi) - F(x-J\xi)}{t}\right]g(w)\dc w \dc\xi .\]

Por \cite[(4.7)]{1} temos que 
existem $k \in \N$ e uma constante $c_k \in \R^+$,
tais que $\|L_G\| \leq c_k\|G\|_{2k}$, para toda $G \in\CBA$. 
Sejam $H_t(y)= \frac{F(y  + te_j) -F(y)}{t}$
e $\partial_{e_j}F(y) = \lim_{t \rightarrow 0} H_t(y)$. Como $F \in CB^\infty(\R^n,A)$, temos
que para todo $l \in \N$, $\lim_{t \rightarrow 0}\|H_t - \partial_{e_j}F\|_l = 0.$
De fato, dado $l\in\N$, $\|H_t-\partial_{e_j}F\|_l=\dps\sup_{|\alpha|\leq l}
\|\partial^\alpha(H_t-\partial_{e_j}F)\|_\infty.$ Pelo Teorema Fundamental do C\'alculo,
\ref{A.13}, 
\[H_t(y)-\partial_{e_j}F(y) = \frac{F(y+te_j)-F(y)}{t} - \partial_{e_j}F(y) = \]
\[ = \dps\int_{[0,t]}\left[\frac{\partial_{e_j}F(y+he_j)}{t}-\partial_{e_j}F(y)\right]dh = 
\dps\int_{[0,t]}\frac{1}{t}[\partial_{e_j}F(y+he_j)-\partial_{e_j}F(y)]dh =\]
\[ = \dps\int_{[0,t]}\frac{1}{t}\dps\int_{[0,h]}\partial^2_{e_j}F(y+se_j)dsdh.\] Logo,
$\|H_t(y) - \partial_{e_j}F(y)\|\leq t\|\partial^2_{e_j}F\|_\infty, \hspace{0.5cm} \forall y\in\R^n.$

\medskip De modo an\'alogo, temos que \[\|\partial^\alpha(H_t-\partial_{e_j}F)\|_\infty \leq 
t\|\partial^2_{e_j}\partial^\alpha F\|_\infty, \hspace{0,5cm} |\alpha|\leq l,\] e, assim,
$\dps\lim_{t\rightarrow 0}\|H_t-\partial_{e_j}F\|_l = 0, \hspace{0,5cm} \forall l\in\N$.

\medskip
Deste modo, $\|L_{H_t} - L_{\partial_{e_j}F}\| \leq c_k\|H_t -  \partial_{e_j}F\|_{2k} $ e,
portanto, $\lim_{t \rightarrow 0} \|L_{H_t} - L_{\partial_{e_j}F}\| = 0$,  o
que implica que existe a derivada parcial
$\partial_{z_j}(L_F)_{z,\zeta}=L_{\partial_{e_j}F}$. Como obtivemos um operador
do mesmo tipo, o processo pode ser repetido indefinidamente. Logo, $L_F$ \'e Heisenberg-suave. $\ai$\end{remark}

\begin{preprop}\label{1.1} {\rm Dadas $F, G \in \CBA$, quaisquer, 
$[L_F,R_G] = 0$.} \end{preprop}

\begin{remark}  Consideremos $u \in \sa$. \[L_FR_Gu = F \times_J R_Gu = F \times_J(u \times_JG)\]
\[R_GL_Fu = L_Fu \times_JG = (F \times_J u)\times_JG .\]  Como o producto deformado \'e associativo (\cite[2.14]{1}), temos $[L_F,R_G] = 0$. $\ai$ \end{remark} 

\begin{prelema}\label{1ummeio} {\rm Deste modo,  dado $F \in \CBA$, temos que $L_F$ \'e um operador Heisenberg-suave que comuta com $R_G$ para todo $G \in \CBA$. }\end{prelema}

\begin{preobs}\label{1umcuarto} {\rm Pelo visto na proposi\c c\~ao \ref{asterisco},
$L_F g(x) = \dps\int e^{iu(x-v)}F(x-Ju)g(v)\dc v\dc u$}\end{preobs}

\bigskip
Demonstraremos a seguir a conjectura proposta por Rieffel 
para o caso em que a $C^*$-\'algebra \'e o conjunto dos n\'umeros complexos.

\bigskip
Observemos que se $A = \C$, $S^{\C}(\R^n)$ \'e o conjunto das fun\c c\~oes de Schwartz $\s$, a norma $\dps \|  \cdot \dps \|_2$ \'e a norma de $L^2(\R^n)$ e $\overline{\s} = L^2(\R^n)$.

\medskip 
\noindent 
\begin{prenotac}${\cal H} = L^2(\R^n)$. \end{prenotac}

\medskip

\begin{theorem}
\label{theo1.1}.  Dadas $A \in B({\cal H})$ e $J$, matriz $n \times n$ anti-sim\'etrica, $A$ \'e 
Heisenberg-suave e $[A,R_G]=0 \,\,\,  \forall G \in \CB$ se e somente se existe $ F \in \CB  $ tal que $A = L_F$. 
\end{theorem}

Cordes fez uma caracteriza\c c\~ao dos operadores suaves pela a\c c\~ao do grupo de Heisenberg (\cite{17}): $A \in B({\cal H})$ \'e Heisenberg-suave se, e somente se, existe uma fun\c c\~ao $a \in CB^\infty(\R^{2n})$ tal que $Au(x) = 
 \int e^{i(x-y)\xi}a(x,\xi)u(y) \dc y \dc \xi$, onde $\dc y  = (2\pi)^{-n/2}dy$ e $u \in C^{\infty}_{0}(\R^n)$. A fun\c c\~ao $a$ \'e chamada de {\em s\'\i mbolo} de $A$ e o operador $A$ \'e um operador pseudo-diferencial.  Ressaltamos que a integral acima \'e uma
integral oscilat\'oria. Al\'em disso, prova-se facilmente que $Au \in \s$. O teorema de
Calder\'on-Vaillancourt \cite{10} mostra que $A$ 
pode ser estendido a um operador em $\cal{B}(H)$, \cite[Cap\'\i tulo 8, sec\c c\~ao 1]{2}.

 Pelo visto na observa\c c\~ao \ref{1umcuarto}, 
\'e f\'acil ver que o operador $L_F$, onde $F \in \CB$, 
\'e um operador pseudo-diferencial com 
s\'\i mbolo $\tilde{F} \in CB^{\infty}(\R^{2n})$ com $\tilde{F}(x,\xi) = F(x-J\xi ).$ 

%
%


\bigskip
\bigskip
Vejamos agora um resultado que estabelece uma rela\c c\~ao entre um operador pseudo-diferencial e seu s\'\i mbolo.

\medskip
Observemos primeiramente que podemos ver $\s$ como um sub-conjunto de $S'(\R^n)$, o conjunto das distribui\c c\~oes temperadas, como segue: cada $f \in \s$ determina um funcional linear  cont\'\i nuo em $\s$, $F$, dado por:

\[ <F,g> = \int fg \,\,\,\,  para \,\,\, g \in \s \] (\cite[cap.8,sec.5(i)]{3}).

\begin{predef}\label{1.5meio} {\rm Uma fun\c c\~ao $f \in C(\R^n)$ \'e dita de {\em crescimento polinomial} se  existe $N\geq 0$ tal que $(1+|x|)^{-N} f$ \'e limitada. }\end{predef}

\begin{preprop}\label{1.3} {\rm Existe um isomorfismo $T$ em $S'(\R^{2n})$ 
tal que todo $A \in \cal{B}(H)$, opera~-dor pseudo-diferencial com s\'\i mbolo $a \in CB^\infty(\R^{2n})$, satisfaz 

\[<T(a),u \otimes v> = <Av,u> \hspace{1cm} u, v \in \s ,\] }\end{preprop} 

\noindent onde $u \otimes v$
denota a fun\c c\~ao $u \otimes v(x,y) = u(x)v(y).$

\medskip
\begin{remark} Sejam $\varphi, \psi: S(\R^{2n}) \rightarrow S(\R^{2n})$ dadas por 

$$\varphi(f)(x,y) = f(x,x-y)$$
$$\psi(f)(x,y) =  (2\pi)^{-n} \dps \int e^{iy\xi}f(x,\xi)d \xi .$$ $\varphi, \psi$  s\~ao isomorfismos.

\medskip Seja $T = \varphi \circ \psi .$

\[ <T(f),g> = \dps \int \varphi \circ \psi (f)(x,y)g(x,y) dxdy = \]
\[ =  \dps \int \psi(f)(x,x-y)g(x,y) dxdy = \dps \int \psi(f)(x,z)g(x,x-z)dxdz = \]
\[ = \dps \int (2\pi)^{-n} \dps \left( \dps \int e^{iz \cdot \xi}f(x,\xi)d\xi \dps \right) \varphi (g)(x,z) \dps  dxdz = \]
\[ = \dps \int f(x,\xi ) \dps \left( (2\pi)^{-n} \dps \int e^{iz\xi} \varphi (g)(x,z)dz \dps \right) dx d\xi = \]
\[ = \dps \int f(x,\xi) \psi \circ \varphi (g) (x,\xi)dx d \xi = <f,\psi \circ \varphi (g)> .\]

Se $\tilde{T} = \psi \circ \varphi$, temos $<T(f),g>  = <f,\tilde{T}(g)>$ . Ent\~ao, $T$ se estende a um isomorfismo em $S'(\R ^{2n})$ (ver \cite[cap.8,sec.5]{3}), que tamb\'em chamamos $T$,
definido como segue. Para $Q \in S'(\R^{2n}), $
\[ <T(Q), g > = <Q,\tilde{T}(g)>   .\] Notemos que $T$ \'e um
operador cont\'\i nuo em $S'(\R^{2n})$ com a topologia fraca *.

\bigskip Dada $a \in CB^{\infty}(\R ^{2n})$, $a$ tem crescimento polinomial;
logo $a$ define uma distribui\c c\~ao 
em $S'(\R ^{2n})$ (\cite[V.3]{4}).

\bigskip Para $u,v \in \s$,

\[ <T(a), u \otimes v> = <a, \tilde{T}(u \otimes v)> .\]

 E \[   \tilde{T}(u \otimes v)(x,y) = \psi \circ \varphi (u \otimes v)(x,y) = (2\pi)^{-n} \dps \int e^{i y\xi}\varphi (u \otimes v)(x,\xi)d\xi = \]
  \[ = (2\pi)^{-n} \dps \int e^{i y \xi}(u \otimes v)(x,x- \xi)d\xi =
(2\pi)^{-n} \dps \int e^{i y(x-\xi)}(u \otimes v)(x,\xi)d\xi  = \] \[ =
 (2\pi)^{-n}e^{i yx}u(x) \dps \int e^{-i y \xi}v(\xi)d\xi .\]
Ent\~ao \[ <a,\tilde{T}(u \otimes v)> = \dps \int a(x,y)(2\pi)^{-n}e^{iyx}u(x) \dps 
\int e^{-iy \xi}v(\xi)d\xi dx dy = \]
\[ =  \dps \int \left(\dps \int e^{iy(x-\xi)}a(x,y)v(\xi)\dc  \xi \dc y \right) u(x)dx = <Av,u> .  \ai \] 
\end{remark}

\begin{precorollary}\label{cor:2.3quarto}
{\rm  Dado $A \in \cal{B}(H)$, operador pseudo-diferencial com s\'\i mbolo $a \in \CB $, 
a aplica\c c\~ao $a \rightarrow A$ \'e injetora.}    \end{precorollary}

\begin{remark}  Pela proposi\c c\~ao \citeprop{1.3}, temos que existe um
isomorfismo $T$ em $S'(\R^{2n})$ tal que, para $u, v \in \s$ vale
$<T(a),u \otimes v> = <Av,u>.$ Logo, se $A = 0$,
$<T(a), u \otimes v> = 0$.  Como $\s \otimes  \s$ \'e denso em $S(\R^{2n})$, 
\cite[teorema 51.6]{8}, ent\~ao temos que $T(a) = 0$, o que implica $a = 0$.  $\ai$
\end{remark}

\begin{prenota}\label{una} {\rm Citaremos a seguir, na proposi\c c\~ao \citeprop{1.3ummeio}, 
o lema 2 de 
\cite{14}, que afirma que existe uma seq\"u\^encia $(e_k)_{k \in \N}$ em $\s$ tal que
$L_{e_k}$ converge para a identidade na topologia forte de operadores. 
}\end{prenota}


%
%

\begin{prenotac}  {\rm Seja $\cal A$ o conjunto dos operadores Heisenberg-suaves, $A$,
tais que $$[A,R_G] = 0  \hspace{1cm} \forall G \in \CB .$$ Pela caracteriza\c c\~ao
de Cordes, todo $A\in{\cal A}$ \'e da forma $A=a(x,D)$ para alguma $a\in CB^\infty(\R^{2n})$.
} \end{prenotac}

\begin{preprop}\label{1.3ummeio}  {\rm Dado $A \in \cal{A}$, existe uma seq\"u\^encia $(f_k)_{k \in \N}$ em $\s$ de modo que $L_{f_k}$ convege para $A$ na
topologia forte de operadores. } \end{preprop}

\begin{remark} Seja $(e_k)_{k \in \N}$ uma seq\"u\^encia em $\s$ tal que $L_{e_k}$ converge para a identidade na topologia forte de operadores \cite[Lema 2]{14}. Como $A$
mant\'em $\s$ invariante, $A(e_k) \in \s$. Seja $f_k = A(e_k)$. 

Dada $u \in \s$,
\[L_{f_k}u = f_k \times_J u = Ae_k \times_Ju = R_uAe_k = AR_ue_k = AL_{e_k}u .\]

Assim, $L_{f_k} = AL_{e_k}$, pois $\s$ \'e denso em $\cal{H}$; e, como $L_{e_k}$
converge para $I \in \cal{B}(H)$ na topologia forte de operadores, ent\~ao
$L_{f_k}$ converge para $A$ nesta topologia.  $\ai$ \end{remark}

\begin{preprop}\label{1.4} {\rm Dado $A \in \cal{A}$, seja $(F_k)$ uma seq\"u\^encia
de fun\c c\~oes em $\CB$ de modo que $A$ \'e limite de $L_{F_k}$ na 
topologia forte de operadores. Se $a$ \'e o s\'\i mbolo de $A$ e $a_k$ \'e o s\'\i mbolo de 
$L_{F_k}$, ent\~ao temos $a_k \longrightarrow a$ em $S'(\R^{2n})$. }\end{preprop} 

\begin{remark} Por hip\'otese vale que $$L_{F_k}v \longrightarrow Av \hspace{1cm} \forall v \in H .$$  Logo, $$<L_{F_k}v,u> \: \longrightarrow \: <Av,u> \hspace{1cm} \forall u,v \in \s .$$  Logo, pela proposi\c c\~ao \citeprop{1.3}
$$<T(a_k),u \otimes v> \: \longrightarrow \: <T(a),u \otimes v>. $$
ou seja, $<T(a_k)-T(a),u \otimes v> \longrightarrow_{k \rightarrow \infty} 0$. 

Queremos provar que, para toda $h \in S(\R^{2n})$, $<T(a_k)-T(a),h> \longrightarrow_{k \rightarrow \infty}0.$ 

Primeiramente, observemos que, pelo Princ\'\i pio da Limita\c c\~ao Uniforme, $\supp \|L_{F_k}\| < \infty$ (ver \cite[teorema 2.6]{6}). 

Al\'em disso, vemos em 
\cite[teorema 51.6]{8} que $\s \otimes \s$ \'e denso em $S(\R^{2n})$. Pelo corol\'ario
anterior ao teorema 51.6, temos que $\s$ \'e nuclear; logo, pelo feito na segunda p\'agina
do cap\'\i tulo 50, temos que $\s \otimes \s = \s \otimes_{\pi} \s$, onde $\otimes_\pi$ \'e
definido em \cite[cap\'\i tulo 43, defini\c c\~ao 43.2]{8}. 

Assim, pelo teorema 45.1 de \cite{8}, temos que dada
$h \in S(\R^{2n})$, existem seq\"u\^encias $(u_j)_{j \in \N}$ e $(v_j)_{j \in \N}$ em $S(\R^{2n})$,
convergindo a zero, e $(\lambda_j)_{j \in \N}, \, \lambda_j \in \C$, tal que 
$\sum_j|\lambda_j| < \infty$, que satisfazem $$h = \sum_j\lambda_ju_j\otimes v_j .$$

\noindent Como visto acima, temos que, para todo $j$, 
$$\lim_{k \rightarrow \infty} <T(a_k)-T(a),u_j\otimes v_j> \,\, = \,\, 0.$$ E, considerando que:

\begin{enumerate}
\item $<T(a_k)-T(a),u_j\otimes v_j> = <(L_{F_k}-A)v_j,u_j>$,
\item $(u_j)$ e $(v_j)$ convergem a zero em $\s$ e, portanto, convergem a zero em $L^2(\R^n)$,
\item $\supp_{k}\|L_{F_k}\| < \infty$,  \end{enumerate} temos que 
$\supp _{k,j}\left|<T(a_k)-T(a),u_j \otimes v_j>\right| < \infty$. Portanto, pelo Teorema da Converg\^encia Dominada para s\'eries, $\lim_{k \rightarrow \infty} \sum_j\lambda_j<T(a_k)-T(a),u_j\otimes v_j> \,\, = \,\, 0.$  

Isto \'e, $\lim_{k \rightarrow \infty} <T(a_k)-T(a),h> \,\, = \,\, 0.$ Como $T$ \'e isomorfismo
em $S'(\R^{2n})$, temos ent\~ao que $(a_k)$ converge para $a$ em $S'(\R^{2n})$. $\ai$
\end{remark}

\begin{definition}\label{1.1} Dadas fun\c c\~oes $a, b \in C^{\infty}(\R^{2n})$, definimos seu
{\em colchete de Poisson} como segue:

\[ \{a,b\} = \dps \sum^{n}_{j=1} \left( \frac{\partial a}{\partial x_j} 
\frac{\partial b}{\partial \xi_j} - \frac{\partial a}{\partial \xi_j}
\frac{\partial b}{\partial x_j} \right) .\] \end{definition}

\begin{preprop}\label{1.5} {\rm Dadas $F, G \in \CB$ e $J$ matriz anti-sim\'etrica $n \times n$, sejam $$a(x,\xi) = F(x-J\xi)$$ e $$ b(x,\xi) = G(x+J\xi) .$$ Ent\~ao, $\{a,b\} = 0$. }\end{preprop}

\begin{remark} Fazendo contas, temos que:

\[ \frac{\partial a}{\partial x_i} = \partial_{i}F \hspace{0.7cm} \frac{\partial a}{\partial \xi_i} = \dps \sum^{n}_{j=1} J_{ij}\partial_j F \]

\[ \frac{\partial b}{\partial x_i} = \partial_i G  \hspace{0.7cm} \frac{\partial b}{\partial \xi_i} = - \dps
\sum^{n}_{j=1}J_{ij}\partial_j G \]

\[ \frac{\partial a}{\partial x_i}\frac{\partial b}{\partial \xi_i } - 
\frac{\partial a}{\partial \xi_i} \frac{\partial b}{\partial x_i} = 
 - \partial_i F \dps \sum^{n}_{j=1} J_{ij}\partial_j G - 
\partial_i G \dps \sum^{n}_{j=1} J_{ij}\partial_j F = 
 \partial_i F \dps \sum^{n}_{j=1} J_{ji} \partial_j G + \partial_i G \dps \sum^{n}_{j=1} J_{ji}
\partial_j F .\]

Assim, $$\{ a,b \} = \dps \sum^{n}_{j=1} \dps \sum^{n}_{j=1}(J_{ji} \partial_i F\partial_j G + J_{ji} \partial_i G \partial_j F)$$ e $\{a,b\} = 0$ pois $J$ \'e anti-sim\'etrica. $\ai$ \end{remark}

\begin{prenotac} {\rm Seja ${\cal P}_B(\R ^{n})$ o conjunto das fun\c c\~oes em $C^\infty(\R^n)$ tais que elas e todas as suas
derivadas s\~ao de crescimento polinomial.} \end{prenotac}

\begin{preprop}\label{1.6} {\rm Dadas $A \in \cal{A}$, com s\'\i mbolo $a \in CB^\infty (\R ^{n})$, 
e a fun\c c\~ao $b(x,\xi) = G(x + J\xi)$, com $G \in {\cal{P}}_B(\R^{2n})$, ent\~ao
$\{a,b\}= 0$.}
\end{preprop}

\begin{remark} Dado $A \in \cal{A}$, existe uma seq\"u\^encia$(A_k)_{k \in \N}$, onde $A_k = L_{F_k}$, com $F_k \in \s$, de modo que $(A_k)_{k \in \N}$ converge para $A$ na topologia forte de operadores, \ref{1.3ummeio}. Se $a_k$ \'e o s\'\i mbolo de $A_k$, vimos que $(a_k)_{k \in \N}$ converge para $a$ em $S'(\R ^{2n})$  \citeprop{1.4}. Como o s\'\i mbolo de $A_k$ \'e $F_k (x-J\xi)$, temos pela proposi\c c\~ao anterior que $\{a_k,b\} = 0, \hspace{0,5cm} \forall k$. 

\medskip
Como $b \in {\cal P}_{B}(\R ^{2n})$, $\{a_k,b\} \in S'(\R ^{2n}), \, \forall k$. 
Como $a_k \longrightarrow a$ em $S'(\R ^{2n})$, 
e como a multiplica\c c\~ao por uma derivada qualquer da fun\c c\~ao $b$ \'e uma opera\c c\~ao
cont\'\i nua em $S'(\R^{2n})$ (\cite[V.3]{4}), temos que :

\[ \frac{\partial b}{\partial \xi_i}\frac{\partial a_k}{\partial x_i} \: \longrightarrow \: \frac{\partial b}{\partial \xi_i} \frac{\partial a}{\partial x_i}\] 

e

\[\frac{\partial b}{\partial x_i} \frac{\partial a_k}{\partial \xi_i} \: \longrightarrow \: 
\frac{\partial b}{\partial x_i}\frac{\partial a}{\partial \xi_i}  \] em $S'(\R ^{2n}) $,  e vale que $\{a_k,b\} \: \longrightarrow \: \{a,b\}$ em $S'(\R ^{2n})$. 
Como $\{a_k,b\} = 0 \hspace{0,5cm} \forall k$ temos que $\{a,b\} = 0.$ $\ai$ \end{remark}

\begin{lemma}\label{1.7} Seja  $a \in CB^{\infty}(\R ^{2n})$ que satisfaz $\{a,b\} = 0$ para qualquer $b \in {\cal P}_B (\R ^{2n})$ com $b(x,\xi) = G(x + J\xi)$. Ent\~ao existe $F \in  \CB$ tal que $a(x,\xi) = F(x- J\xi)$. \end{lemma}

\begin{remark} Seja $b_i$ a proje\c c\~ao de $x + J\xi$ na coordenada $i$,

\[ b_{i}(x,\xi) = x_i +\dps \sum_{k=1}^{n} J_{ik}\xi_k .\] Assim definida, $b_i \in {\cal P}_B(\R ^{2n})$ e, pela proposi\c c\~ao anterior, temos $\{a,b_i\} = 0$, para todo $i = 1, \cdots , n$. Logo,

\[ \dps \sum_{j=1}^{n}(\frac{\partial a}{\partial x_j} J_{ij} - \frac{\partial a}{\partial \xi_j} \delta_{ij}) = \dps \sum_{j=1}^{n} J_{ij}\frac{\partial a}{\partial x_j} - \frac{\partial a}{\partial\xi_i} = 0 .\] 

Ou seja, \[\frac{\partial a}{\partial \xi_i} = \dps \sum_{j=1}^{n} J_{ij} \frac{\partial a}{\partial x_j}  \hspace{1cm} i = 1, \cdots , n  .   \hspace{1cm}  (*) \]

\noindent Portanto, aplicando as id\'eias de \cite[2.1]{16}, temos que \[ a(x,\xi) = a((x,\xi)+s(-J_{i1}, \cdots
, -J_{in},e_i)) , \hspace{1cm} (**) \] para todo $(x,\xi) \in \R^{2n}$, para todo $s \in \R$, para
qualquer $i \in \{1,2,\cdots ,n\}$, onde $\{e_i\}$ \'e uma base ortonormal de $\R^n$, $e_i = (0,\cdots ,0,1,0,\cdots ,0)$.

Agora, em $(**)$, consideremos $s_1$ em lugar de $s$, $i=1$, e $(x.\xi) = (x_0,0)$:

\[ a(x_0,0) = a((x_0,0) + s_1(-J_{11}, \cdots ,-J_{1n},e_1)) = a((x_0,0)+(-J_{11}s_1,\cdots ,-J_{1n}s_1,s_1e_1)) .\]  Apliquemos $(**)$ novamente para $s_2$ no lugar de $s$, $i=2$  e  $(x,\xi) = (x_0-(J_{11}s_1,\cdots ,J_{1n}s_1),$ $s_1e_1)$. Assim, temos

\[ a(x_0,0) = a(x_0 + (-J_{11}s_1,\cdots ,-J_{1n}s_1) + (-J_{21}s_2,\cdots ,-J_{2n}s_2),
s_1e_1 + s_2e_2) .\]

Repetindo o racioc\'\i nio para $i = 3,\cdots ,n$,
temos 
\[ a(x_0,0) = a(x_0 -(\sum^n_{i=1} J_{i1}s_i, \cdots , \sum^n_{i=1}J_{in}s_n),\sum^n_{i=1}
s_ie_i). \] Como $J$  \'e anti-sim\'etrica, tomando $\xi = (s_1,\cdots ,s_n)$, temos:

\[ a(x_0,0) = a(x_0 + J\xi,\xi). \]

Seja $F(x) = a(x,0).$  Ent\~ao $a(x,\xi) = F(x-J\xi)$.
$\ai$ \end{remark}

\begin{precorollary}\label{cor:1.8} {\rm Dado $A \in \cal{A}$, existe $F \in \CB$ tal que $A = L_F$. } \end{precorollary}

\begin{remark} Dado $A \in {\cal{A}}$, pela caracteriza\c c\~ao feita por Cordes dos operadores Heisenberg-suaves, $A$ \'e um operador pseudo-diferencial. Seja $a \in CB^\infty(\R^{2n})$ seu s\'\i mbolo. Pela proposi\c c\~ao \citeprop{1.6}, dada $b(x,\xi) = G(x +J\xi)$, com $G \in {\cal{P}}_B(\R^n)$, temos $\{a,b\}=0$.
Assim, pelo lema \citelemma{1.7}, existe $F \in \CB$ tal que $a(x.\xi) = F(x - J\xi)$. Temos, ent\~ao, que o operador $A$ tem s\'\i mbolo $F(x - J\xi)$, ou seja, $A = L_F$. $\ai$\end{remark}

\bigskip
Com isto fica provado o Teorema \citetheo{theo1.1} e portanto a quest\~ao proposta por Rieffel, no
caso da $C^*$-\'algebra ser o conjunto dos n\'umeros complexos.

\vfill\break

\chapter{ }


\vspace{1cm}

Seja $A$ uma C*-\'algebra separ\'avel ao longo deste cap\'\i tulo. 
Consideremos,
novamente, a a\c c\~ao do grupo de Heisenberg em $\op$ , 
\[ (z,\zeta) \rightarrow 
E_{z,\zeta}^{-1}TE_{z,\zeta} \] onde $T \in \op$ e $E_{z,\zeta}u(x) = 
e^{i\zeta x}u(x-z)$ se
$u \in \sa $.

\bigskip \noindent
\begin{prenotac}\label{2.0} $T_{z,\zeta}= E_{z,\zeta}^{-1}TE_{z,\zeta}$. \end{prenotac}

\begin{predef} {\rm Dizemos que $T \in \op$ \'e {\it Heisenberg-suave} se a 
aplica\c c\~ao $(z,\zeta) \mapsto T_{z,\zeta}$ \'e $C^\infty$. Dizemos
que $T \in \op$ \'e {\it suave por transla\c c\~oes} se a aplica\c c\~ao
$z \mapsto T_{z,0}$ \'e $C^\infty$.} \end{predef}

Neste cap\'\i tulo  veremos que dado um operador $T$ em $\op$ que comuta com qualquer operador $R_g$ para $g$ em $\sa $ (ver cap\'\i tulo 1), se $T$ \'e suave por transla\c c\~oes ent\~ao \'e
Heisenberg-suave.

\begin{prelemma}\label{2.1}  {\rm Seja $f: \R^n \rightarrow A$ uma fun\c c\~ao cont\'\i nua 
tal que $\sup_{x \in \R^n} |x|^l \|f(x) \| < \infty$ para todo inteiro positivo $l$. Ent\~ao, dada $J$, matriz anti-sim\'etrica $n \times n$, para cada $l \in \Z^+$ e para todo $ \epsi >0$, existe $\delta > 0$ de modo que \[|x|^l \|f(x+J \xi)e^{ix \xi} - f(x) \| < 
\epsi\] para quaisquer $x,\xi \in \R^n$ com $|\xi | < \delta $.} \end{prelemma}

\begin{remark} Dados $\epsi > 0, \,\, l \in \Z^+$, existe $R_l \in \N$ tal que se $| x | > R_l$,
$$| x |^{l} \| f(x) \| < \frac{\epsi}{2}.  \hspace{1cm} (1)$$ 

Seja $d_J, \, \, 0 < d_J \leq  1 $,  tal que se $|\xi |< d_J$,
$|J\xi | < 1$. 
Consideremos
\[ X_{f,l} = \{ (x,\xi ) \in \R ^{2n}, |x| \leq R_l+1, | \xi | \leq d_J \} .\] Embora a nota\c c\~ao
fique um pouco "carregada", consideramos que \'e necess\'ario para situar mais
facilmente o leitor.

\medskip
\noindent \underline{1$^{o}$ passo:} \, $l = 0$. 

\bigskip
Se $g(x,\xi) = f(x+J \xi)e^{ix \xi}$, $g$ \'e uniformemente cont\'\i nua em $X_{f,0}$, 
pois $X_{f,0}$ \'e
compacto. Logo, dado $\epsi > 0$ existe $\delta > 0, \,\,\delta \leq d_J$,
tal que se $|(x,\xi) - (x',\xi ')| < \delta$ e $(x,\xi),(x',\xi') \in X_{f,0}$, temos 
$\|g(x,\xi) - g(x',\xi ') \| < \epsi$.
Seja, agora, $(x',\xi ') = (x,0)$. Assim $|(x,\xi ) - (x',\xi ')| = |\xi |$. Portanto, se $|\xi | <
\delta, \,\,\|f(x+J\xi )e^{ix \xi } - f(x) \| < \epsi$ para todo $x \in \R^n$ com
$|x| \leq R_0+1$. 

Se $|x| > R_0+1$ e $|\xi | < \delta$, $|x+J\xi | > R_0$ pois
$|J\xi| < 1$, e temos, por (1):
$$\| f(x+J\xi )e^{ix \xi} - f(x) \| \leq \|f(x+J\xi ) \| + \|f(x)\| < \frac{\epsi}{2} + \frac{\epsi}
{2} = \epsi .$$

\bigskip
\noindent \underline{2$^{o}$ passo:} \, $l$ qualquer.

\bigskip
Seja $h(x) = |x|^l f(x)$. Pelo primeiro passo, aplicado a $h$, temos que existe 
$\sigma_l \leq d_J$ tal que, se $|\xi| < \sigma_l$, ent\~ao


%
%
%
%
%
%

\[ \| |x+J\xi|^lf(x+J\xi)e^{ix\xi} - |x|^lf(x)\| < \epsi . \hspace{1cm} (2) \]

E notemos que, por $(2)$, se $|\xi| < \sigma_l$, 
\[ |x|^l\|f(x+J\xi)e^{ix\xi} - f(x) \| \leq \| |x+J\xi|^lf(x+J\xi)e^{ix\xi}-|x|^lf(x)\| + 
 \||x|^l f(x+J\xi)e^{ix\xi} - \] \[|x+J\xi|^lf(x+J\xi)e^{ix\xi} \| < \epsi + ||x|^l - |x+J\xi|^l|\|f(x+J\xi)\|. \hspace{1cm}
(3)\] 

\medskip
Se $(x,\xi) \in X_{f,l}$, \,\, $|x+J\xi| \leq R_l + 2$, pois $|J\xi| < 1$.

Seja $K = \supp_{|y|\leq R_l+2}\|f(y)\|$.
Se $p(x,\xi) = |x+J\xi|^{l}$, $p$ \'e uniformemente cont\'\i nua em $X_{f,l}$. Assim, dado $\epsi > 0$, existe $\delta _l> 0$, $\delta_l \leq \sigma_l$ ($\delta_l$ depende de $l$), tal que se 
$(x,\xi), (x',\xi') \in X_{f,l}$ e $|(x,\xi) - (x',\xi ')| < \delta_l$ ent\~ao $|p(x,\xi)-p(x',\xi ')| < 
\frac{\epsi}{K}$. Tomando $(x',\xi ') = (x,0)$, temos  $| |x+J\xi|^{l}-|x|^{l}| < 
\frac{\epsi}{K}$ para $|\xi| < \delta_l$ e $|x| \leq R_l + 1$. 

Logo, se $|\xi| < \delta_l$, temos $|x|^{l}\|f(x+J\xi)e^{ix\xi}
- f(x)\| < 2\epsi$ para $|x| \leq R_l+1$. 

Se $|x| > R_l+1$, \,\, e $|\xi| < \delta_l \leq \sigma_l \leq d_J,$
temos que $|x+J\xi| > R_l$  pois $|J\xi| < 1$.
Sem perda de generalidade, consideremos $R_l > 1$ e teremos: $|J\xi| < 1 < R_l < |x+J\xi|$.

Assim 

\[ |x| = |x+J\xi-J\xi| \leq |x+J\xi| + |J\xi| \leq 2|x+J\xi|. \hspace{1cm} (4) \]

O lado direito de (3) \'e, portanto, menor ou igual a

%
%
\[  \epsi + |x|^l\|f(x+J\xi)\| + |x+J\xi|^l\|f(x+J\xi)\| \leq \]

\[ \leq \epsi + 2^l|x+J\xi|^l\|f(x+J\xi)\| + |x+J\xi|^l\|f(x+J\xi)\| < \epsi + (2^l+1)\frac{\epsi}{2}, \] por 
$(1)$, j\'a que $|x+J\xi| > R_l$.

\medskip Ou seja, dado $\epsi > 0$, existe $\delta_l > 0$ tal que
\[ |x|^l\|f(x+J\xi)\|e^{ix\xi}-f(x)\| < \epsi(2^{l-1}+\frac{3}{2}) \] para $x,\xi \in \R^n$ com 
$|x| > R_l + 1$ e  $|\xi| < \delta_l$.

\medskip Para todo $\epsi' >0$, seja $\epsi >0$ de modo que $\epsi(2^{l-1} + \frac{3}{2}) 
< \epsi'$, e $2\epsi < \epsi'$. 

Assim, temos que, para qualquer $\epsi'>0$, existe $\delta_l > 0$ tal que 
$|x|^l\|f(x+J\xi)e^{ix\xi} - f(x)\| < \epsi'$, para todo $x \in \R^n$ e $\xi \in \R^n$ com $|\xi| < \delta_l$, tendo
fixado $l \in \N$.
$\ai$ \end{remark}

\begin{preobs}\label{2.1ummeio} {\rm   Usaremos o lema \citelem{2.1} na demonstra\c c\~ao da
proposi\c c\~ao \citeprop{2.2} do seguinte modo: estaremos considerando a fun\c c\~ao 
$s(x) = (1+|x|)^n$ e vamos provar que, dado $\epsi_1 > 0$, existe $\delta > 0$
tal que, se $|\xi| < \delta$, $s(x)\|f(x+J\xi)e^{ix\xi} - f(x)\| < \epsi_1$. Notemos que $s(x)$ \'e uma
soma finita de termos do tipo $c_l|x|^l$, onde $c_l \in \Z^+$  \'e uma constante, 
e $l \in \{0,1,\cdots ,n\}$. 
Dado qualquer $\epsi_1 > 0$, escolhemos $\epsi >0$ de modo que 
$\dps\sum c_l\epsi < \epsi_1$.
Para cada valor de $l$, aplicamos o lema \citelem{2.1}. Isto demonstra a afirma\c c\~ao.}
\end{preobs}

%

\begin{proposition}\label{2.2}  Dada $J$, matriz anti-sim\'etrica $n \times n$,
existe uma seq\"u\^encia $(e_k)$ em $\sa$ tal que, para toda $f \in \sa$, temos $\lim_{k \rightarrow 0} L_{e_{k}}(f) = f$ em $\sa$ com $\| \, \|_2$. \end{proposition}

\begin{remark} Como $A$ \'e separ\'avel, admite uma aproxima\c c\~ao da unidade que \'e uma 
seq\"u\^encia (\cite[Remark(3.11)]{5}).
Seja $(u_k)_{k \in \N}$ uma aproxima\c c\~ao da unidade para $A$, com $u_k 
\geq 0$ e $\| u_k \| \leq 1$ para todo $k \in \N $. 

Seja $\psi \in C^{\infty}_{c}(\R^n) , \psi \geq 0$ 
tal que o suporte de $\psi$ est\'a contido na bola unit\'aria, $({\rm supp}(\psi) \subseteq B(0,1))$, e
$\dps \int_{\R^n} \psi(\xi)d \xi = 1$. Para $k \in \Z^+$, seja $\psi_k(\xi) =
k^n\psi(k\xi)$; assim ${\rm supp}(\psi_k) \subseteq B(0,\frac{1}{k})$. Seja $\phi_k = \psi_k u_k$.
Assim ${\rm supp}(\phi_k) \subseteq B(0,\frac{1}{k})$ e $\dps \int_{\R^n}\phi_k(\xi)d \xi  = u_k$.
Seja $\varphi_k$ dada por $\dps \int _{\R^n}e^{-iy\xi}\varphi_k(y) \dc y = \phi_k(\xi)$.
Como $\phi_k \in C^\infty_{c}(\R^n,A)$, $\varphi_k \in \sa$. 
Como feito na proposi\c c\~ao, \ref{asterisco}, 
e, por \cite[(1.13)]{1}, para $f \in \sa$ temos que $R_f(\varphi_k)(x)=\dps\int e^{ix\xi}\phi_k(\xi)f(x+J\xi)\dc\xi$. 

\medskip Como $\dps \int \phi_k(\xi)f(x)d \xi  = u_kf(x)$, temos: 

\[ \|(2\pi)^{\frac{n}{2}}R_f(\varphi_k)(x)-u_kf(x) \|  \leq \dps \int \|e^{ix\xi}\phi_k(\xi)f(x+J\xi) - \phi_k(\xi)f(x) \|d\xi  \leq \]
\[ \leq \dps \int \|\phi_k(\xi)\| \|e^{ix\xi}f(x+J\xi) - f(x) \|d\xi  . \]



Pela observa\c c\~ao \citeobs{2.1ummeio} temos que, para todo $\epsi > 0$ existe $\delta > 0$
tal que, se $k > \frac{1}{\delta}$: 
\[s(x) \|(2\pi)^{\frac{n}{2}}R_f(\varphi_k)(x) - u_kf(x) \| \leq \epsi\dps\int_{B(0,\frac{1}{k})}  
\|\phi_k(\xi)\| d\xi  < \epsi , \] pois
\begin{center} $ \dps\int_{B(0,\frac{1}{k})} \|\phi_k(\xi)\|d\xi  = 
\dps\int_{B(0,\frac{1}{k})}\psi_k(\xi)d\xi \|u_k\| \leq 1,$  para $   k > \frac{1}{\delta}. $
\end{center}

Portanto  \[\|(2\pi)^{\frac{n}{2}}R_f(\varphi_k) - u_kf \|^2_2 \leq \dps\int_{\R^n} \|(2\pi)^{\frac{n}{2}}
R_f(\varphi_k)(x) - u_kf(x)\|^2dx \leq \epsi^2\dps\int_{\R^n}\frac{1}{s(x)^2}dx. \hspace{1cm} (1) \]

Por outro lado, $\|f - u_kf \|^2_2  \leq  \dps \int \|(f-u_kf)(x)\|^2 dx .$
Como $u_k$ \'e aproxima\c c\~ao da unidade em $A$, para cada $x \in \R^n$ temos que
$\lim_{k \rightarrow \infty} \|f(x) - u_kf(x) \| = 0$. Al\'em disso, $\|f(x) - u_kf(x) \| \leq 2\|f(x)\|$. Como
$2\|f(x)\| \in L^1(\R^n)$, pelo Teorema da Converg\^encia Dominada, 
temos que $\lim_{k \rightarrow \infty} \dps \int \|(f - u_kf)(x)\|^2dx = 0, $ o que implica que
\[ \lim_{k \rightarrow \infty} \|f - u_kf \|_2 = 0.  \hspace{1cm} (2) \]

Portanto, por $(1)$ e $(2)$, $\lim_{k \rightarrow \infty} \|(2\pi)^{\frac{n}{2}}R_f(\varphi_k) - f\|_2 = 0.$

Tomando $e_k = (2\pi)^{-\frac{n}{2}}\varphi_k$, temos que $\lim_{k \rightarrow \infty} \|R_f(e_k)
- f \|_2 = 0$.

\medskip  \noindent
Observamos que $R_f(e_k) = e_k \times_J f  = L_{e_k}(f)$. Portanto, temos que para toda
$f \in \sa, \hspace{0,5cm} \lim_{k \rightarrow \infty}L_{e_k}(f) = f$ em $(\sa, \| \cdot \|_2)$.
$\ai$ \end{remark}

\vspace{1cm}Usamos agora argumentos semelhantes a \cite[sec\c c\~ao 1]{9}.

\begin{predef}\label{2.2meio} {\rm Dada $F \in \SA$, para cada $g \in \sa$ definimos
$L_F(g) = R_g(F)$. Como $R_g$ \'e cont\'\i nuo em ($\sa$,$\| \cdot \|_2$), se estende a $\SA$, logo $L_F$ est\'a bem definido, como operador de $\sa$ em $\SA$. } \end{predef}

\begin{precorol}\label{2.3} {\rm Dado $T \in \op$ que comuta com $R_g$ para qualquer
$g \in \sa$, existe $(F_k)$ em $\SA$ tal que $\lim_{k \rightarrow \infty}L_{F_k}(g) = T(g),$
 $\forall g \in \sa$.}  \end{precorol}

\begin{remark} Pela proposi\c c\~ao \citeprop{2.2} existe uma 
seq\"u\^encia $(e_k)$ em $\sa$ tal que $\lim_{k \rightarrow \infty}L_{e_k}(g) = g,
\forall g \in \sa$. 

Por outro lado, temos que $L_{Te_k}g = R_gTe_k = TR_ge_k = TL_{e_k}g$. Portanto, $$\lim_{k  \rightarrow \infty}L_{T_{e_k}}(g) = T(g), \,\,\,  \forall g \in \sa .$$
Tomando $F_k = Te_k$ temos a tese. $\ai$ \end{remark}

\begin{preobs} {\rm Notemos que o operador $L_{F_k}$ acima \'e cont\'\i nuo para todo 
$k \in \Z^+$, pois $\|L_{F_k}\| = \|L_{Te_k}\| = \|TL_{e_k}\| \leq \|T\|\|L_{e_k}\| $. }
\end{preobs}

\begin{proposition}\label{2.4} {\rm Sejam $U, V \in \op$ operadores que preservam $\sa$. Supo\-nhamos que para toda $h \in \sa$ vale $U^*L_hU = V^*L_hV$. Ent\~ao,
dada $F \in \SA$, temos $U^*L_FU = V^*L_FV.$ Aqui   \, estamos
vendo os operadores $U^*L_FU$ e $V^*L_FV$ como operadores de $\sa$
em $\SA$, n\~ao necessariamente limitados.}\end{proposition}

\begin{remark} 
Dada $F \in \SA$, existe uma seq\"u\^encia $(h_k)$ em $\sa$  tal que $\lim_{k \rightarrow \infty} h_k = F$  em $\SA$. Dada $g \in \sa$, como $R_g$ \'e cont\'\i nuo em $\SA$,
temos \[ L_F(g) = R_g(F) = \lim_{k \rightarrow \infty}R_g(h_k). \]
Calculemos (lembrando que $U(g) \in \sa$ e $V(g) \in \sa$):
\[ U^*L_FU(g)= U^* \lim_{k \rightarrow \infty} R_{U(g)}(h_k) = \lim_{k \rightarrow \infty} U^*R_{U(g)}(h_k) = \]
\[ =  \lim_{k \rightarrow \infty}U^*L_{h_k}U(g)   = \lim_{k \rightarrow \infty}V^*L_{h_k}V(g) 
= V^*\lim_{k \rightarrow \infty}R_{V(g)}(h_k) = V^*L_F V(g).  \,\, \ai  \] \end{remark}

\begin{prelemma}\label{2.5} {\rm Dada $h \in \sa$ vale que $(L_h)_{z,\zeta} =
(L_h)_{z - J\zeta,0}.$} \end{prelemma}

\begin{remark} Lembremos que $(L_h)_{z,\zeta} = E_{z,\zeta}^{-1}L_hE_{z,\zeta}$, onde
$E_{z,\zeta} = M_{\zeta}T_z$, com $M_{\zeta}u(x) = e^{i\zeta x}u(x)$, $
T_z(u)(x) = u(x-z)$ para $u \in \sa$ (ver cap\'\i tulo 1 e nota\c c\~ao (\ref{2.0})). Al\'em disso,
\[ E_{z,\zeta}^{-1}u(x) = T_{-z}M_{-\zeta}u(x) = e^{-i\zeta(x+z)}u(x+z).\]


Dada $g \in \sa$, 
\[\widehat{E_{z,\zeta}(g)}(\xi) = \dps \int e^{-i\xi y}E_{z,\zeta}g(y)\dc y = \]
\[ =  \dps \int e^{-i\xi y}e^{i\zeta y}g(y-z)\dc y = \dps \int e^{-i(\xi - \zeta)(y+z)}g(y)\dc y = e^{-iz(\xi - \zeta)}\hat{g}(\xi - \zeta) . \]

\medskip \noindent Portanto, por \cite[(3.1)]{1}
\[ (L_h)_{z,\zeta}g(x) = e^{-i\zeta (x+z)}L_hE_{z,\zeta}g(x+z) 
 = e^{-i\zeta(x+z)}\dps \int e^{i(x+z)\xi}h(x+z - J\xi)\widehat{E_{z,\zeta}g}(\xi) = \]
\[ = e^{-i\zeta(x+z)}\dps \int e^{ix\xi} e^{i\zeta z}h(x+z - J\xi)\hat{g}(\xi - \zeta) \dc \xi  = e^{-i\zeta x}\dps \int e^{ix(\xi +\zeta)}h(x+z - J(\xi + \zeta))\hat{g}(\xi) \dc \xi = \]
\[ = \dps \int e^{ix\xi}h\left( x + z - J\zeta - J\xi\right) \hat{g}(\xi)\dc \xi .\]
Isto \'e, \[ (L_h)_{z,\zeta}g(x) = \dps \int e^{ix\xi}h(x+z-J\zeta-J\xi)\hat{g}(\xi)\dc \xi \] para qualquer $(z,\zeta) \in \R^{2n}.$

\noindent Logo,\[ (L_h)_{z-J\zeta,0}g(x) = \dps \int e^{ix\xi}h(x+(z-J\zeta)-J\xi)\hat{g}(\xi)\dc \xi .\] 

Portanto, para qualquer $g \in \sa$, $(L_h)_{z,\zeta}(g) = (L_h)_{z-J\zeta,0}(g).$
Como $L_h$ \'e cont\'\i nua em $\SA$, temos $(L_h)_{z,\zeta} = (L_h)_{z-J\zeta,0}
. \ai $ \end{remark}

\begin{precorol}\label{2.6} {\rm Dada $F \in \SA $, temos $(L_F)_{z,\zeta} =
(L_F)_{z - J\zeta,0}$, como operadores de $\sa$ em $\SA$. }  \end{precorol} 

\begin{remark} Sejam $U = E_{z,\zeta}$ e $V = E_{z - J\zeta,0}$.
Notemos que, dada $h \in \sa$, $U(h)$ e $V(h)$ pertencem a $\sa$. Logo, por \citeprop{2.4} e \citelem{2.5}, temos \[(L_F)_{z,\zeta}= (L_F)_{z - J\zeta,0}. \,\,\,  \ai \] \end{remark}

\begin{pretheo}\label{2.7} {\rm Seja $T \in \op$ com $TR_g = R_gT$ para toda $g \in \sa$.
Se $T$ \'e suave por transla\c c\~oes, ent\~ao $T$ \'e Heisenberg-suave.} \end{pretheo}

\begin{remark} Pelo corol\'ario \citecorol{2.3}, temos que existe uma seq\"u\^encia
 $(F_k)$ em $\SA$ tal que $\lim_{k \rightarrow \infty}L_{F_k}(g) = T(g)$, 
\hspace{0,1cm} $\forall g \in \sa$. Isto implica que para $(z,\zeta) \in \R^{2n}$
temos $\lim_{k \rightarrow \infty}(L_{F_k})_{z,\zeta}(g) = T_{z,\zeta}(g)$, \hspace{0,1cm} 
$\forall g \in \sa$. Al\'em disso, pelo \citecorol{2.6}, temos que
 $\lim_{k \rightarrow \infty} (L_{F_k})_{z,\zeta}(g) = \lim_{k \rightarrow 
\infty}(L_{F_k})_{z - J\zeta,0}(g)$, \hspace{0,1cm} $\forall g \in \sa$.
Assim, $T_{z,\zeta}(g) = T_{z - J\zeta,0}(g)$, \hspace{0,1cm} $\forall g \in \sa$ e, como $T$ \'e cont\'\i nua, temos $T_{z,\zeta} = T_{z - J\zeta,0}$.

Como $T$ \'e suave por transla\c c\~oes, a aplica\c c\~ao 
\[ (z,\zeta) \longrightarrow z - J\zeta \longrightarrow T_{z - J\zeta,0} \]
\'e $C^\infty$,  ou seja $(z,\zeta) \longrightarrow T_{z,\zeta}$ \'e $C^\infty. \ai$ \end{remark}

\chapter{ }


\bigskip
\hspace{0,8cm} Neste cap\'\i tulo veremos uma generaliza\c
c\~ao do teorema de
Calder\'on-Vaillancourt
\cite{10} para operadores
pseudo-diferenciais com s\'\i mbolos tomando valores numa
$C^*$-\'algebra separ\'avel, baseada na demonstra\c
c\~ao
feita por Seiler em \cite{11}. 

\bigskip Dadas uma
$C^*$-\'algebra separ\'avel $A$ e
uma
fun\c c\~ao $a$ em $CB^\infty(\R^{2n},A)$,
veremos que  o operador $a(x,D)$
dado por 

\begin{center} (I) \hspace{1cm} $a(x,D)u(x) = \dps \int
e^{i(x-y)\xi}a(x,\xi)u(y)\dc y \dc \xi$ para $u$ em $\sa$ \end{center}

\noindent \'e cont\'\i nuo com a norma
$\| \, \|_2$; $\|u\|_2 = \| \int
u(x)^*u(x)dx\|^\frac{1}{2}$.

\medskip
A integral em (I) deve ser entendida como uma integral  oscilat\'oria. Mais precisamente,
para cada $x\in\R^n$, $a(x,D)u(x) = \dps\lim_{k,m\rightarrow\infty}\dps\int\dps\int e^{i(x-y)\xi}
a(x,\xi)u(y)\psi_m(y)\psi'_k(\xi)\dc y\dc\xi$, onde $\psi_m, \psi'_k$ s\~ao as seq\"u\^encias introduzidas na segunda p\'agina do cap\'\i tulo 1. Demonstra-se, ent\~ao, que a 
integral oscilat\'oria \'e igual \`a integral iterada (n\~ao necessariamente
absolutamente convergente) $\dps\int e^{ix\xi}a(x,\xi)\left[\dps\int e^{-iy\xi}u(y)\dc y\right]\dc\xi.$

\begin{preobs}\label{3ypoco}{\rm  Notemos que o operador dado por
\[a(x,D)_Ru(x) = \dps\int e^{i(x-y)\xi}u(y)a(x,\xi)\dc y\dc\xi, \hspace{0,5cm} u\in\sa,\]
n\~ao 
\'e um homomorfismo de m\'odulos \`a direita em $\SA$, se $A$ n\~ao \'e comutativa, 
pois,
dados $b\in A, u\in\sa$, $a(x,D)_R(ub)(x)=\int e^{i(x-y)\xi}u(y)ba(x,\xi)\dc y\dc\xi$
e $a(x,D)_R (u)$ $b(x)=\int e^{i(x-y)\xi}u(y)a(x,\xi)b\dc y\dc\xi.$ Ou seja, em geral,
$a(x,D)_R(ub)$ e $a(x,D)_R(u)b$ n\~ao s\~ao iguais. }\end{preobs}

\begin{predef}\label{3umquarto}  {\rm
Definimos a {\it transformada de
Fourier} de $u$ em $\sa$ como segue: \[
\hat{u}(\xi) = \dps \int
e^{-iy\xi}u(y)\dc y .\] } \end{predef}

\begin{preobs}\label{3ummeio} {\rm Temos
(como no caso escalar) que a
transformada de
Fourier \'e um operador cont\'\i nuo e
invers\'\i vel em $\sa$. Al\'em
disso, \'e uma isometria
com rela\c c\~ao \`a norma $\|\, \|_2$, (corol\'ario \ref{B.4}).
Assim, $a(x,D)u(x) =\dps \int
e^{ix\xi}a(x,\xi)\hat{u}(\xi) \dc \xi $. }\end{preobs}

\begin{preobs}\label{3um3/4} {\rm Notemos que, pelo Teorema da Converg\^encia Dominada (\ref{A.7}), $a(x,D)$ $u$
\'e uma fun\c c\~ao cont\'\i nua, pois $\|e^{ix\xi}a(x,\xi)\hat{u}(\xi)\|
\leq \|a\|_\infty\|\hat{u}(\xi)\| \in L^1(\R^n)$.}
\end{preobs}

\begin{predef}\label{3.1} {\rm Dada $a \in
CB^{\infty}(\R^{2n},A)$, seja
$\pi(a) = \supp \{\|\partial^\beta_y
\partial^\gamma_\eta a \|_\infty,
\beta,\gamma
\leq (1,1,\cdots ,1) \}. $} \end{predef}

\begin{prenotac}\label{3.tresmeio} $ \alpha
= (1,1, \cdots ,1) \in \N^n$.
\end{prenotac}

\begin{prelemma}\label{3.2} {\rm Seja $\phi
\in C^\infty_c(\R^{2n})$ com
$\phi \equiv 1$ perto do zero. Dada $a \in
CB^\infty(\R^{2n},A)$, seja
$a_\varepsilon (y,\eta) =
\phi(\varepsilon y,\varepsilon
\eta)a(y,\eta)$, $0 < \varepsilon \leq 1$.
Ent\~ao existe uma constante $c_1 \in \R^+$,
independente de $a$ e
$\varepsilon$, tal que $\pi(a_\varepsilon)
\leq c_1\pi(a)$.}\end{prelemma}

\begin{remark} Dados os multi-\'\i ndices
$\beta, \gamma \leq \alpha$,
\[ \partial^\beta_y\partial^\gamma_\eta
a_\varepsilon(y,\eta) =
\partial^\beta_y[\partial^\gamma_\eta
\phi(\varepsilon y,\varepsilon
\eta)a(y,\eta)] =
 \partial^\beta_y \left[ \dps
\sum_{\gamma_1+\gamma_2=\gamma}\partial^{\gamma_1}_\eta
\phi(\varepsilon
y,\varepsilon \eta) \partial^{\gamma_2}_\eta
a(y,\eta) \right]= \]
\[ = \dps \sum_{\beta_1+\beta_2=\beta} \dps
\sum_{\gamma_1+\gamma_2=\gamma}
\partial^{\beta_1}_y\partial^{\gamma_1}_\eta
\phi(\varepsilon y,\varepsilon
\eta)
\partial^{\beta_2}_y\partial^{\gamma_2}_\eta
a(y,\eta). \]

\medskip
Portanto, se $\phi_\epsi(x,y) = \phi(\epsi y,\epsi x)$, temos: 
\[ \|\partial^\beta_y\partial^\gamma_\eta a_\epsi\|_\infty  \leq  \dps\sum_{\beta_1\leq\beta}
\dps\sum_{\gamma_1\leq\gamma} \|\partial^{\beta_1}_y\partial^{\gamma_1}_\eta 
\phi_\epsi \|_\infty \pi(a). \]

Assim, \[ \pi(a_\varepsilon) \leq \dps
\sum_{\beta,\gamma \leq
\alpha}\varepsilon^{|\beta
+\gamma|}\|\partial^\beta_y\partial^\gamma_\eta
\phi\|_\infty \pi(a) \leq \dps
\sum_{\beta,\gamma \leq
\alpha}\|\partial^\beta_y\partial^\gamma_\eta
\phi \|_\infty \pi(a).
\hspace{0,3cm} \ai \]\end{remark}

\begin{prenota}\label{3.2umquarto} {\rm
Dadas $f,g \in L^2 (\R^n,A)$ (ver \ref{A.8}),
$(f|g) = \dps \int f(x)^*g(x) dx$ est\'a bem
definido, pois $f^*g \in L^1(\R^n,A)$ (ver \ref{A.8}) por H\"older \cite[6.2]{3}.
E temos que $(\cdot|\cdot):L^2(\R^n,A)\times L^2(\R^n,A)
\rightarrow A$ \'e cont\'\i nua como fun\c c\~ao
das duas vari\'aveis, pois $\|(f|g)\|\leq\|f\|_{L^2}\|g\|_{L^2}.$

\noindent Notemos que, se $h_1,h_2\in\sa, \,\, (h_1|h_2)=<h_1,h_2>$. E
temos que $\|h\|_2 \leq \|h\|_{L_2}$, para $h\in\sa$.}
%
%
%
\end{prenota}

\begin{preprop}\label{3.3} {\rm $L^2(\R^n,A)
\subseteq \SA $ e a inclus\~ao \'e cont\'\i nua.}\end{preprop}

\begin{remark} Definamos $\iota :(\sa ,\|\cdot\|_{L_2}) \rightarrow \SA$ por
  $\iota (h)=h,\, h \in \sa$. Notemos que
  $\|\iota (h)\|_{2} = \|h\|_{2} \leq \|h\|_{L_2}$. Portanto, $\iota$ \'e cont\'inua.

Seja ${\cal I}:L^2(\R^n,A) \rightarrow \SA$ a \'unica extens\~ao
  cont\'\i nua de $\iota$. Provemos que ${\cal I}$ \'e injetora; ou seja que se
  ${\cal I}(f)=0 \,\Rightarrow\, f=0$, para $f \in L^2(\R^n,A)$.

Dada $f \in L^2(\R^n,A)$, seja $(f_m)$ uma seq\"u\^encia em $\sa$
  tal que $\lim_{m \rightarrow \infty}\|f_m-f \|_{L_2} = 0$ (ver
  \ref{A.22}). Para $g \in L^2(\R^n,A)$ temos que (por H\"older): $\|(f_m-
  f|g)\| \leq \dps\int\|(f_m-f)(x)^*g(x)\|dx \leq
  \|f_m-f\|_{L_2}\|g\|_{L_2}$. Logo,

\[ (1) \hspace{1cm} (f_m|g)\,\,  \rightarrow
\,\, (f|g), \hspace{0,3cm}
\forall g \in L^2(\R^n,A). \]

Se ${\cal I}(f) = 0$, como ${\cal I}$ \`e cont\`inua, ${\cal I}(f_m)
\rightarrow 0$, ou seja, $\|f_m\|_2 \rightarrow 0$. Assim, dada $h \in \sa$,
$\|<f_m ,h>\| \leq  \|f_m\|_2\|h\|_2 \rightarrow 0$ por Cauchy-Schwartz
\cite[lemma 1.1.2]{6}. Ou seja,

\begin{center} (2) \hspace{1cm}  ${\cal
I}(f) = 0 \, \Rightarrow \,\,
<f_m,h> \rightarrow 0, \hspace{0,3cm}
\forall h \in \sa $.    \end{center}

Suponhamos ${\cal I}(f) = 0$. Ent\~ao, dada
$h \in \sa$, por (1) e (2) temos
$(f|h) = 0$;
em particular, para todo $m \in \N$,
$(f|f_m) = 0.$ Novamente por (1),
tomando $g = f$,
temos que $(f|f_m) \rightarrow (f|f)$ e,
portanto, $(f|f) = 0$, ou seja, $f
= 0$,      (\ref{A.15})     e ${\cal{I}}$ \'e
injetora.$\ai$ 
\end{remark}

\begin{precorol}\label{3.3umquinto}{\rm   Se $f, g\in L^2(\R^n,A)$,
ent\~ao $<f,g> = (f|g)$. Em particular, $\|f\|_2=\left\|\dps\int f(x)^*f(x)dx\right\|^\frac{1}{2}.$                                   }\end{precorol}

\begin{remark}   Sejam $(f_m)$ e $(g_k)$ seq\"u\^encias em $\sa$
que convergem em $L^2(\R^n,A)$ para $f$ e $g$, respectivamente.
Como a inclus\~ao $L^2(\R^n,A) \subseteq \SA$ \'e cont\'\i nua, (pela proposi\c c\~ao \ref{3.3}),
$(f_m)$ e $(g_k)$ convergem
em $\SA$ para $f$ e $g$, respectivamente. Portanto, (ver \ref{3.2umquarto}),
\[ (f|g)  \longleftarrow (f_m|g_k) = <f_m,g_k> \longrightarrow <f,g>.\] Isto \'e
$<f,g> = \dps\int f(x)^*g(x)dx.$       $\ai$              \end{remark}

\begin{preobs}\label{3.3umquarto} Se {\rm
$D_{y_j} = -i\partial_{ y_j}$, vale
que $(i+x)^\alpha e^{ixy} =
(i+D_y)^\alpha e^{ixy}$, \hspace{0,3cm} $x,y
\in \R^n$. 

\noindent (Escrevemos $(i+x)$ querendo dizer $(i\alpha +x)$).}
\end{preobs}

\begin{proposition}\label{3.4} Dadas $a \in
CB^\infty (\R^{2n},A)$ e $u \in
\sa$, ent\~ao
$a(x,D)(u) \in L^2(\R^n,A)$. \end{proposition}

\begin{remark} Notemos que, pela observa\c c\~ao \ref{3um3/4},
$a(x,D)(u)(x)$ \'e uma fun\c c\~ao cont\'\i nua de $\R^n$ em $A$ e,
portanto, mensur\'avel. 

\[a(x,D)(u)(x) = \dps \int
e^{ix\eta}a(x,\eta)\hat{u}(\eta)\dc
\eta = \dps \int
(i+x)^{-\alpha}[(i+x)^\alpha
e^{ix\eta}]a(x,\eta)\hat{u}(\eta)\dc \eta =
\]
\[ = (i+x)^{-\alpha} \dps
\int[(i+D_\eta)^\alpha
e^{ix\eta}]a(x,\eta)\hat{u}(\eta) \dc \eta =
(i+x)^{-\alpha} \dps \int
e^{ix\eta} [(i-D_\eta)^\alpha
a(x,\eta)\hat{u}(\eta)]\dc \eta \]
(para a \'ultima igualdade fizemos integra\c
c\~ao por partes, ver \ref{A.14}).

Como $a \in CB^\infty(\R^{2n},A)$ e $\hat{u}
\in \sa$, \citeobs{3ummeio}, a
\'ultima integral \'e limitada; digamos que,
em norma, \'e menor ou igual
que $c_2 \in \R^+$. Veremos na observa\c c\~ao \citeobs{3.5} que
$c_2$ depende de $a$ e de $u$, mas n\~ao depende de $x$. Assim, temos
$\|a(x,D)(u)\|^2_{L_2} \leq c^2_2 \dps \int
|i+x|^{-2|\alpha|}dx$ e portanto
$a(x,D)(u) \in L^2(\R^n,A)$. $\ai$
\end{remark}

\begin{precorol}\label{3.8ummeio} {\rm Pela Proposi\c c\~ao \citeprop{3.3}
$a(x,D)(u) \in \SA$. }\end{precorol}

\begin{preobs}\label{3.5} {\rm Sobre a
constante $c_2$ da proposi\c c\~ao
\citeprop{3.4}.
\[ \left\| \dps \int e^{ix\eta}[(i-D_\eta)^\alpha 
a(x,\eta)\hat{u}(\eta)]\dc \eta \right\| \leq
\dps \int \left\| (i-D_\eta)^\alpha
a(x,\eta)\hat{u}(\eta) \right\|\dc \eta \] e temos
que  $(i-D_\eta)^\alpha
a(x,\eta)\hat{u}(\eta)$ \'e soma
de termos do tipo $D_\eta^{\alpha_1}
a(x,\eta)D^{\alpha_2}_\eta\hat{u}(\eta)$
(multiplicados
por uma pot\^encia de $i$) e vale

\[ \|D^{\alpha_1}_\eta
a(x,\eta)D_\eta^{\alpha_2}\hat{u}(\eta) \|
\leq
\pi(a) \|D^{\alpha_2}_\eta \hat{u}(\eta)
\|,\]

\noindent onde $\alpha_1,\alpha_2 \leq
\alpha$. Assim, $\|\dps \int
e^{ix\eta}[(i-D_\eta)^\alpha
a(x,\eta)\hat{u}(\eta)]\dc \eta\| \leq
M\pi(a)$, onde $M \in \R^+$ \'e uma
constante (que depende de $n$ e de $u$) que tem a ver com a soma citada
acima (n\'umero de termos e
m\'odulo das constantes) e com a soma das
integrais do tipo $\dps \int
\|D^{\alpha_2}_\eta
\hat{u}(\eta)\|\dc \eta$. Logo, a constante
$c_2$ mencionada acima \'e
$M\pi(a)$, onde $M$ n\~ao depende de $a$ nem de $x$.}
\end{preobs}

\bigskip
A seguinte proposi\c c\~ao foi baseada no lema
3.13 de \cite{11}.

\begin{proposition}\label{3.6} Dados $a$ e
$a_\epsi$ como no lema
  \citelema{3.2} e $u, v \in \sa$, vale que
\[<v,a_\varepsilon(x,D)u> \,\,
\underset{\varepsilon \rightarrow 0}
\longrightarrow \,\,
<v,a(x,D)u>
  \textrm{ em } A.\] \end{proposition}
\begin{remark} Sejam $b_\epsi = a_\epsi - a
\hspace{0,3cm} (b_\epsi \in
CB^\infty(\R^{2n},A))$
e $B_\epsi = b_\epsi(x,D)$. Provemos que
$<v,B_\epsi u> \,\, 
\underset{\epsi \rightarrow 0}
\longrightarrow \,\,
0$
$ B_\epsi u(x) = \dps \int
e^{ix\eta}b_\epsi(x,\eta)\hat{u}(\eta)\dc
\eta$, e
 temos que
$\|e^{ix\eta}b_\epsi(x,\eta)\hat{u}(\eta)\|
\leq
 \pi(b_\epsi)\|\hat{u}(\eta)\|$.

 Observemos que, pelo lema \citelema{3.2},
existe $c_1 \in \R^+$,
independente de $b_\epsi$, tal que
$\pi(b_\epsi) = \pi(a_\epsi - a) \leq
\pi(a_\epsi) + \pi(a) \leq (c_1 + 1)\pi(a).$
Portanto,
$\|e^{ix\eta}b_\epsi(x,\eta)\hat{u}(\eta)\|
\leq (c_1 +
1)\pi(a)\|\hat{u}(\eta)\|$.

Como
$\|e^{ix\eta}b_\epsi(x,\eta)\hat{u}(\eta)\|
\rightarrow 0$ e $\dps
\int\|\hat{u}(\eta)\|\dc \eta < \infty$,
pelo Teorema da Converg\^encia
Dominada, $\|B_\epsi u(x)\| \rightarrow 0$,
para cada
$x \in \R^n$.

Por outro lado, como visto na proposi\c c\~ao \citeprop{3.4} e na observa\c c\~ao \citeobs{3.5}, existe uma constante $M
\in \R^+$, independente de
$\epsi$, tal que
\[\dps \int \|B_\epsi u(x)\|^2dx \leq
M^2\pi(b_\epsi)^2 \dps \int
|i+x|^{-2|\alpha|}dx \leq
M^2(c_1 +
1)^2\pi(a)^2\dps\int|i+x|^{-2|\alpha|}dx <
\infty . \] Novamente,
pelo Teorema da Converg\^encia Dominada,
$\|B_\epsi u\|_{L_2} \rightarrow 0$.
Isto implica que
$\|B_\epsi u\|_2 \rightarrow 0$ por
\citenota{3.2umquarto}. Finalmente, por
Cauchy-Schwarz,
como $\|<v,B_\epsi u>\| \leq
\|v\|_2\|B_\epsi u\|_2$, temos que
$<v,a_\epsi(x,D)u> \rightarrow
<v,a(x,D)u>$ em $A$. $\ai$ \end{remark}

\begin{prenota}\label{3.6meio} {\rm Dada $a
\in CB^\infty(\R^{2n},A)$, pelo
    visto na proposi\c c\~ao \citeprop{3.4},
$a(x,D)(u) \in L^2(\R^n,A)$, para $u \in \sa$.
Assim, pela nota
\citenota{3.2umquarto}, $<v,a(x,D)u>$
est\'a bem definido para $v \in \sa$. }
\end{prenota}

\begin{prelemma}\label{3.6.1/2}{\rm   Dada $u \in \sa$ seja $f(x,\xi) = 
\dps\int e^{-iy\xi}(i+x-y)^{-\alpha}u(y)dy$. Ent\~ao, $f \in L^2(\R^{2n},A).$  } 
\end{prelemma}

\begin{remark}  Como para qualquer $N \in \N$

\[ \frac{(1 - \Delta_y)^Ne^{-iy\xi}}{(1 + |\xi|^2)^N} = e^{-iy\xi} ,\]    temos

\[ f(x,\xi) = \frac{1}{(1 + |\xi|^2)^N}\dps\int \left[ (1 - \Delta_y)^Ne^{-iy\xi} \right]
(i+x-y)^{-\alpha}u(y)dy .\]
Queremos provar que a integral acima \'e combina\c c\~ao linear de 
termos do tipo

\[ \dps\int e^{-iy\xi}\partial^\beta_yu(y)(i+x-y)^{-\gamma}dy, \hspace{0,5cm} 
\beta \geq 0 \,\, e \,\, \gamma \geq \alpha . \] 

Para isto, calculemos, para $j = 1, \cdots , n, \,\, M, P \leq N$:

\[ \dps\int \partial_{y_j} \left[ (1 - \Delta_y)^M e^{-iy\xi} \right] \left[ 
(1 - \Delta_y)^P u(y)(i+x-y)^{-\alpha} \right] dy = (*) . \]

Para $(x,\xi)$ fixo, o integrando acima \'e uma fun\c c\~ao de 
$L^1(\R^n,A)$, pois $u \in \sa$. Sejam $g(y) = (1 - \Delta_y)^Me^{-iy\xi}$ e
$h(y) = (1 - \Delta_y)^Pu(y)(i+x-y)^{-\alpha}$. 
Pelo corol\'ario \ref{A14ummeio}, temos que 
\[\dps\int\left[\partial_{y_1}g(y)\right](y)dy = - \dps\int g(y)(\partial_{y_1}h(y))dy.\]
Pelo Teorema de Fubini (\ref{A.12}), temos,
tomando $j = 1$ sem perda de generalidade:

\[ (*) = \dps\int_{\R^{n-1}} \left( \dps \int^{+\infty}_{-\infty} \partial_{y_1}g(y)h(y)dy_1 \right)
dy_2 \cdots dy_n = \] 
\[ =   \dps\int_{\R^{n-1}}\left[ \dps\lim_{R \rightarrow +\infty}
g(R, y_2, \cdots , y_n)h(R, y_2, \cdots , y_n) - \right. \] \[ \left. 
\dps\lim_{R \rightarrow -\infty} g(R, y_2, \cdots , y_n)h(R, y_2, \cdots , y_n)\right]
dy_2, \cdots , dy_n   ,       \] pelo Teorema Fundamental do C\'alculo, \ref{A.13}.

Cada um dos limites acima \'e igual a zero.
Portanto,

\[ \dps\int[\partial_{ y_1}g(y)]h(y)dy = - \dps\int g(y)[\partial_{ y_1}h(y)]dy , \] isto \'e,
podemos fazer
a integra\c c\~ao por partes. Deste modo, temos que $f(x,\xi)$ \'e
combina\c c\~ao linear de termos do tipo

\[ \frac{1}{(1 + |\xi|^2)^N}\dps\int e^{-iy\xi}\partial^\beta_yu(y)(i+x-y)^{-\gamma}dy, \hspace{0.5cm}  
\beta \geq 0, \gamma \geq \alpha, \] onde a integral \'e finita pois $u \in \sa$.

Por outro lado, pela desigualdade de Petree (\cite[(3.6)]{15}), temos que

\[ |(i+x-y)^{-\gamma}| \leq 4^n\frac{\sqrt{1 + y_1^2} \cdots \sqrt{1+y_n^2}}
{\sqrt{1+x_1^2} \cdots \sqrt{1+x_n^2}}. \]

E, al\'em disso,

\[ \partial_{ y_1}(i+x-y)^{-\gamma} = -\gamma_1(i+x_1-y_1)^{-\gamma_1-1}(i+x_2-y_2)^{-\gamma_2}
\cdots (i+x_n-y_n)^{-\gamma_n}. \]

Deste modo, $\|f(x,\xi)\|$ \'e menor ou igual que uma combina\c c\~ao linear
de termos do tipo

\[ \frac{4^n}{   (1 + |\xi|^2)
^N\sqrt{1+x_1^2} \cdots \sqrt{1+x_n^2}}\dps\int\|\partial^\beta_yu(y)\|
\sqrt{1+y_1^2} \cdots \sqrt{1+y_n^2}dy, \hspace{0.5cm} \beta \geq 0, \] que
pertence a $L^2(\R^{2n})$ para $N \geq \frac{n}{2}$. Ent\~ao $f \in L^2(\R^{2n},A)$.
$\ai$  \end{remark}

\begin{prelemma}\label{3.7}{\rm Dada $u \in
\sa$, seja $f(x,\xi) = \dps
\int
e^{-iy\xi}(i+x-y)^{-\alpha}u(y)\dc
y.$ Ent\~ao  existe $c_3 \in \R^+$,
independente de $f$ e de $u$, tal que
$\|f\|_2 \leq
c_3\|u\|_2$.}\end{prelemma}

\begin{remark}  Sejam $h(z) =
(i-z)^{-\alpha}$ e $ h_z(y) = h(y-z)$;
assim,
$f(x,\xi) = \dps \int
e^{-iy\xi}h_x(y)u(y)\dc y =
\widehat{h_xu}(\xi).$ 
(Notemos que $h_xu \in \sa$.)
Como
a transformada de Fourier
\'e um operador "unit\'ario" em $\sa$
(\ref{B.3}), temos que, dados $ g, p
\in \sa$, $<\hat{p},\hat{g}> = <p,g>$.

\medskip \noindent Vimos no lema \citelem{3.6.1/2} 
que $f \in L^2(\R^{2n},A)$, logo,
\noindent pelo Teorema de Fubini,
(\ref{A.12}), temos

\[\dps \int f(x,\xi)^* f(x,\xi)dx d\xi
=\dps \int [\widehat{h_xu}(\xi)]^*
\widehat{h_xu}(\xi)d \xi dx  = \dps \int
<\widehat{h_xu},\widehat{h_xu}>dx
=  \] \[ = \dps \int <h_xu,h_xu>dx  =  \dps \int u(\xi)^*\overline{h_x(\xi)}
h_x(\xi)u(\xi)d\xi dx = \dps \int
u(\xi)^* \dps\int
\overline{h_x(\xi)}h_x(\xi)dx u(\xi)d\xi =
\]

\[= \dps\int u(\xi)^*\dps\int
\overline{h(\xi -x)}h(\xi -x) dx u(\xi) d\xi
=
\dps\int \overline{h(x)}h(x) dx \dps\int
u(\xi)^*u(\xi)d\xi .\] Logo, se $c_3
= \|h\|_{L_2}$, (\ref{3.3}), temos $\|f\|_2 \leq c_3
\|u\|_2$. $\ai$
\end{remark}

\bigskip
O pr\'oximo teorema foi enunciado e demonstrado 
baseando-nos no teorema 3.14 de \cite{11}.

\begin{pretheo}\label{3.8} {\rm Dada $a \in
CB^\infty(\R^{2n},A)$, ent\~ao
$a$ induz um operador cont\'\i nuo em $\sa$
como pr\'e-m\'odulo de Hilbert.
Al\'em disso, existe uma constante $k \in
\R^+$, independente de $a$, tal
que $\|a(x,D)\| \leq k\pi(a)$.}\end{pretheo}

\begin{remark} Suponhamos, primeiramente,
que $a$ tem suporte compacto. Seja
$T= a(x,D)$. Calculemos, para $u,v \in \sa$,
$<\hat{v},\widehat{Tu}> =
<v,Tu>.$

\[ Tu(x) = \dps \int
e^{i\xi(x-y)}a(x,\xi)u(y)\dc y\dc \xi \]

\[ \widehat{Tu}(\eta) = \dps\int
e^{-ix\eta}e^{i\xi(x-y)}a(x,\xi)u(y)\dc
y\dc \xi\dc x \]

\[ <\hat{v},\widehat{Tu}> =
\dps \int
e^{-ix\eta}e^{i\xi (x-y)}\hat{v}(\eta)^*a(x,\xi)u(y)\dc
y\dc
\xi\dc x d\eta = \]

\[ = \dps\int
e^{-ix\eta}(i+x-y)^{-\alpha}\hat{v}(\eta)^*[(i+x-y)^\alpha e^{i\xi(x-y)}]a(x,\xi)u(y)\dc
\xi\dc x\dc y d\eta.\]
Como $a$ tem suporte compacto e $u,v \in
\sa$, vale Fubini (\ref{A.12}), o
que justifica as igualdades acima; al\'em
disso, podemos aplicar a integra\c
c\~ao por partes e, pela
observa\c c\~ao \citeobs{3.3umquarto},
temos:

\[ < \hat{v},\widehat{Tu} > = \dps\int
e^{-ix\eta}(i+x-y)^{-\alpha}\hat{v}(\eta)^*
[(i+D_\xi)^{\alpha}e^{i\xi(x-y)}]a(x,\xi)u(y)\dc
\xi\dc x\dc y d\eta =  \]

\[ = \dps\int
e^{-ix\eta}(i+x-y)^{-\alpha}v(\eta)^*e^{i\xi
x}e^{-i\xi
y}[(i-D_\xi)^\alpha a(x,\xi)]u(y)\dc \xi\dc
x\dc y d\eta = \]

\[ = \dps\int
e^{-iy\xi}(i+x-y)^{-\alpha}(i+\xi
-\eta)^{-\alpha}[(i+\xi
-\eta)^\alpha e^{ix(\xi
-\eta)}]\hat{v}(\eta)^*[(i-D_\xi)^\alpha
a(x,\xi)]u(y)\dc \xi \dc x\dc y d\eta =\]

\[ = \dps\int
e^{-iy\xi}(i+x-y)^{-\alpha}(i+\xi
-\eta)^{-\alpha}[(i+D_x)^\alpha e^{ix(\xi
  -\eta)}]\hat{v}(\eta)^*[(i-D_\xi)^\alpha
a(x,\xi)]u(y)\dc x\dc\xi\dc y
d\eta
= \]

\[ = \dps\int e^{-iy\xi}(i+x-y)^{-\alpha}
e^{-ix\eta}(i+\xi -\eta)^{-\alpha}
e^{ix\xi}\hat{v}(\eta)^*[(i-D_{x})^{\alpha}(i-D_{\xi})^{\alpha}a(x,\xi)]u(y)\dc
x\dc\xi\dc y d\eta
= \]

\[ = \dps\int
e^{ix\xi}(e^{-ix\eta}h_\xi(\eta)\hat{v}(\eta)^*[(i-D_x)^{\alpha}(i-D_\xi)^\alpha
a(x,\xi)]
(e^{-iy\xi}h_x(y)u(y))\dc y\dc\eta\dc xd\xi
\] usando a nota\c c\~ao do lema
\citelemma{3.7}.

\medskip
Sejam \[f(x,\xi) = \dps\int
e^{-ix\eta}h_\xi(\eta)\hat{v}(\eta)^*\dc\eta
\]
\[ g(x,\xi) = \dps\int
e^{-iy\xi}h_x(y)u(y)\dc y . \]

Ent\~ao: \[ <\hat{v},\widehat{Tu}> =
\dps\int
e^{ix\xi}f(x,\xi)[(i-D_x)^\alpha(i-D_\xi)^\alpha
a(x,\xi)]g(x,\xi)\dc x d\xi
.\]

\medskip
Abusando da nota\c c\~ao, colocando as
vari\'aveis das fun\c c\~oes, temos:

\[ \|<\hat{v},\widehat{Tu}>\|=
\frac{1}{(2\pi)^{\frac{n}{2}}}
\|<e^{-ix\xi}f^*(x,\xi),[(i-D_x)^\alpha
(i-D_\xi)^\alpha a(x,\xi)]g(x,\xi)
>\|. \]

\medskip Seja $c(x,\xi) =
(i-D_x)^\alpha(i-D_\xi)^\alpha a(x,\xi)$.
Como $a
\in CB^\infty(\R^{2n},A)$, temos que
$\|c\|_\infty = \supp_{(x,\xi) \in
\R^{2n}}\|c(x,\xi)\| <\infty .$
Al\'em disso, pela defini\c c\~ao de $c$,
existe uma constante $l \in \R^+$,
que n\~ao depende de $a$, tal que
$\|c\|_\infty \leq l\pi(a).$ E temos, pela demonstra\c c\~ao do teorema 2.2.5
de \cite{5},

\[ \dps\int
g(x,\xi)^*c(x,\xi)^*c(x,\xi)g(x,\xi) dx d\xi
\leq
\dps\int\|c(x,\xi)\|^2g(x,\xi)^*g(x,\xi)
dx d\xi \leq \]

\[ \leq \| c \|^2_\infty \dps\int
g(x,\xi)^*g(x,\xi) dx d\xi . \]

Portanto, pelo lema \citelemma{3.7}, $\| cg
\|_2 \leq l\pi(a)c_3\|u\|_2 .$

\medskip
Como $f^*(x,\xi) = \dps\int
e^{ix\eta}\overline{h_\xi(\eta)}\hat{v}(\eta)\dc\eta$,
de um modo
an\'alogo ao feito na demonstra\c c\~ao do
lema \citelemma{3.7}, temos
$\|f^*\|_2 \leq c_3\|\hat{v}\|_2 =
c_3\|v\|_2$ (ver observa\c c\~ao
\citeobs{3ummeio}). Logo, por
Cauchy-Schwarz,
 \[ \|< \hat{v},\widehat{Tu} >\| \leq
\frac{lc^2_3}{(2\pi)^{\frac{n}{2}}}
\pi(a) \| u \|_2\| v \|_2. \]

Tomando $k_1 = \dps\frac{lc^2_3}{(2\pi)^{\frac{n}{2}}}$, temos
que
$\|<v,a(x,D)u>\|\leq
k_1\pi(a)\|u\|_2\|v\|_2$, onde  $k_1$ n\~ao
depende de
$a$, e ent\~ao $a(x,D)$ \'e
cont\'\i nuo em $\sa$.

\bigskip
Para o caso geral, consideremos $a_\epsi$
como na proposi\c c\~ao
\citeprop{3.6}.
Ent\~ao, pelo feito acima, temos:

\[ \| <v,a(x,D)u> \|
\underset{\varepsilon \rightarrow 0}
\longleftarrow
\|<v,a_\epsi(x,D)u>\| \leq
k_1\pi(a_\epsi)\|u\|_2\|v\|_2.\]

Portanto, pelo lema \citelemma{3.2},
$\|<v,a(x,D)u>\| \leq
k_1c_1\pi(a)\|u\|_2\|v\|_2.$
E, assim, existe $k \in \R^+, k = k_1c_1$,
independente de $a$, tal que
$\|<v,a(x,D)u>\| \leq k\pi(a)\|u\|_2\|v\|_2$.

\medskip
Tomemos $v = a(x,D)(u)$. Como $v \in 
L^2(\R^n,A)$, proposi\c c\~ao \citeprop{3.4}, 
pela nota \citenota{3.2umquarto}, a seguinte 
desigualdade faz sentido:

\[ \|<a(x,D)(u),a(x,D)(u)>\| \leq  k\pi(a)\|u\|_2
\|a(x,D)(u)\|_2, \,\, \forall u \in \sa .\]

Portanto,  $\|a(x,D)(u)\|_2 \leq k\pi(a)\|u\|_2$,
e temos $\|a(x,D)\| \leq k\pi(a)$.$\ai$
\end{remark}

\begin{precorol}\label{3.9} {\rm Nas condi\c
c\~oes do teorema, vemos que
$a(x,D)$ se estende ao m\'odulo de Hilbert
$\SA$, continuamente. }
\end{precorol}

\chapter{}

\vspace{1cm}
No cap\'\i tulo 1, vimos a demonstra\c c\~ao
da conjectura proposta por Rieffel no final
do 
cap\'\i tulo 4 de \cite{1} para o caso da
$C^*$-\'algebra dos n\'umeros complexos. Foi
necess\'aria a caracteriza\c c\~ao feita por
Cordes (\cite[cap\'\i tulo 8]{2}) dos
operadores suaves
 pela a\c c\~ao do grupo de Heisenberg, os
operadores pseudo-diferenciais. Neste
cap\'\i tulo,
veremos que uma generaliza\c c\~ao da
caracteriza\c c\~ao de Cordes para
uma $C^*$-
\'algebra, $A$, qualquer, implicaria
que a conjectura proposta por Rieffel \'e verdadeira. 

\bigskip
Suponhamos ent\~ao que vale a seguinte
generaliza\c c\~ao da caracteriza\c c\~ao de Cordes:
\begin{prelema}\label{4.1}  {\rm Seja $\cal{HS}$ o conjunto dos operadores
Heisenberg-suaves em ${\cal B}^* (\overline{S^A(\R^n)})$. Ent\~ao
a aplica\c c\~ao $O:CB^\infty(\R^{2n},A)\rightarrow {\cal{HS}}$,
dada por $O(a) =a(x,D)$ \'e uma bije\c c\~ao, onde $a(x,D)$
\'e o operador definido em (I) do cap\'\i tulo 3. }\end{prelema}

\bigskip 
Vejamos que $O$ (acima) est\'a bem definida.

\begin{prelemma}\label{4.2} {\rm  Dado $a\in CB^\infty(\R^{2n},A)$,
de suporte compacto, o operador $O(a)$ possui adjunto; ou seja,
existe um operador $T$ em $\SA$, cont\'\i nuo, tal que 
\[ <O(a)u,v> = <u,Tv>, \,\,\, se \,\, u,v\in\sa.\]                    }\end{prelemma}

\begin{remark}   \[<O(a)u,v> =  \dps\int \left[ \dps\int e^{i(x-y)\xi}a(x,\xi)u(y)\dc y\dc\xi\right]
^*v(x)dx = \] \[= \dps\int\dps\int e^{i(y-x)\xi}u(y)^*a(x,\xi)^*\dc y\dc\xi v(x)dx = \dps\int u(y)^*
\left[\dps\int e^{i(y-x)\xi}a(x,\xi)^*v(x)\dc x\dc\xi \right]dy,\]   para a
\'ultima igualdade aplicamos o Teorema de Fubini, \ref{A.12}, pois $a$ tem suporte
compacto e $u,v\in\sa$.

\medskip  Como $a$ tem suporte compacto, $a\in S^A(\R^{2n})$.
Seja $F_2$ a transformada de Fourier,
na segunda vari\'avel, definido em $S^A(\R^{2n})$. Com esta nota\c c\~ao, seja

\[ c(y,z) = \dps\int e^{iz\xi}a(y-z,\xi)^*\dc\xi = F^{-1}_2a^*(y-z,z) \hspace{0,5cm} (1).\]     e
temos $c\in S^A(\R^{2n})$.

\medskip  Seja $ p(y,z) = F_2 c(y,z)$, \hspace{0,3cm} ($p\in S^A(\R^{2n})$, ver in\'\i cio do Ap\^endice B), ent\~ao 

\[c(y,z) = F^{-1}_2p(y,z) \hspace{1,5cm} (2) \] Por (1), temos $c(y,y-x) =
\dps\int e^{i(y-x)\xi}a(x,\xi)^*\dc\xi$ e, por (2), $c(y,y-x) = \dps\int e^{i(y-x)\xi}p(y,\xi)\dc\xi.$
Portanto, $\dps\int e^{i(y-x)\xi}a(x,\xi)^*v(x)\dc\xi\dc x = \dps\int e^{i(y-x)\xi}p(y,\xi)v(x)
\dc\xi\dc x.$

\bigskip Assim,
\[\dps\int u(y)^*\left[\dps\int e^{i(y-x)\xi}a(x,\xi)^*v(x)\dc x\dc\xi\right]dy = \dps\int u(y)^*\left[
\dps\int e^{i(y-x)}p(y,\xi)v(x)\dc x\dc\xi\right]dy\] ou seja, $<O(a)u,v> = <u,O(p)v>$,
\hspace{0,3cm} $\forall u,v\in\sa$.  $\ai$
  \end{remark}

\begin{preobs}\label{4.3} {\rm Sobre a fun\c c\~ao $p$  do lema \ref{4.2},
\[p(y,\xi) = F_2c(y,\xi) = \dps\int e^{-i\xi z}c(y,z)\dc z = \dps\int e^{-i\xi z}\dps\int e^{iz\zeta}a(y-z,\zeta)^*\dc\zeta\dc z =\] \[= \dps\int e^{-iz(\xi-\zeta)}a(y-z,\zeta)^*\dc\zeta\dc z  = \dps\int
e^{-iz\eta}a(y-z,\xi-\eta)^*\dc\eta\dc z.\]  }                       \end{preobs}

\begin{preobs}\label{4.3ummeio} {\rm Notemos que, como visto no cap\'\i tulo 1, $
CB^\infty (\R^{2n},A)$ \'e um espa\c co de Fr\'echet, com as semi-normas 

\[\|F\|_j = \dps\sup_{|\alpha|\leq j}\|\partial^\alpha F\|_\infty,\] onde
$F\in CB^\infty (\R^{2n},A), \,\,\, j\in \N, \,\,\, \alpha$ \'e um multi-\'\i ndice e
$\|\partial^\alpha F\|_\infty = \dps\sup_{(x,\xi)\in\R^{2n}}\|\partial^\alpha F(x,\xi)\|$.

Temos tamb\'em que $CB^\infty (\R^{2n},CB^\infty(\R^{2n},A))$ \'e um
espa\c co de Fr\'echet com as semi-normas 

\[\|G\|_{jk} = \dps\sup_{i \leq j}\dps\sup_{|\beta|\leq k}\|\partial^\beta G\|_i,\] $G\in
 CB^\infty(\R^{2n}, CB^\infty(\R^{2n},A))$, onde 
\[ \|G\|_i = \dps\sup_{(z,\eta)\in \R^{2n}}\|G(z,\eta)\|_i;\] lembremos que,
para cada $(z,\eta)\in\R^{2n}$, $G(z,\eta)$ \'e uma 
fun\c c\~ao em $CB^\infty(\R^{2n},A)$.}\end{preobs}

\bigskip
Aplicaremos no pr\'oximo lema os resultados de \cite{1} sobre integrais
oscilat\'orias, para fun\c c\~oes tomando valores no espa\c co de Fr\'echet $CB^\infty(\R^{2n},A)$.

\begin{prelemma}\label{4.4}  {\rm  Dada $a\in CB^\infty(\R^{2n},A)$, se $p(y,\xi) =
\dps\int e^{-iz\eta}a(y-z,\xi-\eta)^*\dc z\dc\eta$, ent\~ao
$p\in CB^\infty(\R^{2n},A)$. }\end{prelemma}

\begin{remark}  Seja $G:\R^{2n}\rightarrow CB^\infty(\R^{2n},A)$ dada por $G(z,\eta)(y,\xi)
=a(y-z,\xi-\eta).$ Provemos que $G\in CB^\infty(\R^{2n},CB^\infty(\R^{2n},A))$.  

\medskip Dado $\epsi>0$, existe $\delta>0$ tal que, se $\|(z,\eta)-(z',\eta')\|<\delta$,
$\left(\Leftrightarrow \|(y-z,\xi-\eta)-(y-\right.$ $\left. z', \xi-\eta')\|<\delta, \,\,\, \forall (y,\xi)\in\R^{2n}\right),$ ent\~ao $\|(\partial^\alpha a)(y-z,\xi-\eta)-(\partial^\alpha a)(y-z',\xi-\eta')\|<\epsi$, pois 
$\partial^\alpha a$ \'e uniformemente cont\'\i nua.

\medskip                                   
Como, para cada $(z,\eta)\in\R^{2n}, \, G(z,\eta)$ \'e uma fun\c c\~ao de vari\'avel
$(y,\xi)\in\R^{2n}$, escrevemos $\partial^\alpha(G(z,\eta))$ para denotar
a derivada de multi-\'\i ndice $\alpha$ da fun\c c\~ao $G(z,\eta)\in CB^\infty(\R^{2n},A)$, para $(z,\eta)$ fixado. E, como
$ \partial^\alpha (G(z,\eta))(y,\xi) = (\partial^\alpha a)(y-z,\xi-\eta),$ 
temos $\|\partial^\alpha (G(z,\eta))-\partial^\alpha (G(z',\eta'))\|_\infty
<\epsi$. Logo, dados $j\in\N$ e $\epsi>0$,
existe $\delta>0$
tal que, se $\|(z,\eta)-(z',\eta')\|<\delta$, $\|G(z,\eta)-G(z',\eta')\|_j<\epsi$,   e $G$ \'e cont\'\i nua.

\medskip Al\'em disso, \'e f\'acil ver que
$\forall j\in\N,\forall(z,\eta)\in\R^{2n}, \|G(z,\eta)\|_j\leq\|a\|_j,$
i.e. $G$ \'e limitada.

\medskip Vejamos agora que $G$ \'e deriv\'avel. Para isto, provemos
que $G$ possui derivadas parciais de primeira ordem.

\medskip Se $e_j = (0,\cdots,1,0,\cdots,0)\in \R^n$, provemos que 
\[\dps\lim_{t\rightarrow 0} \frac{G(z+te_j,\eta)-G(z,\eta)}{t} = H(z,\eta), \hspace{1cm} 
(z,\eta)\in\R^{2n},\] onde 
\[H(z,\eta)(y,\xi) = \dps\lim_{k\rightarrow 0}\frac{a(y-z-ke_j,\xi-\eta)-a(y-z,\xi-\eta)}
{k}.\] Notemos que $H(z,\eta)\in CB^\infty(\R^{2n},A)$.

\medskip Consideraremos $t \geq 0$, o caso $t\leq 0$ \'e an\'alogo.

Dado $(y,\xi)\in\R^{2n}$, pelo Teorema Fundamental do C\'alculo \ref{A.13},
temos:

\[\frac{G(z+te_j,\eta)(y,\xi) - G(z,\eta)(y,\xi)}{t} - H(z,\eta)(y,\xi) = \]
\[ = \frac{a(y-z-te_j,\xi-\eta) - a(y-z,\xi-\eta)}{t}  + \partial_j a(y-z,\xi-\eta) = \]
\[ = \dps\int_{[0,t]}\frac{-\partial_j a(y-z-he_j,\xi-\eta)}{t} + \partial_j a(y-z,\xi-\eta)dh =\]
\[= -\dps\int_{[0,t]}\frac{\partial_j a(y-z-he_j,\xi-\eta) - \partial_j a(y-z,\xi-\eta)}{t}dh = \]
\[ = \dps\int_{[0,t]}\dps\int_{[0,h]}\frac{\partial^2_j a(y-z-se_j,\xi-\eta)}{t}ds dh, \hspace{1cm}
 h\in[0,t].\] Aqui estamos usando a nota\c c\~ao
$\partial_j a(y,\zeta) = \dps\lim_{k\rightarrow 0} \frac{a(y+ke_j,\xi)-a(y,\xi)}{k}$.

\noindent Logo, para todo $(y,\xi)\in\R^{2n}$, 
\[ \left\|\frac{G(z+te_j,\eta)(y,\xi) - G(z,\eta)(y,\xi)}{t} - H(z,\eta)(y,\xi)\right\|  \leq \]
\[ \leq \dps\int_{[0,t]}\dps\int_{[0,h]}\frac{\|a\|_2}{t} ds dh \leq t\|a\|_2,\]

o que implica que, para qualquer $(z,\eta)\in\R^{2n},$
\[ \dps\sup_{(y,\xi)\in\R^{2n}}\left\|\frac{G(z+te_j,\eta)(y,\xi)-G(z,\eta)(y,\xi)}{t} 
- H(z,\eta)(y,\xi)\right\| \leq t\|a\|_2.\]

Aplicando o mesmo resultado para $\partial^\alpha(G(z,\eta)), \,\,\, |\alpha|\leq j$,
a derivada de  ordem $\alpha$ da fun\c c\~ao $G(z,\eta)$, temos:
\[\dps\lim_{t\rightarrow 0}\left\|\frac{ G(z+te_j,\eta)-
G(z,\eta)}{t} - H(z,\eta)\right\|_j \leq \]
\[\leq \dps\lim_{t\rightarrow 0}\dps\sup_{|\alpha|\leq j}\left\|\frac{\partial^\alpha (G(
z+te_j,\eta))-\partial^\alpha (G(z,\eta))}{t} - \partial^\alpha(H(z,\eta))\right\|_\infty \leq \]
\[ \leq \dps\lim_{t\rightarrow 0}\dps\sup_{|\alpha|\leq j}\|\partial^\alpha a\|_2.t = 0, \hspace{1cm}
 \forall j \in \N. \]

Do mesmo modo, vemos que existem todas as derivadas parciais
de $G$ e portanto $G$ \'e $C^\infty$. Como antes, vemos que estas derivadas
s\~ao limitadas e, portanto, $G\in CB^\infty(\R^{2n},CB^\infty(\R^{2n},A))$.

\medskip 
Assim, pela defini\c c\~ao de integral oscilat\'oria (\cite[(1.3)]{1}), 
\[\dps\int e^{iz\eta}
G(z,\eta)\dc z\dc\eta \in CB^\infty(\R^{2n},A)\] e, portanto, $\left(\dps\int e^{iz\eta}G(z,\eta)
\dc z\dc\eta\right)^* \in CB^\infty(\R^{2n},A)$, ou seja \[p(y,\xi) = \dps\int e^{-iz\eta}
a(y-z,\xi-\eta)^*\dc z\dc\eta \in CB^\infty(\R^{2n},A). \hspace{0,3cm} \ai\]
  \end{remark}  

\begin{preprop}\label{4.5} {\rm    Dada $a\in CB^\infty(\R^{2n},A)$, existe $p\in CB^\infty(\R^{2n},A)$ tal que, se $u,v\in\sa$, $<O(a)u,v> = <u,O(p)v>.$                                    }\end{preprop}                                           

\begin{remark} Seja $p(y,\xi) = \dps\int e^{-iz\eta}a(y-z,\xi-\eta)^*\dc z\dc\eta$.
Vimos no lema \ref{4.4} que $p\in CB^\infty(\R^{2n},A)$. A  integral acima \'e uma
integral oscilat\'oria.  

Sejam $\psi_m,\psi'_k \in C^\infty_c(\R^n)$, como no cap\'\i tulo 1
e sejam $a_{m,k}(x,w) = \psi_m(x)\psi'_k(w)a(x,w)$ e $p_{m,k}(y,\xi)
= \dps\int e^{-iz\eta}a_{m,k}(y-z,\xi-\eta)^*\dc z\dc\eta$.
Segue da proposi\c c\~ao 1.6 de \cite{1} que, para cada $(y,\xi)\in\R^{2n}$, $\dps\lim_{m,k}p_{m,k}(y,\xi) = p(y,\xi)$.

\medskip Dadas $u,v\in\sa$, \[<u,O(p)v>= \dps\int u(y)^*\left[\dps\int
e^{iy\xi}p(y,\xi)\hat{v}(\xi)\dc\xi\right]dy = \] \[ = \dps\int u(y)^*\left[\dps\int e^{iy\xi}\dps\lim_{m,k}p_{m,k}(y,\xi)\hat{v}(\xi)\dc\xi\right]dy\] 

\'E f\'acil ver que decorre da desigualdade b\'asica 1.4 de \cite{1} (e da defini\c c\~ao
de $\psi_m$ e $\psi'_k$) que existem $j\in\N$ e uma constante $M\in\R^+$, que
n\~ao depende de $a$ nem de $\psi_m$ ou $\psi'_k$ (\'e poss\'\i vel escolher $\psi_m$
tais que $\dps\supp_{\alpha,m}|\psi^{(\alpha)}_m|<\infty$,
e analogamente para $\psi'_k$), tais que $\|p_{mk}(y,\xi)\|
\leq M\|a\|_{2j}.$

\medskip Assim, $\|u(y)^* p_{m,k}(y,\xi)\hat{v}(\xi)\| \leq
M\|a\|_{2j}\|u(y)\|\|\hat{v}(\xi)\| \in L^1(\R^{2n})$.

\medskip Pelo Teorema de Fubini (\ref{A.14}) temos:
\[\dps\int u(y)^*\left[\dps\int e^{iy\xi}p_{m,k}(y,\xi)\hat{v}(\xi)\dc\xi\right] dy = \dps\int\dps\int
u(y)^*e^{iy\xi}p_{m,k}(y,\xi)\hat{v}(\xi)\dc\xi dy.\] Pelo Teorema
da Converg\^encia Dominada (\ref{A.7}), temos:

\[\dps\lim_{m,k}\dps\int\dps\int u(y)^*e^{iy\xi}p_{m,k}(y,\xi)\hat{v}(\xi)\dc\xi dy 
 = \dps\int \dps\int u(y)^* e^{iy\xi}\dps\lim_{m,k}p_{m,k}(y,\xi)\hat{v}(\xi)\dc\xi dy =\]
\[ = \dps\int u(y)^*\left[\dps\int e^{iy\xi}\dps\lim_{m,k}p_{m,k}(y,\xi)\hat{v}(\xi)\dc\xi\right]dy\] (novamente pelo Teorema de Fubini). 

Portanto, $<u,O(p)v> = \dps\lim_{m,k}<u,O(p_{m,k})v>.$

\medskip Pelo lema \ref{4.2} e pela observa\c c\~ao \ref{4.3},
 $<O(a_{m,k})u,v> = <u,O(p_{m,k})v>.$
Al\'em disso, como $a_{m,k}(x,\xi) = \psi_m(x)
\psi'_k(\xi)a(x,\xi) $, temos $\|a_{m,k}(x,\xi)\|\leq\|a(x,\xi)||, \,\, \forall m,\forall k$.
Portanto, $\|\hat{u}(\xi)^* a(x,\xi)^* v(x)\|\leq\|a\|_\infty \|\hat{u}(\xi)\|\|v(x)\| \in L^1 (\R^{2n})$. 

\medskip Aplicando Fubini e o Teorema da Converg\^encia Dominada,
temos
\[\dps\lim_{m,k} <O(a_{m,k})u,v> = 
\dps\lim_{m,k}\dps\int\left[\dps\int e^{-ix\xi}\hat{u}(\xi)^*a_{m,k}(x,\xi)^*\dc\xi\right]v(x)dx
= \] \[ = \dps\int\left[\dps\int e^{-ix\xi}\hat{u}(\xi)^*a(x,\xi)^*\dc\xi\right]v(x)dx = <O(a)u,v>.\]

\noindent Ou seja, $<O(a)u,v> = <u,O(p)v>.$ \,\,\, $\ai$\end{remark}

\begin{preprop}\label{4.6} {\rm Dada $a \in CB^\infty(\R^{2n},A)$, $O(a)$ \'e 
Heisenberg-suave.}\end{preprop}

\begin{remark} 

Seja $T=a(x,D)$ e $T_{z,\zeta} = E^{-1}_{z,\zeta}TE_{z,\zeta}$
(como no cap\'\i tulo 1). Vemos que o s\'\i mbolo de $T_{z,\zeta}$ \'e $a_{z,\zeta}$, com $a_{z,\zeta}(x,\xi) = a(x+z,\xi +\zeta)$, (da mesma maneira como
no caso $A=\C$ \cite[2.20]{2}). Queremos provar que
a aplica\c c\~ao $(z,\zeta) \mapsto T_{z,\zeta}$ \'e $C^\infty$. 

\medskip Vejamos que $(z,\zeta) \mapsto T_{z,\zeta}$ \'e cont\'\i nua
em $(z_0,\zeta_0)$, i.e., que $\forall \epsi >0, \exists \delta >0$, tal que, se
$\|(z,\zeta)-(z_o.\zeta_o)\|<\delta$, ent\~ao $\|T_{z,\zeta}-T_{z_0.\zeta_0}\|<\epsi$.

\medskip Vimos no cap\'\i tulo 3 (\ref{3.8}) que existe uma constante $k \in \R^+$,
independente de $a$, tal que $\|O(a)\| \leq k\pi(a)$, onde $\pi(a) = \dps\sup
\{\|\partial^\beta_x\partial^\gamma_\xi a\|_\infty, \hspace{0,3cm} \beta,\gamma
\leq(1,\cdots ,1)\}$.

\medskip Como $T_{z,\zeta}-T_{z_0,\zeta_0}$ tem
s\'\i mbolo $a_{z,\zeta}-a_{z_0,\zeta_0}$, temos que $\|T_{z,\zeta}-T_{z_0,\zeta_0}\|
\leq k\pi(a_{z,\zeta}-a_{z_0,\zeta_0})$. 

\medskip Dados $\beta,\gamma\leq(1,\cdots,1)$,
\[\|\partial^\beta_x\partial^\gamma_\xi(a_{z,\zeta}-a_{z_0,\zeta_0})\|_\infty = 
\dps\sup_{(x,\xi)\in\R^{2n}}\|\partial^\beta_x\partial^\gamma_\xi
a_{z,\zeta}(x,\xi)-\partial^\beta_x\partial^\gamma_\xi a_{z_0,\zeta_0}(x,\xi)\| =\]
\[= \dps\sup_{(x,\xi)\in\R^{2n}}\|(\partial^\beta_x\partial^\gamma_\xi a)(x+z.\xi+\zeta) -
(\partial^\beta_x\partial^\gamma_\xi)(x+z_0,\xi+\zeta_0)\|.\]

\noindent Como $\partial^\beta_x\partial^\gamma_\xi a$ \'e uniformemente cont\'\i nua, 
$\forall \epsi>0,\exists\sigma>0$ tal que, se
\[ \|(y,\eta)-(y',\eta')\|<\sigma \Rightarrow \|\partial^\beta_x\partial^\gamma_\xi a(y,\eta) -
\partial^\beta_x\partial^\gamma_\xi a(y',\eta')\| < \epsi \] se

\[\|(z,\zeta)-(z_0,\zeta_0)\|<\sigma\Rightarrow\|(z+x,\zeta+\xi)-(z_0+x,\zeta_0+\xi)\|<\sigma,
\hspace{0,3cm} \,\,\, \forall(x,\xi)\in\R^{2n} \Rightarrow\]
\[\Rightarrow\|(\partial^\beta_x\partial^\gamma_\xi a)
(x+z,\xi+\zeta)-(\partial^\beta_x\partial^\gamma_\xi a)(x+z_0,\xi+\zeta_0)\|<\epsi, \,\,\, \forall(x,\xi)\in\R^{2n}.\] Portanto, $\|\partial^\beta_x\partial^\gamma_\xi(a_{z,\zeta}-a_{
z_0,\zeta_0})\|_\infty<\epsi$, \,\,\, se $\|(z,\zeta)-(z_0,\zeta_0)\|<\sigma.$

\medskip Notemos que $\sigma$  depende de $\beta$ e $\gamma$. Como $\beta,\gamma$
s\~ao finitos, tomamos $\delta>0$ como sendo o m\'\i nimo dos $\sigma$ acima e
temos que $\forall\epsi>0$, $\exists\delta>0$ tal que,  $\|(z,\zeta)-(z_0,\zeta_0)\|<\delta
\Rightarrow \pi(a_{z,\zeta}-a_{z_0,\zeta_0})<\epsi.$ Logo, $\|T_{z,\zeta}
-T_{z_0,\zeta_0}\|<\epsi$ e $(z,\zeta)\mapsto T_{z,\zeta}$ \'e cont\'\i nua em $\R^{2n}$.

\medskip
Por outro lado, se $\{e_i\}$ \'e a base can\^onica de $\R^n$, seja

\[\partial_i T_{z,\zeta} = \dps\lim_{h\rightarrow 0}\frac{T_{z+he_i,\zeta}-T_{z,\zeta}}{h} =
\dps\lim_{h\rightarrow 0}\frac{O(a_{z+he_i,\zeta})-O(a_{z,\zeta})}{h}.\]

Provemos que $\partial_i T_{z,\zeta} = O(\partial_i a_{z,\zeta})$, onde
$\partial_i a_{z,\zeta}(x,\xi) =\dps\lim_{h\rightarrow 0}\frac{a_{z,\zeta}(x+he_i,\xi)
-a_{z,\zeta}(x,\xi)}{h}$ e $a_{z,\zeta}(x,\xi)=a(x+z,\xi+\zeta)$.

Temos, como antes, que
\[\left\|\frac{T_{z+he_i,\zeta}-T_{z,\zeta}}{h}-O(\partial_i a_{z,\zeta})\right\|=\left\|\frac{O(a_{z+he_i,\zeta}
-a_{z,\zeta}-h\partial_i a_{z,\zeta})}{h}\right\|\]\[\leq \frac{k}{h}\pi(a_{z+he_i,\zeta}-a_{z,\zeta}-h\partial_i a_{z,\zeta}).\]
Dados $(x,\xi)\in\R^{2n}$, $h\geq 0$ (o caso $h\leq 0$ \'e an\'alogo), pelo
Teorema Fundamental do C\'alculo, (\ref{A.13}),
\[a(x+he_i+z,\xi+\zeta)-a(x+z,\xi+\zeta)-h\partial_i a(x+z,\xi+\zeta) = \]
\[= \dps\int_{[0,h]}\partial_i a(x+z+te_i,\xi+\zeta)dt - h\partial_i a(x+z,\xi+\zeta) =\]
\[ = \dps\int_{[0,h]} \left[\partial_i a(x+z+te_i,\xi+\zeta)-\partial_i a(x+z,\xi+\zeta)\right]dt =\]
\[ = \dps\int_{[0,h]}\dps\int_{[0,t]}\partial_i(\partial_i a)(x+z+se_i,\xi+\zeta)dsdt = \]
\[ = \dps\int_{[0,h]}\dps\int_{[0,t]}\partial^2_i a(x+z+se_i,\xi+\zeta)dsdt.\]
A \'ultima igualdade determina o significado da nota\c c\~ao $\partial^2_i a$. 

Logo,
\[\|a(x+he_i+z,\xi+\zeta)-a(x+z,\xi+\zeta)-h\partial_i a(x+x,\xi+\zeta)\|\leq\]
\[\leq h^2\dps\sup_{(x,\xi)\in\R^{2n}}\|\partial^2_i a_{z,\zeta}(x,\xi)\|.\]

De modo an\'alogo, dados $\beta,\gamma\leq(1,1,\cdots,1)$, temos que
\[\|(\partial^\beta_x\partial^\gamma_\xi a)(x+z+he_i,\xi+\zeta)-\partial^\beta_x\partial^\gamma_\xi a)(x+z,\xi+\zeta) - \]
\[ -h(\partial^\beta_x\partial^\gamma_\xi a)(x+z,\xi+\zeta)\| \leq\]
\[  h^2\dps\sup_{(x,\xi)\in\R^{2n}}\|\partial^2_i\partial^\beta_x\partial^\gamma_\xi a_{z,\zeta}(x,\xi)\| \leq h^2\pi(\partial^2_i a_{z,\zeta}).\]
Portanto,
\[\dps\lim_{h\rightarrow 0}\left\|\frac{O(a_{z+he_i,\zeta})-O(a_{z,\zeta})}{h} -
O(\partial_i a_{z,\zeta})\right\| \leq \dps\lim_{h\rightarrow 0} k.h\pi(\partial^2_i a_{z,\zeta})
= 0.\] ou seja, $\partial_i T_{z,\zeta}= O(\partial_i a_{z,\zeta})$.

%

\medskip  Do mesmo modo, temos que $\partial^\beta_z\partial^\gamma_\zeta
T_{z,\zeta} = O\left(\partial^\beta_z\partial^\gamma_\zeta a_{z,\zeta}\right)$, e temos, como antes,
que $(z,\zeta)\mapsto \partial^\beta_z\partial^\gamma_\zeta T_{z,\zeta}$ \'e
cont\'\i nua, i.e., $(z,\zeta)\mapsto T_{z,\zeta}$ \'e $C^\infty$. $\ai$ \end{remark}

\begin{precorol}\label{4.7}{\rm   Dada $a\in CB^\infty(\R^{2n},A), \,\,\, O(a)\in {\cal{HS}}$.                                            }\end{precorol}

\bigskip Nosso objetivo agora \'e provar que, se vale \ref{4.1}, dado um operador $T\in {\cal{HS}}$
tal que para toda $G\in CB^\infty(\R^n,A)$ comuta com $R_G$, ent\~ao existe $F\in CB^\infty(\R^n,A)$ de modo que $T = L_F$.

\medskip Na verdade, provaremos o resultado acima com uma hip\'otese mais fraca,
que o operador $T\in\cal{HS}$ comuta com $R_g$, para toda $g\in\sa$.

\bigskip Consideraremos primeiramente o caso $n=1$. 

\begin{prelemma}\label{4.8} {\rm  Dada $f\in CB^\infty(\R,A)$ temos que
\[\dps\int \gamma(t)(1-\partial_t)^2 f(t)dt = f(0),\]  onde
\begin{center}$\gamma(t) = \left\{ \begin{array}{l} te^{-t} \hspace{0,5cm} {\rm se} \,\, t\geq 0 \\
0 \hspace{1cm} {\rm se} \,\, t<0 \end{array} \right. $ \end{center} .   }\end{prelemma}

\begin{remark}    Como $f\in CB^\infty(\R,A)$, $\gamma(t)(1-\partial_t)^2 f(t) \in L^1(\R,A)$. 
\[\dps\int\gamma(t)(1-\partial_t)^2 f(t)dt = \dps\int_{\R^+}te^{-t}(f(t) - 2\partial_t f(t) + \partial^2_t
f(t))dt.\]  Fazendo integra\c c\~ao por partes (ver \ref{A.14}), para $R\in\R^+$,   
\[\dps\int_{[0,R]} te^{-t}\partial_t f(t)dt = Re^{-R}f(R) - \dps\int_{[0,R]}\partial_t(te^{-t})f(t)dt.\]
\[\dps\int_{[0,R]}te^{-t}\partial^2_t f(t)dt = Re^{-R}\partial_t f(R) - \dps\int_{[0,R]}\partial_t(te^{-t})
\partial_t f(t)dt = \] \[ = Re^{-R}\partial_t f(R) - \left[(e^{-R} - Re^{-R})f(R) - f(0) - \dps\int_{[0,R]}\partial^2_t(te^{-t})f(t)dt\right].\] Assim, como $\dps\int_{R^+}te^{-t}(1-\partial_t)^2 f(t)dt = \dps\lim_{R\rightarrow\infty}\dps\int_{[0,R]}te^{-t}(1-\partial_t)^2 f(t)dt$, pela
observa\c c\~ao \ref{A.6},  temos:
\[ \dps\int\gamma(t)(1-\partial_t)^2 f(t)dt = \dps\int_{R^+}\left[te^{-t} + 2\partial_t(te^{-t}) 
+ \partial^2_t(te^{-t})\right]f(t)dt + f(0) = f(0). \ai\]  \end{remark}                                    

\begin{prenota}\label{4.8ummeio}{\rm  Para o caso geral, consideramos $\overline{\gamma}
(x) = \dps\prod^n_{j=1}\gamma(x_j)$, onde $x=(x_1,x_2,\cdots,x_n)$ e temos,
de modo an\'alogo, que     $\dps\int\overline{\gamma}(x)\prod^n_{j=1}(1-\partial_{x_j})^2
f(x)dx = f(0)$, se $f\in CB^\infty(\R^n,A)$. 

De fato, pelo Teorema de Fubini (\ref{A.12}), temos
\[\dps\int\overline{\gamma}(x)\dps\prod^n_{j=1}(1-\partial_{x_j})^2f(x)dx  =\]
\[ = \dps\int\dps\prod^n_{j=2}\gamma(x_j)\left(\dps\int\gamma(x_1)(1-\partial_{x_1})^2\left[
\dps\prod^n_{j=2}(1-\partial_{x_j})^2f(x_1,x_2,\cdots,x_n)\right]dx_1\right)dx_2\cdots dx_n =\] 
\[ = \dps\int\prod^n_{j=3}\gamma(x_j)\left(\dps\int\gamma_{x_2}(1-\partial_{x_2})^2
\left[\prod^n_{j=3}(1-\partial_{x_j})^2 f(0,x_2,\cdots,x_n)\right]dx_2\right)dx_3\cdots dx_n = \]
\[ = \cdots  = f(0,\cdots,0) = f(0).\]      }\end{prenota}

\bigskip 
\bigskip
A proposi\c c\~ao seguinte, para o caso $A=\C$, foi provada por Cordes em 
\cite[corol\'ario 2.4]{2}.

\bigskip
\begin{preprop}\label{4.9}{\rm  Dadas $a\in CB^\infty(\R^2,A)$ e $b = (1+\partial_y)^2(1+\partial_\eta)^2 a$    ent\~ao temos que
\[a(x,\xi) = \dps\int\gamma(y)\gamma(\eta)b(x-y,\xi-\eta)dyd\eta.\]     }\end{preprop}

\begin{remark}  Notemos que
\[b(x-y,\xi-\eta) = [(1+\partial_y)^2(1+\partial_\eta)^2 a](x-y,\xi-\eta) = (1-\partial_y)^2(1-\partial_\eta)^2[a(x-y,\xi-\eta)].\] Logo,
\[\dps\int\gamma(y)\gamma(\eta)b(x-y,\xi-\eta)dyd\eta = \dps\int\gamma(y)\gamma(\eta)(1-\partial_y)^2(1-\partial_\eta)^2a(x-y,\xi-\eta)dyd\eta.\]
Como $\gamma(y)\gamma(\eta)(1-\partial_y)^2(1-\partial_\eta)^2 a(x-y,\xi-\eta) \in
L^1(\R^2,A)$, podemos aplicar o Teorema de Fubini e, considerando que $(1-\partial_\eta)^2 
a(x-y,\xi-\eta)\in CB^\infty(\R,A)$ para $x, \xi, \eta$  fixados, aplicamos o lema \ref{4.8}:
\[ \dps\int \gamma(\eta) \dps\int \gamma(y)(1-\partial_y)^2 (1-\partial_\eta)^2 a(x-y,\xi-\eta)dyd\eta = \dps\int \gamma(\eta)(1-\partial_y)^2 a(x,\xi-\eta)d\eta.\]  Novamente,
para $x,\xi$ fixos, pelo lema \ref{4.8}, $\dps\int\gamma(\eta)(1-\partial_\eta)^2 a(x,\xi-\eta)d\eta = a(x,\xi)$                                            $\ai$ \end{remark}

\begin{precorol}\label{4.10} {\rm   Nas hip\'oteses da proposi\c c\~ao
    \ref{4.9}, temos

\begin{center}$a(z,\zeta) = \dps\int \gamma(-x)\gamma(-\xi)b(x+z,\xi+\zeta)dxd\xi $.\end{center}                                           } \end{precorol}

\begin{prenota}\label{4.10ummeio} {\rm   Para o caso geral, dada $a\in CB^\infty(\R^{2n},A)$, seja $b = \dps\prod^n_{j=1}(1+\partial_{y_j})^2(1+\partial_{\eta_j})^2 a$, e
temos, de modo an\'alogo, que $a(z,\zeta)=\dps\int\overline{\gamma}(-x)
\overline{\gamma}(-\xi)b(x+z,\xi+\zeta)dx d\xi$, usando a nota\c c\~ao
da nota \ref{4.8ummeio}.  }   \end{prenota}

\begin{prelemma}\label{4.11}{\rm    Sejam $F$ a transformada de Fourier e $E_{z,\zeta}$
definido como no cap\'\i tulo 1. Ent\~ao temos:
\begin{enumerate} 
\item $FE_{z,\zeta} = (E_{-\zeta,z})^{-1}F$
\item $F(E_{z,\zeta})^{-1} = E_{-\zeta,z}F$.  \end{enumerate}}\end{prelemma}

\begin{remark} 

Dada $u\in S^A(\R^n)$,

 \begin{center} {\bf 1.}  \hspace{2,2cm} $FE_{z,\zeta}u(x) = \dps\int e^{-ixy}E_{z,\zeta}u(y)\dc y = 
\dps\int e^{-ixy}e^{i\zeta y}u(y-z)\dc y = $\end{center}
\begin{center}$ = \dps\int e^{-i(y+z)(x-\zeta)}u(y)\dc y = e^{-iz(x-\zeta)}Fu(x-\zeta) = 
T_\zeta M_{-z}Fu(x)
= (E_{-\zeta,z})^{-1}Fu(x).$\end{center}

\bigskip
\begin{center} {\bf 2.} \hspace{0,7cm}  $FE_{z,\zeta}^{-1}u(x) = \dps\int e^{-ixy}T_{-z}M_{-\zeta}u(y)\dc y =
\dps\int e^{-ixy}e^{-i\zeta(y+z)}u(y+z)\dc y = $\end{center}
\begin{center} $\dps\int e^{-ix(y-z)}e^{-i\zeta y}u(y)\dc y = 
e^{ixz}\dps\int e^{-iy(x+\zeta)}u(y)\dc y = e^{ixz}Fu(x+\zeta) =
M_zT_{-\zeta}Fu(x) = E_{-\zeta,z}Fu(x). \,\,\, \ai $\end{center}\end{remark}

\begin{preobs}\label{4.11ummeio} {\rm  Lembremos que, como vimos na demonstra\c c\~ao da proposi\c c\~ao \ref{2.2}, dadas $g,u\in\sa$, temos que $R_g u(x) = \dps\int
 e^{ix\xi}Fu(\xi)g(x+J\xi)\dc\xi$}\end{preobs}

\begin{prelemma}\label{4.12} {\rm Dada $g\in S(\R^n)$, vale que
\begin{enumerate}
\item $E_{z,\zeta}^{-1}R_g = R_{T_{-z - J\zeta}g}E_{z,\zeta}^{-1}$
\item $E_{z,\zeta}R_g = R_{T_{z + J\zeta}g}E_{z,\zeta}.$
\end{enumerate}} \end{prelemma}

\begin{remark} 

\bigskip

{\bf 1.} Dada $u\in S(\R^n)$, pelo lema \ref{4.11},

\[E_{z,\zeta}^{-1}R_g u(x) = e^{-i\zeta(x+z)}R_g u(x+z) = e^{-i\zeta(x+z)}\dps\int
e^{i(x+z)\xi}Fu(\xi)g(x+z+J\xi)\dc\xi.\]

\[R_{T_{-z - J\zeta}g}E_{z,\zeta}^{-1}u(x) = 
\dps\int e^{ix\xi}F(E_{z,\zeta}^{-1}u)(\xi)T_{-z - J\zeta}g(x + J\xi)
\dc\xi = \] \[ = \dps\int e^{ix\xi}E_{-\zeta,z}Fu(\xi)g(x+z + J(\xi+\zeta)\dc\xi = 
 \dps\int e^{ix\xi}e^{iz\xi}Fu(\xi+\zeta)g(x+z + J(\xi+\zeta))\dc\xi = \]
\[ = \dps\int e^{i(x+z)(\xi-\zeta)}Fu(\xi)g(x+z + J\xi)\dc\xi = E_{z,\zeta}^{-1}R_g u(x).\]

\begin{center}{\bf 2.} \hspace{0,5cm} $E_{z,\zeta}R_g u(x) = e^{i\zeta x}R_g u(x-z) = 
e^{i\zeta x}\dps\int e^{i(x-z)\xi}Fu(\xi)g(x-z + J\zeta)\dc\xi.$\end{center}

\[R_{T_{z+J\zeta}g}E_{z,\zeta}u(x) = 
\dps\int e^{ix\xi}F(E_{z,\zeta}u)(\xi)T_{z + J\zeta}g(x + J\xi)\dc\xi = \]
\[ =\dps\int e^{ix\xi}(E_{-\zeta,z})^{-1}Fu(\xi)g(x-z + J(\xi-\zeta))\dc\xi =
\dps\int e^{ix\xi}e^{-iz(\xi-\zeta)}Fu(\xi-\zeta)g(x-z + J(\xi-\zeta))\dc\xi =\]
\[ = \dps\int e^{ix(\xi+\zeta)}e^{-iz\xi}Fu(\xi)g(x-z+J\xi)\dc\xi =
E_{z,\zeta}R_g u(x). \,\,\, \ai \]\end{remark}

\begin{preprop}\label{4.13} {\rm  Dado $T\in {\cal B}^*(\overline{S^A(\R^n)})$, se $[T,R_g]=0, \forall g\in S^A(\R^n)$, ent\~ao $[T_{z,\zeta},R_g]=0, \forall g\in S^A(\R^n)$. }\end{preprop}

\begin{remark} Pelo lema \ref{4.12}, temos 
\[T_{z,\zeta}R_g = E_{z,\zeta}^{-1}TE_{z,\zeta}R_g = 
E_{z,\zeta}^{-1}TR_{T_{z+J\zeta}g}E_{z,\zeta} = 
E_{z,\zeta}^{-1}R_{T_{z+J\zeta}g}TE_{z,\zeta} = R_g E_{z,\zeta}^{-1}TE_{z,\zeta} = R_g T_{z,\zeta}. \,\,\, \ai\]\end{remark}

\begin{preprop}\label{4.14} {\rm  Dadas $a\in CB^\infty(\R^{2n},A)$ e $b = \dps\prod^n_{j=1}(1+\partial_{x_j})^2(1+\partial_{\xi_j})^2 a$, se $O(a)\in \cal{HS}$ e $[O(a),R_g]=0$, para todo $ g\in S^A(\R^n)$,
ent\~ao $O(b)\in \cal{HS}$ e $[O(b),R_g]=0, \forall g\in S^A(\R^n)$.      }\end{preprop}

\begin{remark} Notemos que, como visto na demonstra\c c\~ao de \ref{4.6}, se 
$T = O(a)$ e $B = O(b)$, $B_{z,\zeta} = \dps\prod^n_{j=1} (1+\partial_{z_j})^2
(1+\partial_{\zeta_j})^2 T_{z,\zeta}.$

\medskip Logo, se $(z,\zeta) \mapsto T_{z,\zeta}$ \'e $C^\infty$, $(z,\zeta)
\mapsto B_{z,\zeta}$ \'e $C^\infty$. 

\medskip Vimos na proposi\c c\~ao \ref{4.13}
que se $[T,R_g]=0, \forall g\in S^A(\R^n)$, ent\~ao $[T_{z,\zeta},R_g]=0,
\forall g\in S^A(\R^n)$. Assim, se $\{e_j\}$ \'e a base can\^onica de $\R^n$,
\[(\partial_{z_j} T_{z,\zeta})R_g = \dps\lim_{h\rightarrow 0}\frac{T_{z+he_j,\zeta}
-T_{z,\zeta}}{h}R_g = \]\[\dps\lim_{h\rightarrow 0}\frac{T_{z+he_j,
\zeta}R_g - T_{z,\zeta}R_g}{h} = \dps\lim_{h\rightarrow 0}R_g\frac{T_{z+he_j,\zeta}-
T_{z,\zeta}}{h} = 
R_g(\partial_{z_j} T_{z,\zeta}).\] 
De modo an\'alogo, provamos que $[R_g,B]=0, \forall g\in S^A(\R^n)$. \,\,\, $\ai$\end{remark}

\begin{pretheo}\label{4.15} {\rm Dado $T\in {\cal B}^*(\overline{S^A(\R^n)})$, Heisenberg-suave, tal que 
para toda $ g\in \sa, [T,R_g]=0$. Supondo que vale \ref{4.1}, então 
existe $F\in CB^\infty(\R^n,A)$ de modo que $T = L_F$.    }\end{pretheo}

\begin{remark}  Por \ref{4.1} temos que existe $a\in CB^\infty(\R^{2n},A)$ 
de modo que $T=O(a)$. Sejam $b = \dps\prod^n_{j=1}(1+\partial_{x_j})^2
(1+\partial_{\xi_j})^2 a$ e
$B = O(b)$. Pela proposi\c c\~ao \ref{4.14}, temos que $[B,R_g]=0, \forall g\in S^A(\R^n)$.
Logo, $[B_{z,\zeta},R_g]=0, \forall g\in S^A(\R^n)$, pela proposi\c c\~ao \ref{4.13}.
Assim, pelo visto na demonstra\c c\~ao do teorema \ref{2.7}, $B_{z,\zeta} = 
B_{z-J\zeta,0}.$

\medskip Por outro lado, o s\'\i mbolo de $B_{z,\zeta}$ \'e $b(x+z,\xi+\zeta)$ e
o s\'\i mbolo de $B_{z-J\zeta,0}$ \'e $b(x+z-J\zeta,\xi)$. Portanto,
por \ref{4.1}, temos que $b(x+z,\xi+\zeta) = b(x+z-J\zeta,\xi)$.
Assim, pela nota  \ref{4.10ummeio}, temos:
\[a(z,\zeta) = \dps\int\overline{\gamma}(-x)\overline{\gamma}(-\xi)b(x+z,\xi+\zeta)dxd\xi = \]
\[= \dps\int\overline{\gamma}(-x)\overline{\gamma}(-\xi)b(x+z-J\zeta,\xi)dxd\xi =
a(z - J\zeta,0). \]
Tomando $F(z) = a(z,0)$, temos $F\in CB^\infty(\R^n,A)$ e 
$a(z,\zeta)=F(z-J\zeta)$, e, pela observa\c c\~ao \ref{1umcuarto}, portanto,
$T = L_F$. \,\,\, $\ai$\end{remark}

\vspace{1cm}


\renewcommand{\chaptername}{Apêndice} 
\renewcommand{\thechapter}{\Alph{chapter}}
\setcounter{chapter}{0} 




\chapter{}






\vspace{1cm}
Veremos agora alguns resultados sobre integrais de fun\c c\~oes que chegam num espa\c co de Banach, Integral de Bochner. Para isto, nos baseamos fundamentalmente no ap\^endice
E de \cite{13}.

\vspace{1.0cm}

Seja $E$ um espa\c co de Banach separ\'avel e seja ${\cal B}(E)$ a $\sigma$-\'algebra
gerada pelos sub-conjuntos abertos de $E$.

\begin{predef}\label{A.1} {\rm Dizemos que uma fun\c c\~ao $f:\R^n \rightarrow E$ \'e
{\em mensur\'avel} se \'e mensur\'avel com respeito a ${\cal B}(\R^n)$ e a ${\cal B}(E)$.}
\end{predef}

\'E f\'acil ver que, se $f$ \'e mensur\'avel, ent\~ao $x \mapsto \|f(x)\|$ \'e mensur\'avel.

\begin{predef}\label{A.2} {\rm Uma fun\c c\~ao mensur\'avel $f:\R^n \rightarrow E$ \'e 
dita {\em simples} se $f(\R^n)$ \'e um subconjnto finito de $E$. Na verdade, se
$f(\R^n) = \{a_1, \cdots , a_p\}$, seja $J_l = f^{-1}(a_l)$; ent\~ao,
$f (x) = \dps\sum^p_{l=1}a_l\chi_{J_l}(x)$, onde $\chi_{J_l}$ \'e
a fun\c c\~ao caracter\'\i stica de $J_l$ e onde os conjuntos $J_l$ s\~ao disjuntos
e pertencem a ${\cal B}(\R^n)$.} \end{predef}

\begin{prelemma}\label{A.3} {\rm Dada $f:\R^n \rightarrow E$, seja $D \subseteq  E$
enumer\'avel denso em $f(\R^n)$. Ent\~ao, $\forall\epsi > 0, \forall a \in f(\R^n)$, existem
$d \in D$ e $r \in \Q $,  tais que $\|rd\| \leq \|a\|$ e $\|rd- a\| < \epsi$.}\end{prelemma}

\begin{remark} Dado $1>\epsi>0$, como $D$ \'e denso em $f(\R^n)$, existe $d \in D$
tal que $\|a - d\| < \epsi$. Seja $r\in\Q^+$ tal que $1-\epsi<r<1$ e $r\leq\frac{\|a\|}{\|d\|}$.
Assim, $\|rd\| \leq \|a\|$ e temos 
$$\|rd-a\|\leq\|d-rd\|+\|a-d\|\leq\|d\|(1-r)+\epsi\leq\frac{\|a\|}{r}\epsi +\epsi\leq\|a\|\frac{\epsi}{1-\epsi}+\epsi.$$ 

Dado $\epsi_1 > 0$, escolhemos $\epsi > 0$ de
modo que $\|a\|\frac{\epsi}{1-\epsi} + \epsi \leq \epsi_1$ e temos $\|rd - a \| < \epsi_1$  e 
$\|d\| \leq \|a\|$     $\ai$\end{remark}

\begin{preprop}\label{A.4} {\rm Dada $f:\R^n \rightarrow E$, mensur\'avel, existe uma seq\"u\^encia $(f_k)_{k \in N}$ de fun\c c\~oes simples tal que $\forall x \in \R^n,  f(x) = \dps\lim_{k \rightarrow \infty}f_k(x), $ e $\|f_k(x)\| \leq \|f(x)\|,  \forall  k \in \N.$  }\end{preprop}

\begin{remark}  Suponhamos que $f$ \'e n\~ao nula. 

Seja $D \subseteq E$ enumer\'avel, denso em $f(\R^n)$.  Seja $D'$  o conjunto dos m\'ultiplos racionais de $D$,
$D'= \{c_k\}_{k \in \N}$. Podemos assumir que $c_1 = 0$. 

Para cada $x \in \R^n$, e cada $k \in \N$, seja $E_k(x) \subseteq E$ dado por:

\[ E_k(x) = \{ c_j \in D', j \leq k \,\, {\rm e} \,\, \|c_j\| \leq \|f(x)\| \} .\]

Como $c_1=0$, $E_k(x)$ \'e n\~ao vazio.

Vamos construir a seq\"u\^encia $(f_k)_{k \in \N}$, tomando
$f_k(x)$ o elemento de $E_k(x)$ que est\'a mais
pr\'oximo de $f(x)$.

Consideremos
\[A_k(x) = \{c_i \in E_k(x) / \|f(x) - c_i\| = \dps\inf \{\|f(x) - c_j\|, c_j \in
E_k(x)\}\}.\]

Seja $I_k(x)$ o conjunto dos \'\i ndices dos $c_i \in A_k(x)$ e seja
$i_0 = \dps\min I_k(x)$ e $f_k(x) = c_{i_0}$. 

Como os conjuntos $\{x \in \R^n, f_k(x) = c_{i_0}\}$
podem ser descritos por meio de desigualdades 
envolvendo $\|f(x)\|, \|c_j\|, j = 1, \cdots , k, \,\, {\rm e} \,\, \|f(x)-c_j\|, j = 1, \cdots , k$,
cada $f_k$ \'e mensur\'avel.  \'E f\'acil ver, ent\~ao, que $f_k$ \'e simples, 
$\forall k \in\N$. Al\'em disso, por constru\c c\~ao de $f_k$ e de $E_k(x)$, 
temos que $\|f_k(x)\| \leq \|f(x)\|, \forall x \in\R, \forall k \in\N$.

Vejamos agora que $(f_k)_{k \in \N}$ converge pontualmente para $f$. Consideremos:


\begin{enumerate}

\item Pelo lema \ref{A.3}, para todo $ \epsi > 0,$ para todo $ x \in \R^n$,  existe $c_{h_0}
\in D'$ tal que $\|c_{h_0}\| \leq \|f(x)\|$ (isto implica $c_{h_0} \in E_{h_0}(x)$)  e
$\|c_{h_0} - f(x)\| < \epsi $.

\item Se $k \geq h_0$, $E_{h_0}(x) \subseteq E_k(x).$

\item $\|f(x) - f_k(x)\| \leq \|f(x) - c_k\|$, para qualquer $c_k \in E_k(x)$.

\end{enumerate}

Assim, dado $\epsi > 0$, seja $h_0\in\N$ como em (1). Ent\~ao temos, 
$\forall k \geq h_0$, $\|f(x) - f_k(x)\| \leq \|f(x)-c_{h_0}\| 
< \epsi$.
$\ai$\end{remark}

\begin{predef}\label{A.5} {\rm  Seja $f:\R^n \rightarrow E$, mensur\'avel. Dizemos 
que $f$ \'e {\em Bochner integr\'avel (integr\'avel)} se a fun\c c\~ao
$x \mapsto \|f(x)\|$ \'e integr\'avel.       }\end{predef}

Dada $f:\R^n \rightarrow E$, simples e integr\'avel, com $f(x) = 
\dps\sum^p_{l=1}a_l\chi_{J_l}, \,\,\, a_l\in E \,\, {\rm e} \,\, J_l\in{\cal{B}}(\R^n)$,
definimos a integral de $f$ como segue:

\[ \dps\int f(x)dx = \dps\sum^p_{l=1}a_l\mu(J_l),\] onde $\mu$ \'e a medida
de Lebesgue em $\R^n$. \'E claro que $f$ integr\'avel implica que $\mu(J_l)
< \infty,\,\,  l=1, \cdots , p.$

\medskip
\'E f\'acil ver que \[\left\|\dps\int f(x)dx\right\| \leq \dps\int\|f(x)\|dx\] e que, 
se $\alpha,\beta \in \C$ e $g$ uma fun\c c\~ao simples integr\'avel,

\[ \dps\int (\alpha f + \beta g)(x)dx = \alpha\dps\int f(x)dx + \beta\dps\int g(x)dx .\]

Dada $f:\R^n \rightarrow E$ integr\'avel, seja $(f_k)_{k \in\N}$ uma 
seq\"u\^encia de fun\c c\~oes simples tal que $f(x) = \dps\lim_{k \rightarrow 
\infty}f_k(x), \forall x\in\R^n$ e tal que $x \mapsto \dps\sup_{k}\|f_k(x)\| = \|f(x)\|$
\'e integr\'avel (ver \ref{A.4}). Como $\|f_k(x)-f(x)\|\leq 2\|f(x)\|$, pelo Teorema 
da Converg\^encia Dominada, $\dps\int\|f_k(x)-f(x)\|dx \underset{x\rightarrow\infty}\longrightarrow 0.$ Logo, a seq\"u\^encia 
$(\int f_k)_{k\in\N}$ \'e de Cauchy em $E$ e portanto converge.

\bigskip Definimos

\[ \dps\int f(x)dx = \dps\lim_{k\rightarrow\infty}\dps\int f_k(x)dx.\] \'E f\'acil
ver que esta defini\c c\~ao n\~ao depende da escolha da seq\"u\^encia $(f_k)_{k\in\N}$
\cite[pg.352]{13}.

\medskip

\begin{preprop}\label{A.6} {\rm   Seja $f:\R^n \rightarrow E$, integr\'avel. Ent\~ao

\[ \left\| \dps\int f(x)dx \right\| \leq \dps\int \|f(x)\|dx.\]                   }\end{preprop}

\begin{remark}   Seja $(f_k)_{k\in\N}$ uma seq\"u\^encia de fun\c c\~oes simples (como
na proposi\c c\~ao \ref{A.5}),  tal que $\forall x\in\R^n, \,\,\, \dps\sup_{k}\|f_k(x)\|
\leq \|f(x)\|$ e $f(x) = \dps\lim_kf_k(x).$

%

\medskip \noindent Logo, $\left\|\dps\int f_k(x)dx\right\| \leq \dps\int \|f_k(x)\|dx \leq \dps\int\|f(x)\|dx$. Como $\dps\int f(x)dx = \dps\lim_{k\rightarrow\infty}\int f_k(x)dx$, vale a tese.
$\ai$\end{remark}

\begin{pretheoCD}\label{A.7}
{\rm   Sejam $f, f_1, f_2, \cdots $   fun\c c\~oes de $\R^n$ em $E$, mensur\'aveis, 
tais que $f(x) = \dps\lim_{k\rightarrow\infty}f_k(x)$ para quase todo $x\in\R^n$ 
e $\|f_k(x)\| \leq g(x)$, q.s., para todo $k\in\N$, onde $g:\R^n\rightarrow 
[0,+\infty]$, integr\'avel. Ent\~ao, $f, f_1, \cdots $ s\~ao integr\'aveis e
$\dps\int f(x)dx = \dps\lim_{k\rightarrow\infty}\dps\int f_k(x)dx.$                             }\end{pretheoCD}

\bigskip
\begin{remark} Como $\|f_k(x) \| \leq g(x)$ q.s., $\forall k\in\N$, a integrabilidade 
de $f_1, f_2, \cdots$
\'e imediata. Al\'em disso, temos $\|f(x)\|  \leq g(x)$ q.s. e $f$ \'e
integr\'avel. Portanto, $\|(f_k-f)(x)\| \leq
2g(x)$ q.s.. Pelo Teorema da Converg\^encia Dominada (caso escalar) temos ent\~ao 
que $\dps\lim_{k\rightarrow\infty}\dps\int \|(f_k-f)(x)\|dx = 0$, e,
pela proposi\c c\~ao \ref{A.6}, temos $\|\dps\int(f_k-f)(x)dx\| \leq \dps\int\|(f_k-f)(x)\|dx
\underset{k\rightarrow\infty}\rightarrow 0$. $\ai$\end{remark}

\medskip 
\begin{preobs}\label{A5ummeio} {\rm Se $f:\R^n\rightarrow E$ 
\'e integr\'avel e $J\in{\cal{B}}(\R^n)$,
definimos

\[ \dps\int_J f(x)dx = \dps\int_{\R^n} f(x)\chi_J(x)dx.\] Se $g:\R\rightarrow E$ 
\'e integr\'avel, notemos que 
$\dps\int_{\R} g(t)dt = \dps\lim_{R\rightarrow\infty}
\dps\int_{[-R,R]}g(t)dt$, pois, se $(R_k)_{k\in\N}$ \'e uma seq\"u\^encia
em $\R^+$, crescente, com $\dps\lim_{k\rightarrow\infty}R_k = +\infty$, e 
se $g_k(t)=g(t)\chi_{[-R_k,R_k]}(t)$, temos que, para cada $t\in\R$, $\dps\lim_{k\rightarrow\infty}g_k(t)=g(t)$ e $\dps\sup_k\|g_k(t)\|\leq\|g(t)\|$;
portanto, pelo Teorema da Converg\^encia Dominada \ref{A.7}
\[\dps\int_{\R} g(t)dt = \dps\lim_{k\rightarrow\infty}\dps\int_{\R} g_k(t)dt =
\dps\lim_{k\rightarrow\infty}\dps\int_{[-R_k,R_k]}g(t)dt.\] }\end{preobs}

\begin{predef}\label{A.8} {\rm  Seja ${\cal{L}}^1(\R^n,E)$  o conjunto das
fun\c c\~oes de $\R^n$ em $E$ integr\'aveis. Defini\-mos

\begin{center}$ \|f\|_{L^1} = \dps\int\|f(x)\|dx $, para $f\in{\cal{L}}^1(\R^n,E).$\end{center}
}\end{predef}

\'E claro que $\|\cdot\|_{L^1}$ \'e semi-norma. \'E f\'acil ver que o conjunto
$L^1(\R^n,E)$ das classes de equival\^encia de ${\cal{L}}^1(\R^n,E)$,
com respeito \`a igualdade quase sempre, \'e um espa\c co
vetorial normado onde $\|\cdot\|_{L^1}$ \'e a norma (\cite[proposi\c c\~ao E.4]{13}).
Tamb\'em podemos provar que $L^1(\R^n,E)$ \'e completo adaptando
a demonstra\c c\~ao de \cite[teorema 3.4.1]{13}. 

\medskip De um modo an\'alogo, definimos os espa\c cos de Banach $L^p(\R^n,E)$
das classes de equival\^encia de fun\c c\~oes iguais q.s. que tem a propriedade
$\dps\int \|f(x)\|^pdx < \infty$. A norma \'e dada por $\|f\|_{L^p} = \left(\dps\int\|f(x)\|^pdx
\right)^{\frac{1}{p}}$, para $f\in L^p(\R^n,E)$.

\begin{preprop}\label{A.9} {\rm  As fun\c c\~oes simples do tipo 
$\dps\sum^q_{l=1}a_l\chi_{J_l}, \,\, \mu(J_l) < \infty, a_l \in E$, 
s\~ao densas em $L^p(\R^n,E)$.}\end{preprop}

\begin{remark} \'E f\'acil ver, pela defini\c c\~ao \ref{A.5}, que as fun\c c\~oes 
simples do enunciado pertencem a $L^p(\R^n,E)$. Dada uma fun\c c\~ao
$f \in L^p(\R^n,E)$, seja uma seq\"u\^encia de fun\c c\~oes  simples
$(f_k)_{k\in\N}$, como na proposi\c c\~ao \ref{A.4}, tal que, para quase todo $x \in \R^n$,
vale que $f(x) = \dps\lim_{k\rightarrow\infty}f_k(x)$ e $\|f_k(x)\|\leq\|f(x)\|$.
Assim, $\|(f_k-f)(x)\|^p\leq 2^p\|f(x)\|^p$. Logo, pelo Teorema da
Converg\^encia Dominada, (\ref{A.7}), $\dps\lim_{k\rightarrow\infty}\dps\int \|(f_k-f)(x)\|^pdx
= 0$, o que implica $\|f-f_k\|_{L^p} \underset{k\rightarrow\infty}\rightarrow 0$.                    $\ai$\end{remark}      

\bigskip Um resultado \'util, mas que n\~ao provaremos aqui, \'e o teorema E.9
de \cite{13}, que diz o seguinte:

\begin{theorem}\label{A.10} {\rm Uma fun\c c\~ao $f:\R^n \rightarrow E$ \'e 
mensur\'avel se e somente se, para todo funcional linear cont\'\i nuo $\varphi$ em
$E^*$, a fun\c c\~ao $\varphi \circ f$ \'e mensur\'aval.}\end{theorem}

\begin{preprop}\label{A.11} {\rm Dada $f:\R^n \rightarrow E$, integr\'avel, ent\~ao, 
para qualquer funcional linear cont\'\i nuo $\varphi \in E^*$ vale que 

\[\dps\int \varphi\circ
f(x)dx = \varphi\left(\dps\int f(x)dx\right) \hspace{0,5cm} (*) \] Al\'em disso,
se existe $a_0 \in E$ tal que $\forall \varphi \in E^*, \varphi(a_0) = 
\dps\int \varphi\circ f(x)dx$, ent\~ao $a_0 = \dps\int f(x)dx$.}
\end{preprop}
\begin{remark} Dado qualquer $\varphi \in E^*$, por (\ref{A.10}),
sabemos que $\varphi\circ f$ \'e mensur\'avel. Vejamos  que 
$\varphi\circ f$ \'e integr\'avel.

\[\dps\int|\varphi\circ f(x)|dx \leq \dps\int \|\varphi\|\|f(x)\|dx =
\|\varphi\|\dps\int\|f(x)\|dx.\] Como $f$ \'e integr\'avel, 
$\varphi\circ f$ \'e integr\'avel.

\medskip Suponhamos que $f$ \'e simples, 
$f(x) = \dps\sum^p_{l=1}a_l\chi_{J_l}(x), \hspace{0,5cm} a_l\in E$ \,\, {\rm e}
$\mu(J_l)<\infty$, onde $J_l$ s\~ao disjuntos, $l = 1, \cdots , p$.

\[\varphi\left(\dps\int f(x)dx \right) = \varphi\left(\dps\sum^p_{l=1}a_l\mu(J_l)\right)
= \dps\sum^p_{l=1}\varphi(a_l)\mu(J_l).\]
\[\dps\int \varphi\circ f(x)dx = \dps\int\varphi\left(\dps\sum^p_{l=1}a_l\chi_{J_l}(x)\right)dx
= \dps\int\dps\sum^p_{l=1}\varphi(a_l)\chi_{J_l}(x)dx = \dps\sum^p_{l=1}\varphi(a_l)\mu(J_l).\]  Portanto,
vale $\dps (*)$ para fun\c c\~oes simples.

\medskip Seja agora $f:\R^n \rightarrow E$ integr\'avel e $(f_k)_{k\in\N}$
uma seq\"u\^encia de fun\c c\~oes simples tal que, para quase todo $x \in\R^n$,
$f(x) = \dps\lim_{k\rightarrow\infty}f_k(x)$ e $\dps\sup_{k}\|f_k(x)\|\leq\|f(x)\|$.
Desta maneira, \'e claro que $f_k$ \'e integr\'avel para todo $k\in\N$
e temos que 

\[\varphi\left(\dps\int f_k(x)dx\right) = \dps\int \varphi\circ f_k(x)dx.   \hspace{0,5cm} (1) \]    

\noindent Al\'em disso, $\|f_k(x)\| \leq \|f(x)\| \in L^1(\R^n)$, para quase todo
$x\in\R^n$; ent\~ao, pelo Teorema da Converg\^encia Dominada (\ref{A.7}),
$\dps\lim_{k\rightarrow\infty}\dps\int f_k(x)dx = \dps\int f(x)dx$ e,
portanto, 

\[\dps\lim_{k\rightarrow\infty}\varphi\left(\dps\int f_k(x)dx \right) =
\varphi\left(\dps\int f(x)dx\right).  \hspace{0,5cm} (2) \]

Por outro lado, $|\varphi\circ f_k(x)|\leq \|\varphi\|\|f_k(x)\|\leq
\|\varphi\|\|f(x)\| \in L^1(\R^n)$,
para quase todo $x\in\R^n$ e $\forall k\in\N$. Como $\dps\lim_{k\rightarrow\infty} 
\varphi\circ f_k(x) = \varphi\circ f(x)$ q.s., pelo Teorema da
Converg\^encia Dominada, $\varphi\left(\dps\int f(x)dx\right) \underset{(2)}=
 \dps\lim_{k\rightarrow\infty} \varphi\left(\dps\int f_k(x)dx\right) \underset{(1)}= 
\dps\lim_{k\rightarrow\infty}\dps\int\varphi\circ f_k(x)dx = \dps\int\varphi\circ f(x)dx$,
e vale $\dps (*)$.

\bigskip
Suponhamos que existe $a_0 \in E$ tal que $\varphi(a_0) = \dps\int\varphi\circ fdx,
 \,\, \forall\varphi\in E^*$. Deste modo, temos que $\varphi\left(a_0 - \dps\int f(x)dx\right)
= 0, \,\, \forall\varphi\in E^*. $ Logo, pelo Teorema de Hahn-Banach, (ver tamb\'em \cite[corol\'ario E.8]{13}), existe $\psi \in E^*$ tal que $\psi\left(
a_0 - \dps\int f(x)dx\right) = \left\|a_0 - \dps\int f(x)dx\right\|,$ ou seja, $\left\|a_0 -
\dps\int f(x)dx\right\| = 0$ e portanto $a_0 = \dps\int f(x)dx.$ $\ai$\end{remark}

\bigskip Tendo este resultado, podemos ver que v\'arias propriedades que s\~ao
v\'alidas para integra\c c\~ao de fun\c c\~oes complexas s\~ao  v\'alidas
tamb\'em para fun\c c\~oes que tomam valores num espa\c co de Banach separ\'avel.

\begin{pretheoF}\label{A.12} {\rm  Se $f \in L^1(\R^{2n},E),$
ent\~ao $y \mapsto f(x,y)$ pertence a $L^1(\R^n,E)$ para quase todo $x\in\R^n$  e 
$x \mapsto f(x,y)$ pertence a $L^1(\R^n,E)$ para quase todo $y\in\R^n$.
Al\'em disso, as fun\c c\~oes $h(x) = \dps\int f(x,y)dy$ e $k(y) = \dps\int 
f(x,y)dx$ definidas q.s. pertencem a $L^1(\R^n,E)$ e vale:

\[\dps\int f(x,y)dxdy = \dps\int\left[ \dps\int f(x,y)dy\right]dx = \dps\int\left[\dps\int f(x,y)dx\right]dy.\]}\end{pretheoF}

\begin{remark} Por defini\c c\~ao,  temos que $\|f\| \in L^1(\R^{2n})$.
Logo, pelo Teorema de Fubini, \cite[(2.37)]{13}, temos que,
para quase todo $x\in\R^n$, $y\mapsto\|f(x,y)\|$ pertence a $L^1(\R^n)$,
isto \'e, $y\mapsto f(x,y) \in L^1(\R^n,E)$ q.s.. Analogamente para
$x\mapsto f(x,y)$. 

\medskip Tamb\'em temos que, se $\overline{h}(x) = \dps\int\|f(x,y)\|dy$,
$\overline{h} \in L^1(\R^n)$. E, como $\|h(x)\|\leq \overline{h}(x)$, $h \in L^1(\R^n,E)$.
Do mesmo modo, $k\in L^1(\R^n,E)$.

\medskip Por outro lado, dado qualquer $\varphi\in E^*$, temos 
$\varphi\circ f\in L^1(\R^{2n})$, pois $\dps\int\|\varphi\circ f(x,y)\|dxdy
\leq \|\varphi\|\dps\int\|f(x,y)\|dxdy < \infty$. Tamb\'em $\varphi\circ k, \,\, \varphi\circ h \in L^1(\R^n)$. Ent\~ao, apliquemos \cite[2.37]{3}:

\[\dps\int \varphi\circ h(x)dx = \dps\int\left[\varphi\left(\dps\int f(x,y)dy\right)\right]dx =
\dps\int\left[\dps\int\varphi\circ f(x,y)dy\right]dx = \] \[ =\dps\int\left[\dps\int\varphi\circ f(x,y)dx\right]dy = \dps\int \varphi\circ k(y)dy.\]             

Logo, pela proposi\c c\~ao \ref{A.11}, $\varphi\left(\dps\int h(x)dx\right) = \varphi\left(
\dps\int k(y)dy\right), \,\, \forall\varphi\in E^*$, 
e ent\~ao $\dps\int h(x)dx = \dps\int k(y)dy$, ou seja, 

\[\dps\int\left[\dps\int f(x,y)dy\right]dx = \dps\int\left[\dps\int f(x,y)dx\right]dy.\]     

Do mesmo modo, temos

\[\dps\int \varphi\circ f(x,y)dxdy = \dps\int \varphi\circ k(y)dy, \,\, \forall\varphi\in E^*
\Rightarrow  \] \[\Rightarrow \dps\int f(x,y)dxdy = \dps\int k(y)dy = \dps\int\left[\dps\int f(x,y)dx\right]dy.  \ai  \]\end{remark}

\vspace{1cm}
\begin{pretheoFC}\label{A.13} {\rm Seja
$f:[a,b] \rightarrow E$ uma fun\c c\~ao cont\'\i nua e diferenci\'avel tal que $f'$ (a fun\c c\~ao derivada de $f$) \'e limitada. Ent\~ao:

\[ \dps\int_{[a,b]} f'(t)dt = f(b) - f(a) .\]}

\end{pretheoFC}
\begin{remark} Dado $\varphi\in E^*$, qualquer, consideremos a fun\c c\~ao
$\varphi\circ f:[a,b] \rightarrow \C$. Seja $F(t) = \varphi\circ f(t)$. Notemos que

\[F'(t_0) = \dps\lim_{h\rightarrow 0}\frac{F(t_0+h)-F(t_0)}{h} =
\dps\lim_{h\rightarrow 0}\varphi\left(\frac{f(t_0+h)-f(t_0)}{h}\right) =
\varphi\circ f'(t_0).\]    

Logo, $F$ \'e cont\'\i nua, diferenci\'avel, com derivada limitada.  Al\'em disso,
como $F'$ \'e mensur\'avel, isto \'e, $\varphi\circ f'$ \'e mensur\'avel para todo
$\varphi\in E^*$, pelo teorema \ref{A.10} temos que $f'$ \'e mensur\'avel. Assim, por \cite[(3.36)]{3},
temos: $\dps\int^b_a F'(t)dt = F(b) - F(a)$, isto \'e,
$\varphi\left( \dps\int_{[a,b]}f'(t)dt\right) = \varphi(f(b)-f(a)), \,\, \forall\varphi\in E^*$
e, portanto, $\dps\int_{[a,b]}f'(t)dt = f(b) - f(a).$            $\ai$\end{remark}

\begin{precorol}\label{A.14}{\rm\bf Integra\c c\~ao por partes.} {\rm  Dadas $f, g \in C^\infty_c(\R^n,E)$, ent\~ao, para $-\infty<a<b<\infty$, vale que
 
\[ \dps\int_{[a,b]}f(x)g'(x)dx + \dps\int_{[a,b]}f'(x)g(x)dx = f(b)g(b) - f(a)g(a). \]                  }\end{precorol}

\begin{remark}  Notemos que $fg \in C^\infty_c(\R^n,E)$. Al\'em disso,
$(fg)'(x) = f(x)g'(x) + f'(x)g(x)$. Logo, pelo Teorema Fundamental do C\'alculo \ref{A.13}
aplicado \`a fun\c c\~ao $fg$,
\[\dps\int_{[a,b]}f(x)g'(x)dx + \dps\int_{[a,b]}f'(x)g(x)dx = f(b)g(b) - f(a)g(a).  \hspace{0,5cm}  \ai \] \end{remark}

\begin{precorol}\label{A14ummeio}{\rm   Seja $A$ uma $C^*$-\'algebra separ\'avel.
Dadas $f\in C^\infty(\R^n,A)$, com $f(x) = x^\gamma e^{ix\xi}, \gamma\in\N^n $ dado, e $g\in\sa$, ent\~ao temos que 
\[ \dps\int \left[\partial_{x_j}f(x)\right]g(x)dx = - \dps\int f(x)\left[\partial_{x_j}g(x)\right]dx.\]                          }\end{precorol}

\begin{remark}    Suponhamos, sem perda da generalidade, que $j=1$. Como $fg\in  L^1(\R^n,A)$, temos, pelo Teorema de Fubini \ref{A.12}, que
\[ \dps\int \partial_{x_1}\left[ f(x)g(x)\right] dx = \dps\int_{\R^{n-1}}\left( \dps\int_{\R}
\partial_{x_1}\left[ f(x)g(x)\right] dx_1\right) dx_2\cdots dx_n = \]
\[ = \dps\int_{\R^{n-1}}\left( \dps\lim_{R\rightarrow\infty}\dps\int_{[-R,R]}\partial x_1\left[ f(x)g(x)\right] dx_1\right)dx_2\cdots dx_n = \] 

\leftskip=-1cm 

\[ = \dps\int_{\R^{n-1}}\left(
\dps\lim_{R\rightarrow\infty}\left[f(R,x_2,\cdots ,x_n)g(R,x_2,\cdots,x_n) -  f(R,x_2,\cdots,x_n)g(R,x_2,\cdots,x_n)\right]\right)dx_2\cdots dx_n = 0\] 

\leftskip=0cm \noindent (ver Observa\c c\~ao \ref{A5ummeio}), pelo corol\'ario \ref{A.14}.    Portanto, \[\dps\int\left[\partial_{x_1}f(x)\right]
g(x)dx = -\dps\int f(x)\left[\partial_{x_1}g(x)\right]dx. \hspace{1cm} \ai\]\end{remark}

\begin{preprop}\label{A.15} {\rm  Seja $A$ uma C$^*$-\'algebra separ\'avel e seja $f \in L^2(\R^n,A)$. Ent\~ao, se $\dps\int f(x)^*f(x)dx = 0$, $f(x) = 0$, para quase todo $x\in\R^n$.                         }\end{preprop}

\begin{remark}  Dado um funcional positivo, $\psi\in A^*$, temos que $\psi\left(f(x)^*f(x)\right)\geq 0$. Por outro lado, pela proposi\c c\~ao \citeprop{A.11},
$\dps\int \psi\left(f(x)^*f(x)\right)dx = 0$, o que implica $\psi\left(f(x)^*f(x)\right) = 0$
q.s.. Como qualquer funcional linear $\varphi \in A^*$ \'e dado por
$\varphi = (\psi^+_1-\psi^-_1)+i(\psi^+_2-\psi^-_2)$, onde $\psi^+_i,\psi^-_i$ s\~ao
funcionais lineares positivos, $i=1,2$ (\cite[teorema 3.3.10 e coment\'arios acima]{5}), temos que 
$\varphi\left(f^*(x)f(x)\right) = 0,$ q.s., $\forall\varphi\in A^*.$

Notemos que a bola unit\'aria, $\mathcal{B}$, de $A^*$ \'e compacta na
topologia fraca$^*$ pelo teorema de Banach-Alaoglu (\cite[3.15]{6}) e
portanto, como $A$ \'e separ\'avel, $\mathcal{B}$ \'e metriz\'avel na
topologia fraca$^*$ (\cite[3.16]{6}). Logo a bola unit\'aria de $A^*$ \'e
separ\'avel (todo espa\c co m\'etrico compacto \'e separ\'avel). Seja
$\mathcal{D}\subseteq \mathcal{B}$ um subconjunto enumer\'avel denso na
topologia fraca$^*$.

Vimos que para qualquer $\varphi_m \in \mathcal{D}$,
$\varphi_m(f^*(x)f(x))=0$ q.s.. Portanto, o conjunto $E=\{x \in \R^n \, : \,
  \varphi_m(f^*(x)f(x))\neq 0, \textrm{ para algum } \varphi_m \in \mathcal{D}\}$
tem medida nula. Al\'em disso, se $x \notin E$ temos que
$\varphi_m(f^*(x)f(x))=0$ para todo $\varphi_m \in \mathcal{D}$ e
ent\~ao $\phi(f^*(x)f(x))=0$ para todo $\phi \in \mathcal{B}$. Assim, se $x
\notin E$, $\varphi(f^*(x)f(x))=0, \, \forall \, \varphi \in A^*$. Logo, se $x
\notin E$, $f(x)=0$; ou seja, $f(x)=0$ para quase todo $x \in \R^n$. $\ai$\end{remark}



\begin{preprop}\label{A.16}{\bf Mudan\c ca de vari\'avel.}
{\rm Sejam $T:\R^n\rightarrow \R^n$ um difeomorfismo $C^1$  e $D_x T$ a 
transforma\c c\~ao linear dada pela matriz $\left(\frac{\partial T_i}{\partial x_j}(x)\right)_{ij}$;
e seja $\dps\det D_x T$ o determinante da matriz acima.
Se $f\in L^1(\R^n,A)$, temos que 
\[\dps\int f(x)dx = \dps\int f\circ T(x)|\dps\det D_xT|dx.\]        }\end{preprop}

\begin{remark}  Dado qualquer $\varphi\in A^*, \,\, \varphi\circ f\in L^1(\R^n)$. Logo,
por \cite[(2.47)]{3}, $\dps\int \varphi\circ f(x)dx = \dps\int\varphi\circ f\circ 
T|\dps\det D_x T|dx$, o que implica que, $\forall\varphi\in A^*$, $\varphi\left(\dps\int
f(x)dx \right) = \varphi\left(\dps\int f\circ T|\dps\det D_xT|dx\right)$  e portanto
$\dps\int f(x)dx = \dps\int f\circ T|\dps\det D_xT|dx.$          $\ai$\end{remark}

\vspace{1,5cm} 
\begin{preprop}\label{A.22}{\rm  Dada uma $C^*$-\'algebra separ\'avel, $A$, temos que $\sa$
\'e denso em $L^2(\R^n,A)$. }\end{preprop}
\begin{remark} Vemos em \cite[(8.17)]{3} a demonstra\c c\~ao de que $\s$ 
\'e denso em $L^p(\R^n)$. A demonstra\c c\~ao de que $\sa$ \'e denso em $L^2(\R^n,A)$ 
\'e feita seguindo os mesmos passos. Primeiro, provamos que as fun\c 
c\~oes cont\'\i nuas de suporte compacto s\~ao densas em $L^2(\R^n,A)$. Depois,
a demonstra\c c\~ao de que as fun\c c\~oes $C^\infty$ de suporte compacto
s\~ao densas em $L^2(\R^n,A)$ segue por convolu\c c\~ao como no caso escalar.$\ai$\end{remark}





\newpage

\chapter{}

Usaremos aqui v\'arios resultados do ap\^endice A e, portanto, consideraremos 
$A$ uma $C^*$-\'algebra separ\'avel.

\vspace{1cm}
\begin{predef}\label{B.1} {\rm Dada $u \in \sa$, definimos a {\em Transformada de Fourier,
 $F$,  de $u$}, como segue:
\begin{center} $F(u)(x) = \dps\int e^{-ix\xi}u(\xi)\dc\xi,$ onde $\dc\xi = \dps\frac{1}{(2\pi)^{\frac{n}{2}}}
d\xi$. \end{center} }\end{predef}

Temos, como no caso escalar (\cite[teorema V.1.1]{12}), que $F$ \'e um operador cont\'\i nuo
de $\sa$ em $\sa$.

\begin{preprop}\label{B.2} {\rm $F$ \'e continuamente invers\'\i vel de $\sa$ em $\sa$ e
vale que
\[F^{-1}(u)(x) = \dps\int e^{ix\xi}u(\xi)\dc\xi, \hspace{0,5cm} u\in\sa. \hspace{0,5cm} (*)\]
}\end{preprop}

\begin{remark} Usaremos a fun\c c\~ao auxiliar $\phi(x) = e^{-\frac{\,|x|^2}{2}}$.
Notemos que $\hat{\phi}(\xi) = \phi(\xi)$ (\cite[exemplo V.1.2]{12}), onde
$\hat{\phi}$ \'e a transformada de Fourier (usual) de $\phi$. Tomemos, em
$\dps (*),  u=F(v)$, para $v\in\sa$.
\[v(x) = \dps\int e^{ix\xi}\left[\dps\int e^{-iy\xi}u(y)\dc y\right]\dc\xi.\]   

Aqui n\~ao podemos trocar a ordem de integra\c c\~ao pois $e^{i(x-y)\xi}v(y)$
n\~ao \'e integr\'avel em $\xi$ se $v(y)\not =0$. Por este motivo,
vamos usar a fun\c c\~ao $\phi$ do seguinte modo:

Seja $\phi_\epsi(x) = \phi(\epsi x)$, para $\epsi > 0$, dado.
\'E f\'acil ver que $\hat{\phi}_\epsi(x) = \epsi^{-n}\phi(\frac{x}{\epsi})$.
Assim, aplicando o Teorema de Fubini (\ref{A.12}), 
\[ \dps\int\phi_\epsi(\xi)e^{ix\xi}\left[\dps\int e^{-iy\xi}v(y)\dc y\right]\dc\xi =
\dps\int e^{-iy\xi}v(y)\left[\dps\int e^{ix\xi}\phi_\epsi(\xi)\dc\xi\right]\dc y  = \]
\[ = \dps\int v(y)\left[\dps\int e^{-i(y-x)\xi}\phi_\epsi(\xi)\dc\xi\right]\dc y =
\dps\int v(y)\hat{\phi_\epsi}(y-x)\dc y =\]
\[ \dps\int v(y+x)\hat{\phi_\epsi}(y)\dc y = 
\dps\int v(y+x)\epsi^{-n}\phi\left(\frac{y}{\epsi}\right)\dc y
= \dps\int v(x+\epsi y)\phi(y)\dc y .\]  E temos $\dps\lim_{\epsi\rightarrow 0}
\phi_\epsi(x) = 1$, $\dps\lim_{\epsi\rightarrow 0}v(x+\epsi y) = v(x)$ e,
pelo visto acima:
\[ \dps\int e^{ix\xi}\phi_\epsi(\xi)F(v)(\xi)\dc\xi = \dps\int v(x+\epsi y)\phi(y)\dc y.\]   

Por outro lado, $\|e^{ix\xi}\phi_\epsi(\xi)F(v)(\xi)\| \leq \|F(v)(\xi)\| \in L^1(\R^n)$ e            
$\|v(x+\epsi y)\phi(y)\| \leq \dps\sup_{x\in\R^n}\|v(x)\|\|\phi(y)\| \in L^1(\R^n)$.

Tomando $\epsi \rightarrow 0$, pelo Teorema da Converg\^encia Dominada (\ref{A.7}),
temos 
$$\dps\int e^{ix\xi}F(v)(\xi)\dc\xi = \dps\int v(x)\phi(y)\dc y .$$
 
Como $\dps\int\phi(y)\dc y = 1$, $v(x) = \dps\int e^{ix\xi}F(v)(\xi)\dc\xi$, 
o que prova $\dps (*)$.

A continuidade de $F^{-1}$ se demonstra de modo an\'alogo \`a continuidade
de $F$, (\cite[teorema V.1.1]{12}).
$\ai$\end{remark}

\begin{preprop}\label{B.3} {\rm  $F$ \'e um operador "unit\'ario" em $\sa$, 
com respeito ao produto interno definido no in\'\i cio do cap\'\i tulo 1.                                   }\end{preprop}

\begin{remark} Dadas $u, v \in \sa$, pelo Teorema de Fubini,
\[<F(u),v> = \dps\int \left[\dps\int e^{-ix\xi}u(\xi)\dc\xi\right]^*v(x)dx =
\dps\int u(\xi)^*\left[\dps\int e^{ix\xi}v(x)\dc x\right]d\xi = <u,F^{-1}(v)>.\] Portanto
$F^* = F^{-1}$ e $F$ \'e unit\'ario em $\sa$.           $\ai$\end{remark}

\begin{precorol}\label{B.4}{\rm Se $u \in \sa$, ent\~ao $\|F(u)\|_2 = \|u\|_2$. }\end{precorol}

\begin{remark}  $<F(u),F(u)> = <u,F^{-1}(F(u))> = <u,u>$, portanto
$\|F(u)\|_2=\|u\|_2$.    $\ai$\end{remark}

\begin{precorol}\label{B.5}{\rm  Dadas $f, g \in L^2(\R^n,A)$, temos que $<F(f),F(g)> = <f,g>$.                                     }\end{precorol}

\begin{remark}    Sabemos, pela proposi\c c\~ao \citeprop{3.3} que $L^2(\R^n,A) 
\subseteq \SA$. Assim, existem seq\"u\^encias $(f_k)_{k\in\N}$ e $(g_l)_{l\in\N}$
de fun\c c\~oes em $\sa$ tais que $\dps\lim_{k\rightarrow\infty}\|f_k-f\|_2 = 0$
e $\dps\lim_{l\rightarrow\infty}\|g_l-g\|_2 = 0$. Pelo corol\'ario \ref{B.4}, $F$ \'e
cont\'\i nua em $\SA$ com a norma $\|\cdot\|_2$. Assim, pela
proposi\c c\~ao \ref{B.3}, $<F(f),F(g)> = \dps\lim_{k,l}<F(f_k),F(g_l)> = \dps\lim_{k,l}
<f_k,g_l> = <f,g>.$                      $\ai$\end{remark}

\end{document}